\documentclass[11pt,english]{article}
\usepackage{microtype}
\usepackage{bbm}
\usepackage{amsthm}
\usepackage{amsmath}
\usepackage{amssymb}
\usepackage{euscript}
\usepackage{subfigure}
\usepackage{enumitem}
\usepackage{graphicx,caption}
\usepackage{geometry}
\usepackage{tikz}
\usetikzlibrary{arrows}
\usepackage{color}
\usepackage[bookmarksnumbered=true]{hyperref} 
\hypersetup{
colorlinks = true,
linkcolor = black,
anchorcolor = blue,
citecolor = blue,
filecolor = blue,
urlcolor = blue
}

\newtheorem{theorem}{Theorem}
\newtheorem{proposition}[theorem]{Proposition}
\newtheorem{corollary}[theorem]{Corollary}
\newtheorem{lemma}{Lemma}
\theoremstyle{definition}
\newtheorem*{remark}{Remark}
\newtheorem{example}{Example}
\newtheorem{definition}{Definition}

\renewcommand{\epsilon}{\varepsilon}
\def\R{\mathbb{R}}
\def\N{\mathbb{N}}
\def\Z{\mathbb{Z}}
\DeclareMathOperator{\Span}{span}
\DeclareMathOperator{\interior}{int}

\def\cH{\EuScript{H}}
\def\cM{\EuScript{M}}

\renewcommand{\epsilon}{\varepsilon}
\def\N{\mathbb{N}}
\def\R{\mathbb{R}}
\def\Z{\mathbb{Z}}
\def\cA{\EuScript{A}}
\def\cB{\EuScript{B}}
\def\cC{\EuScript{C}}
\def\cF{\EuScript{F}}
\def\cG{\EuScript{G}}

\def\cM{\EuScript{M}}

\def\cQ{\EuScript{Q}}

\def\cE{\EuScript{E}}
\def\cP{\EuScript{P}}

\def\Id{\text{\rm Id}}

\begin{document}

\title{\textbf{Nonadditive families of potentials: physical equivalence and some regularity relations}}

\author{Carllos Eduardo Holanda}
\date{%
\small{Institute of Mathematics, Federal University of Alagoas, Macei\'o, 57072-970, Alagoas, Brazil\\
\textit{E-mail address:} \texttt{c.eduarddo@gmail.com}
}}

\maketitle

\begin{abstract}
We show that additive and asymptotically additive families of continuous functions with respect to suspension flows are physically equivalent. In particular, the equivalence result holds for hyperbolic flows and some classes of expansive flows in general. Moreover, we show how this equivalence result can be used to extend the nonadditive thermodynamic formalism and multifractal analysis. In the second part of this work, we obtain a Liv\v{s}ic-like result for nonadditive families of potentials and also address the H\"older and Bowen regularity problem for the physical equivalence relations with respect to hyperbolic symbolic flows. 
\end{abstract}

\tableofcontents

\section{Introduction}

This work is mainly a contribution to the study of thermodynamic formalism and multifractal analysis for flows, taking into consideration the relationships between the additive classical world of potentials and the nonadditive world of families of potentials, which started in \cite{Cun20}. The first part is concerned with the physical equivalence problem for flows and how it affects some results in thermodynamics and multifractal analysis. The second part is mostly inspired by the recent developments in \cite{HS24}, and it deals with regularity issues arising naturally from the problem of how information can be exchanged between additive and nonadditive families through the physical equivalence relation. 

\subsection{Physical equivalence}

Let $X$ be a topological space and $T: X \to X$ a map. A sequence of functions $(f_n)_{n \ge 1}$ is \emph{asymptotically additive} with respect to $T$ if for each $\epsilon > 0$ there exists a function $f: X \to \R$ such that 
\[
\limsup_{n \to \infty}\frac{1}{n}\|f_n - S_nf\|_{\infty} < \epsilon,
\]
where $S_nf := \sum_{k=0}^{n-1}f \circ T^{k}$ and $\|\cdot\|_{\infty}$ is the supremum norm. Notice that the sequence $(S_nf)_{n \ge 1}$ is additive with respect to $T$, that is,
\[
S_{m+n}f(x) = S_mf(x) + S_nf(T^{m}(x)) \quad \textrm{for all $x \in X$ and $m, n \ge 1$}. 
\]

A sequence $\cF = (f_n)_{n \ge 1}$ is \emph{almost additive} with respect to $T$ if there exists $C > 0$ such that 
\[
-C + f_{m}(x) + f_{n}(T^m(x)) \le f_{m+n}(x) \le f_{m}(x) + f_{n}(T^m(x)) + C
\]
for every $x \in X$ and all $m,n \ge 1$. It was showed in \cite{FH10} that almost additive sequences are in fact asymptotically additive. 

Inspired by statistical mechanics, as in \cite{Cun20}, we say that two nonadditive sequences of functions $\cF:= (f_n)_{n \ge 1}$ and $\cG:= (g_n)_{n \ge 1}$ are \emph{physically equivalent}, or $\cF$ is \emph{physically equivalent} to $\cG$, if 
\[
\lim_{n \to \infty}\frac{1}{n}\|f_n - g_n\|_{\infty} = 0. 
\]
Surprisingly, N. Cuneo showed in \cite{Cun20} that asymptotically additive sequences are physically equivalent to additive sequences with respect to a continuous map. This result has a direct impact in the study of nonadditive thermodynamic formalism and multifractal analysis for discrete-time dynamical systems (see \cite{Cun20} and references within).

Motivated by the nonadditive thermodynamic formalism and multifractal analysis for flows, and  inspired by Cuneo's equivalence theorem, in this paper we investigate the same equivalence problem in the case of continuous flows.  

Let $\Phi = (\phi_t)_{t \in \R}$ be a continuous flow on a topological space $X$. A family $a=(a_t)_{t\ge0}$ of functions $a_t\colon X\to \R$ is said to be \textit{almost additive} with respect to $\Phi$ on $X$ if there exists a constant $C>0$ such that
\[
-C + a_{t} + a_{s}\circ\phi_t \le a_{t+s} \le a_{t} + a_{s}\circ\phi_t + C
\]
for every $t,s \ge 0$. Notice that, for each function $b: X \to \R$, the family $(S_tb)_{t \ge 0}$ given by $S_tb = \int_{0}^{t}(b \circ \phi_s)ds$ is additive with respect to $\Phi$. 
 
We also say that a family of functions $a = (a_t)_{t\ge 0}$ is \emph{asymptotically additive} with respect to $\Phi$ on $X$ if for each $\epsilon > 0$ there exists a function $b_{\epsilon}: X \to \R$ such that 
\[
\limsup_{t \to \infty}\frac{1}{t}\bigg\|a_t - \int_{0}^{t}(b_{\epsilon}\circ \phi_s) ds\bigg\|_{\infty} \leq \epsilon.
\]
Proceeding as in \cite{FH10}, one can see that every almost additive family of functions is asymptotically additive.

Following the definition for discrete-time dynamical systems, we say that two families of functions $a = (a_t)_{t\ge 0}$ and $b = (b_t)_{t \ge 0}$ are \emph{physically equivalent}, or $a$ is \emph{physically equivalent} to $b$, with respect to the flow $\Phi$ if 
\[
\lim_{t \to \infty}\frac{1}{t}\|a_t - b_t\|_{\infty} = 0.
\]
We observe that physically equivalent almost additive families have the same topological pressure and share the same equilibrium measures (see \cite{BH21a}). Moreover, they also share the same level sets and the same maximizing measures in the sense of ergodic optimization (see for example \cite{BD09}, \cite{BH21b}, \cite{BHVZ21}, \cite{HLMXZ19} and \cite{MSV20}).  

In the present work, we show that asymptotically additive families of continuous functions are physically equivalent to additive families of continuous functions with respect to suspension flows, and in particular, hyperbolic flows and expansive flows admitting a measure of full support. As in the discrete-time case, this physical equivalence result for flows has the potential to facilitate and simplify many extensions of the nonadditive thermodynamic formalism, multifractal analysis and even ergodic optimization for flows in general (see section \ref{APP}). 

We also address the full problem of equivalence with respect to continuous flows in general, and it turns out that this problem is intrinsically connected to the following dynamical embedding problem, which is also stated as an open question in \cite{BHVZ21}: 
\begin{itemize}
\item \emph{Given a continuous flow $\Phi$ on a metric space $X$ and a continuous function $\widetilde{b}: X \to \R$, is there a continuous function $b:X \to \R$ such that}
\[
\widetilde{b} = \int_{0}^{1}(b \circ \phi_s)ds \quad ? 
\]
\end{itemize}
We stress here that a positive answer for this embedding problem implies a positive answer for the general physical equivalence problem. 

In this work, we are also able to describe a sufficient condition on the function and on the flow where the dynamical embedding problem can be answered affirmatively. On the other hand, we also give a simple counter-example, showing that the embedding problem cannot be solved in full generality (see section \ref{GER}). Moreover, concerning the general physical equivalence problem, we give a positive answer in the case of continuous flows with uniquely ergodic time-one maps. However, as far as we know, the equivalence problem in full generality remains open. 

\subsection{A nonadditive Liv\v{s}ic theorem and some regularity relations}

Based on the physical equivalence relationship between asymptotically additive, almost additive and additive families obtained in the first part of this work (see Theorem \ref{main}), one can naturally consider the problem of identifying the meaningful different levels of regularity that the physical equivalence relation can maintain. In our setup, the most relevant types of regularity are the ones involving Bowen and H\"older functions together with families having the bounded variation condition (see sections \ref{chm} and \ref{mr} for the definitions). In the context of hyperbolic suspension flows and related hyperbolic setups, the space of H\"older continuous functions is contained in the space of Bowen continuous functions. Furthermore, by definition, an additive family generated by a Bowen function has bounded variation with respect to any flow in general.

More specifically, inspired by the recent work \cite{HS24} for maps, we can consider three types of regularity problems: 
\begin{itemize}
\item \underline{Bowen regularity.} \emph{Given any almost additive family of continuous functions $a = (a_t)_{t \ge 0}$ with bounded variation with respect to an hyperbolic suspension flow, is there any Bowen continuous function $b$ such that $(S_tb)_{t \ge 0}$ is physically equivalent to $a$ ? }

\item \underline{Uniform bound.} \emph{Given any almost additive family of continuous functions $a = (a_t)_{t \ge 0}$ with bounded variation with respect to an hyperbolic suspension flow, is there any continuous function $b$ such that $\sup_{t \ge 0}\|a_t - S_tb\|_{\infty}$ ?}

\item \underline{H\"older regularity.} \emph{Given any almost additive family of H\"older continuous functions $a = (a_t)_{t \ge 0}$ with bounded variation with respect to an hyperbolic suspension flow, is there any H\"older continuous function $b$ such that $(S_tb)_{t \ge 0}$ is physically equivalent to $a$~? }
\end{itemize}

Notice that the uniform bound immediately implies Bowen regularity and, a positive answer to the H\"older regularity question in this context also gives an affirmative answer to the Bowen regularity one. It is important mentioning that these regularity issues are also interesting in the more general case of asymptotically additive families.

Apart from the intrinsic interest involved, there are important consequences regarding the regularities problems posted above. In fact, an affirmative answer to the Bowen regularity question would imply that one can obtain the uniqueness of equilibrium measures for almost additive families with bounded variation directly from the (additive) classical uniqueness result for a single potential (see \cite{Fra77}, \cite{BH21a}). On the other hand, a positive answer to the H\"older regularity equivalence problem would either simplify or automatically extend some relevant results in ergodic optimization for almost and asymptotically additive families (see for example \cite{BHVZ21}, \cite{HLMXZ19} and \cite{MSV20}). In addition to that, the H\"older regularity problem is also connected to the existence of a general important class of asymptotically and almost additive families with analytic nonadditive topological pressure function, which also could allow us to obtain higher-regularity of the entropy and dimension spectra for these classes of families of potentials with respect to dynamical systems admitting some hyperbolic behavior (see for example \cite{Rue78}, \cite{BD04}, \cite{BS00}, \cite{PS01}, \cite{BH21b} and \cite{BH22c}).    
	
In the second part of this work we are able to obtain a characterization result for almost additive families, which also has an interesting and significant connection with the regularity problems mentioned above:

(Theorem \ref{BOU}). \emph{Let $\Phi = (\phi_t)_{t \ge 0}$ be a topologically transitive continuous flow on a compact metric space $X$ and satisfying the Closing Lemma. Let $\cB = (b_t)_{t \ge 0}$ be an almost additive family of continuous functions (with respect to $\Phi$) with bounded variation. Then, the following are equivalent:}
\begin{enumerate}
\item $\lim_{t \to \infty}\|b_{t}\|_{\infty}/t = 0$;
	
\item $\sup_{t \ge 0}\|b_t\|_{\infty} < \infty$;
	
\item there exists $K > 0$ such that  $|b_t(p)| \le K$ for all $p \in X$ and $t \ge 0$ with $\phi_t(p) = p$. 
\end{enumerate}
This result is based on the discrete-time counterpart recently showed in \cite{HS24}. It is important mentioning that this theorem is proved here directly in the realm of flows, and without using any of the physical equivalence results for discrete and continuous-time dynamical systems (respectively Theorem 1.2 in \cite{Cun20} and Theorem \ref{main}). Also observe that Theorem \ref{BOU} immediately gives a setup where the uniform bound and Bowen regularity problems are actually equivalent. 

This characterization result holds, in particular, for some types of suspension flows and hyperbolic symbolic flows in general and shows a deeper layer of the physical equivalence relation for these types of nonadditive families (see Corollary \ref{EQ2}). Furthermore, we apply Theorem \ref{BOU} to the study of some matrix cocycles over flows and also to classify the equilibrium states for almost additive families based on their cohomology relations, strongly hinting that, in fact, Theorem \ref{BOU} also works as a nonadditive version of the classical Liv\v{s}ic Theorem for flows (see \cite{Liv71} and \cite{Liv72}). In addition, we also show that Theorem \ref{BOU} has an optimal setup in the sense that it is no longer valid for asymptotically and subadditive families in general. 

Building on some examples showed in \cite{HS24} for full shifts of finite type, we demonstrate how to construct almost and asymptotically additive families of H\"older continuous potentials satisfying the bounded variation condition with respect to hyperbolic symbolic flows and which are not physically equivalent to any additive family generated by a H\"older (Bowen) continuous potential.  As in the discrete-time case, these examples indicate that almost and asymptotically additive families with bounded variation do not always have the same expected optimal regularity properties of H\"older continuous functions in hyperbolic and related scenarios. 

\subsection{Organization of the paper}

We start proving our main equivalence result for suspension flows and, consequently, for hyperbolic flows. After that, we proceed to show how to extend the physical equivalence theorem to expansive flows using the very recent work \cite{GS22}. In the following, we also give natural sources of almost and asymptotically additive families of continuous functions, and also treat the general physical equivalence problem for continuous flows. In the final section of the first part, we conclude with some applications and consequences of the equivalence theorem, extending parts of the nonadditive thermodynamic formalism and multifractal analysis for continuous-time dynamical systems.  

In the second part  we start studying some different notions of cohomology for almost and asymptotically additive families in general. In the next section, we state and prove our nonadditive version of the Liv\v{s}ic theorem for flows and show how it can be applied to the context of matrix cocycles over flows. We proceed to study and compare some related notions of nonadditive Gibbs and weak Gibbs measures with respect to flows and, as another application of Theorem \ref{BOU}, we demonstrate how to classify almost additive families based on cohomology relations and equilibrium measures. In the next section, we start with a simple non-hyperbolic example where the uniform bound problem can always be positively answered, and another example where the equivalences in Theorem \ref{BOU} do not hold. In the next subsection, based on \cite{HS24}, we show how to build examples of almost additive families of H\"older continuous potentials satisfying the bounded variation property but never physically equivalent to any additive family generated by a H\"older continuous potential, giving a negative answer to the H\"older regularity equivalence problem. Proceeding to the next subsection, we show a way of categorizing almost additive families based on the different types of physical equivalence relations with the additive setup, and we conclude our work showing a construction of asymptotically additive families satisfying the bounded variation condition with respect to a hyperbolic symbolic flow, with a unique equilibrium measure but not physically equivalent to any additive family with bounded variation, giving a negative answer to the Bowen regularity problem for the asymptotically additive case. 

\section{Part I: physical equivalence}

In this part, we are going to show that additive and asymptotically additive families of continuous functions are physically equivalent with respect to suspension flows and also with respect to some expansive flows. Moreover, we are going to give some natural examples of asymptotically additive families and also explore the general problem of physical equivalence for continuous flows. In the end, we will show some consequences of the physical equivalence result in thermodynamic formalism and multifractal analysis for flows. 

\subsection{Suspension flows}

Let $X$ be a compact metric space, $T\colon X \to X$ a homeomorphism and $\tau\colon X \to (0, \infty)$ a positive continuous function. Consider the space
\[
W = \bigl\{(x,s) \in X \times \R: 0 \le s \le \tau(x)\bigr\}
\]
and let $Y$ be the set obtained from $W$ identifying $(x,\tau(x))$ with $(T(x),0)$ for each $x \in X$. Then a certain distance introduced by Bowen and Walters in~\cite{BW72} makes $Y$ a compact metric space.
The \emph{suspension flow} over $T$ with height function $\tau$ is the flow $\Phi = (\phi_t)_{t \in \R}$ on $Y$ with the maps $\phi_t\colon Y \to Y$ defined by $\phi_t(x,s) = (x, s+t)$.

Let $\mu$ be a $T$-invariant probability measure on $X$. One can show that $\mu$ induces a $\Phi$-invariant probability measure $\nu$ on $Y$ such that
\[
\int_{\Lambda}g\,d\nu = \frac{\int_{X}I_g\,d\mu}{\int_{X}\tau \,d\mu}
\]
for any continuous function $g\colon Y \to \R$, where $I_g(x) = \int_0^{\tau(x)}(g\circ\phi_s)(x)\,ds$.
Conversely, any $\Phi$-invariant probability measure $\nu$ on $Y$ is of this form for some $T$-invariant probability measure $\mu$ on~$X$. Abramov's entropy formula says that
\[
h_{\nu}(\Phi) = \frac{h_{\mu}(T)}{\int_{X}\tau \,d\mu}.
\]

The next result establishes the physical equivalence between asymptotically additive and additive families of functions with respect to  suspension flows. 

\begin{theorem}\label{main}
Let $\Phi = (\phi_t)_{t \in \R}$ be a suspension flow over a continuous map $T: X \to X$ and $a= (a_t)_{t\ge 0}$ be an asymptotically additive family of continuous functions with respect to $\Phi$. Then, there exists a continuous function $b: Y \to \R$ such that
\[
\lim_{t \to \infty}\frac{1}{t}\bigg\|a_t - \int_{0}^{t}(b \circ \phi_s)ds\bigg\|_{\infty} = 0.
\]
\end{theorem}

Consider the sequence of continuous functions $c= (c_n)_{n \ge 1}$ on $X$ given by $c_n(x) = a_{\tau_n(x)}(x)$,  where 
\[
\tau_n(x) = \sum_{k=0}^{n-1}\tau(T^{k}(x)) \quad \textrm{for every $x \in X$}.
\]

\begin{lemma} \label{LE1}
There exists a continuous function $\xi: X \to \R$ such that
\begin{equation}\label{DCT}
\lim_{n \to \infty}\frac{1}{n}\sup_{x \in X}\bigg|c_n(x) - \sum_{k=0}^{n-1}\xi(T^{k}(x))\bigg| = 0.
\end{equation}
\end{lemma}
\begin{proof}
Since $a$ is asymptotically additive with respect to the flow $\Phi$, for each $\epsilon>0$ there exists a continuous function $b_{\epsilon}: Y \to \R$ such that
\begin{equation}\label{AAF}
\limsup_{t \to \infty}\frac{1}{t}\bigg\|a_t - \int_{0}^{t}(b_{\epsilon}\circ \phi_s)ds\bigg\|_{\infty} \le \epsilon.
\end{equation}

It follows from the proof of Lemma 15 in \cite{BH21b} that 
\[
\int_{0}^{\tau_n(x)}(b_{\epsilon}\circ \phi_s)(x)ds = \sum_{k=0}^{n-1}(I_{b_{\epsilon}}\circ T^{k})(x) \quad \textrm{for every $x \in Y$ and $n \ge 1$}. 
\]
Notice that for each $t > 0$ there exists a unique $n \in \N$ such that $\tau_n(x)\le t \le \tau_{n+1}(x)$ with $t-\tau_n(x) \in [0, \sup\tau)$. Then, in particular from \eqref{AAF}, we have 
\begin{equation}\label{ETQ}
\begin{split}
&\limsup_{n \to \infty}\frac{1}{n}\sup_{x \in X}\bigg|c_n(x) - \sum_{k=0}^{n-1}(I_{b_{\epsilon}}\circ T^{k})(x)\bigg| \\
& = \limsup_{n \to \infty}\frac{1}{\tau_n(x)}\bigg(\frac{\tau_n(x)}{n}\bigg)\sup_{x \in X}\bigg|a_{\tau_n(x)}(x) - \int_{0}^{\tau_n(x)}(b_{\epsilon}\circ \phi_s)(x)ds\bigg| \le \epsilon \sup\tau
\end{split}
\end{equation}
for any $x \in X$. 
Since $\epsilon > 0$ is arbitrarily small, this implies that the sequence $c$ is asymptotically additive with respect to the map $T: X \to X$. Now we can apply Theorem 1.2 in \cite{Cun20} to guarantee the existence of a continuous function $\xi: X \to \R$ satisfying \eqref{DCT}, as desired. 
\end{proof}
\begin{lemma}\label{LE2}
There exists a continuous function $b: Y \to \R$ such that $I_{b}|_{X}= \xi$. 
\end{lemma}
\begin{proof}
Following \cite{BRW04}, we can just define $b: Y \to \R$ as 	
\[
b(\phi_s(x)) = \frac{\xi(x)}{\tau(x)}\psi'\bigg(\frac{s}{\tau(x)}\bigg)
\]
for each $x \in X$ and $s \in [0,\tau(x)]$, where $\psi: [0,1] \to [0,1]$ is any nondecreasing $C^{1}$ function such that $\psi(0) = 0$, $\psi(1) = 1$ and $\psi'(0) = \psi'(1) = 0$. 
\end{proof}

\begin{lemma}\label{LE3} We have
\[
\limsup_{t \to \infty}\frac{1}{t}\sup_{x \in X}|a_t(x) - a_{\tau_n(x)}(x)| = 0.   
\]
\end{lemma}
\begin{proof}
Since $a$ is asymptotically aditive, for each $\epsilon > 0$ there exists $b_{\epsilon}: Y \to \R$ such that
\[
\limsup_{t \to \infty}\frac{1}{t}\bigg\|a_t - \int_{0}^{t}(b_{\epsilon} \circ \phi_s)ds\bigg\|_{\infty}\le \epsilon.
\]
Moreover, following as in the proof of Lemma 15 in \cite{BH21b}, one can check that for each $t > 0$ there exists a unique $n \in \N$ with $t = \tau_n(x) + \kappa$ for some $\kappa \in [0, \sup\tau]$ such that 
\[
\bigg|\int_{0}^{t}(b_{\epsilon}\circ \phi_s)(x)ds - \sum_{k=0}^{n-1}(I_{b_{\epsilon}}\circ T^{k})(x)\bigg| \le \|b_\epsilon\|_{\infty}\sup\tau.
\]
Then, since 
\[
\begin{split}
\sup_{x \in X}|a_t(x) - a_{\tau_n(x)}(x)| &\le \sup_{x \in X}\bigg|a_t(x) - \int_{0}^{t}(b_{\epsilon}\circ \phi_s)(x)ds\bigg|\\
& +  \sup_{x \in X}\bigg|\int_{0}^{t}(b_{\epsilon}\circ \phi_s)(x)ds - \sum_{k=0}^{n-1}(I_{b_{\epsilon}}\circ T^{k})(x)\bigg|\\
& +  \sup_{x \in X}\bigg|c_n(x) - \sum_{k=0}^{n-1}(I_{b_{\epsilon}}\circ T^{k})(x)\bigg|
\end{split}	
\]
for every $t \ge 0$, it follows from \eqref{ETQ} that
\[
\limsup_{t \to \infty}\frac{1}{t}\sup_{x \in X}|a_t(x) - a_{\tau_n(x)}(x)| \le \epsilon (1+\sup\tau) + \limsup_{t \to \infty}\frac{1}{t}(\|b_{\epsilon}\|_{\infty}\sup\tau) = \epsilon (1+\sup\tau).
\]
Hence, letting $\epsilon \to 0$ we obtain
\[
\limsup_{t \to \infty}\frac{1}{t}\sup_{x \in X}|a_t(x) - a_{\tau_n(x)}(x)| = 0,
\]
as desired.
\end{proof}

%\begin{lemma}\label{LE4}
%Let $g = (g_t)_{t\ge0}$ be an asymptotically additive family of continuous functions with respect to the flow $\Phi$. Then:

%\begin{enumerate}
%\item the limit $\lim_{t \to \infty}\frac{1}{t}\sup_{x \in X}|g_t(x)|$ exists;
%\item for each $s \in \R$ and $x \in Y$ we have
%\[
%\limsup_{t \to \infty}\frac{g_t(x)}{t} = \limsup_{t \to \infty}\frac{g_t(\phi_s(x))}{t} \quad \textrm{and} \quad \liminf_{t \to \infty}\frac{g_t(x)}{t} = \liminf_{t \to \infty}\frac{g_t(\phi_s(x))}{t}.
%\]
%\end{enumerate}
%\end{lemma}
%\begin{proof}
	
%\end{proof}

\emph{Proof of Theorem \ref{main}.}
Using Lemma \ref{LE2}, let $b: Y \to \R$ be a function such that $I_{b}|_X = \xi$ and define the family of continuous functions $\Delta:= (\Delta_t)_{t \ge 0}$ as
\begin{equation}\label{DD1}
\Delta_{t}(x):= a_t(x) - \int_{0}^{t}(b\circ \phi_s)(x)ds \quad \textrm{for all $x \in Y$ and $t \ge 0$.}
\end{equation}

For each $x \in X$, we have
\[
\begin{split}
&\limsup_{t \to \infty}\frac{1}{t}\sup_{x \in X}|\Delta_t(x)| \le \limsup_{t \to \infty}\frac{1}{t}\sup_{x \in X}\bigg|a_{\tau_n(x)}(x) - \int_{0}^{\tau_n(x)}(b\circ \phi_s)(x)ds\bigg| \\
&+ \limsup_{t \to \infty}\frac{1}{t}\sup_{x \in X}|a_t(x) - a_{\tau_n(x)}(x)| \\
&+ \limsup_{t \to \infty}\frac{1}{t}\sup_{x \in X}\bigg|\int_{0}^{\tau_n(x)}(b\circ \phi_s)(x)ds - \int_{0}^{t}(b\circ \phi_s)(x)ds\bigg|\\
& \le \biggl(\frac{1}{\inf \tau}\biggr)\limsup_{n \to \infty}\frac{1}{n}\sup_{x \in X}\bigg|c_n(x) - \sum_{k=0}^{n-1}\xi(T^{k}(x))\bigg| \\
&+ \limsup_{t \to \infty}\frac{1}{t}\sup_{x \in X}(|a_t(x) - a_{\tau_n(x)}(x)| +  \|b\|_{\infty}\sup\tau).
\end{split}
\]
Hence, it follows from the lemmas \ref{LE1} and \ref{LE3} that 
\begin{equation}\label{AA0}
\limsup_{t \to \infty}\frac{1}{t}\sup_{x \in X}|\Delta_t(x)| = 0.
\end{equation}

We claim that 
\[
\limsup_{t \to \infty}\frac{1}{t}\sup_{s \in [0,\tau]}\sup_{x \in X}|\Delta_t(\phi_s(x))| = 0.
\]
By the asymptotic additivity of $a$, given $\epsilon > 0$ there exists $b_{\epsilon}: Y \to \R$ such that, in particular, 
\begin{equation}\label{AA1}
\limsup_{t \to \infty}\frac{1}{t}\sup_{x \in X}\bigg|a_t(x) - \int_{0}^{t}(b_{\epsilon}\circ \phi_u)(x)du\bigg| \le \epsilon \quad \textrm{and}
\end{equation}
\begin{equation}\label{AA2}
\limsup_{t \to \infty}\frac{1}{t-s}\sup_{x \in X}\bigg|a_{t-s}(\phi_s(x)) - \int_{0}^{t-s}(b_{\epsilon}\circ \phi_{u+s})(x)du\bigg|\le \epsilon \quad \textrm{for every $s \in [0,\sup\tau]$.}
\end{equation}
Since
\[
\begin{split}
&\sup_{x \in X}|a_{t-s}(\phi_s(x)) - a_t(x)| \le \sup_{x \in X}\bigg|a_{t-s}(\phi_s(x)) - \int_{0}^{t-s}(b_{\epsilon}\circ \phi_{u+s})(x)du\bigg|\\
& + \sup_{x \in X}\bigg|a_t(x) - \int_{0}^{t}(b_{\epsilon}\circ \phi_u)(x)du\bigg|+ \sup_{x \in X}\bigg|\int_{s}^{t}(b_{\epsilon}\circ \phi_u)(x)du - \int_{0}^{t}(b_{\epsilon}\circ \phi_u)(x)du\bigg|\\
&\le \sup_{x \in X}\bigg|a_{t-s}(\phi_s(x)) - \int_{0}^{t-s}(b_{\epsilon}\circ \phi_{u+s})(x)du\bigg| + \sup_{x \in X}\bigg|a_t(x) - \int_{0}^{t}(b_{\epsilon}\circ \phi_u)(x)du\bigg| \\
&+ \|b_{\epsilon}\|_{\infty}\sup\tau
\end{split}
\]
for every $t \ge s$ and $s \in [0,\sup\tau]$, it follows from \eqref{AA1} and \eqref{AA2} that
\[
\frac{1}{t}\sup_{s \in [0,\tau]}\sup_{x \in X}|a_{t-s}(\phi_s(x)) - a_t(x)| \le 2\epsilon + 2 \epsilon + \epsilon = 5\epsilon \quad \textrm{for $t$ sufficiently large.}
\]
Then, the arbitrariness of $\epsilon$ gives that
\begin{equation}\label{BB1}
\limsup_{t \to \infty}\frac{1}{t}\sup_{s \in [0,\tau]}\sup_{x \in X}|a_{t-s}(\phi_s(x)) - a_t(x)| = 0.
\end{equation}

From \eqref{DD1} one can check that
\[
\sup_{s \in [0,\tau]}\sup_{x \in X}|\Delta_{t-s}(\phi_s(x)) - \Delta_t(x)| \le \sup_{s \in [0,\tau]}\sup_{x \in X}|a_{t-s}(\phi_s(x)) - a_t(x)| + \|b\|_{\infty}\sup\tau
\]
for every $t \ge s$.
Hence, it follows from \eqref{BB1} that
\begin{equation}\label{AAB}
\limsup_{t \to \infty}\frac{1}{t}\sup_{s \in [0,\tau]}\sup_{x \in X}|\Delta_{t-s}(\phi_s(x)) - \Delta_t(x)| = 0.
\end{equation}
Now observing that 
\[
\sup_{s \in [0,\tau]}\sup_{x \in X}|\Delta_{t-s}(\phi_s(x))| \le \sup_{s \in [0,\tau]}\sup_{x \in X}|\Delta_{t-s}(\phi_s(x)) - \Delta_t(x)| + \sup_{x \in X}|\Delta_t(x)| 
\]
for every $t \ge s$, from \eqref{AA0} and \eqref{AAB} we get
\[
\limsup_{t \to \infty}\frac{1}{t}\sup_{s \in [0,\tau]}\sup_{x \in X}|\Delta_{t}(\phi_s(x))| = 0,
\]
and the claim is proved. \qed

Now for each $y \in Y$ there exist $x \in X$ and $s \in [0,\sup\tau]$ such that $y = \phi_s(x)$. Then
\[
|\Delta_{t}(y)|  = |\Delta_t(\phi_s(x))|\le \sup_{s \in [0,\tau]}\sup_{x \in X}|\Delta_t(\phi_s(x))| \quad \textrm{for each $t \ge 0$},
\]
which readily implies that 
\[
\limsup_{t \to \infty}\frac{1}{t}\|\Delta_t\|_{\infty} \le \limsup_{t \to \infty}\frac{1}{t}\sup_{s \in [0,\tau]}\sup_{x \in X}|\Delta_t(\phi_s(x))|= 0.
\]
Since, in particular, $\Delta_n$ is asymptotically aditive with respect to the map $\phi_1$ on $Y$, Lemma A.3 in \cite{FH10} guarantees that the limit $\lim_{t \to \infty}\frac{1}{t}\|\Delta_t\|_{\infty}$ exists. Therefore
\[
\lim_{t \to \infty}\frac{1}{t}\|\Delta_t\|_{\infty} = 0,
\]
and the theorem is proved. \qed

\begin{definition}
Let $X$ and $Y$ be topological spaces and consider the flows $\Phi$ on $X$ and $\Psi$ on $Y$. We say that $(X,\Phi)$ is \emph{topologically conjugate} to $(Y,\Psi)$ if there exists an homeomorphism $h: X \to Y$ such that $(h \circ \phi_t)(x) = (\psi_t \circ h)(x)$ for every $t \in \R$ and every $x \in X$. Moreover, we say that $(X,\Phi)$ is $C^{r}$-\emph{conjugate} to $(Y,\Psi)$ if $h$ is $C^{r}$. 
\end{definition}

\begin{corollary}\label{corB}
Let $\Phi = (\phi_t)_{t \in \R}$ be a continuous flow on a compact metric space $M$. Suppose that $\Phi$ is topologically conjugate to a suspension flow. Then, for every asymptotically additive family of continuous functions $a = (a_t)_{t \in \R}$ there exists a continuous function $b: M \to \R$ such that
\[
\lim_{t \to \infty}\frac{1}{t}\bigg\|a_t - \int_{0}^{t}(b \circ \phi_s)ds\bigg\|_{\infty} = 0.
\]
\end{corollary}

\begin{proof}
Let $\Psi$ be a suspension flow on a compact metric space $N$ and $h: M \to N$ the homeomorphism conjugating $(M,\Phi)$ and $(N,\Psi)$. Suppose that $a = (a_t)_{t \ge 0}$ is an asymptotically additive family of continuous functions with respect to $\Phi$. One can easily check that the family $(a_t \circ h^{-1})_{t \ge 0}$ is asymptotically additive with respect to $\Psi$. By Theorem \ref{main}, there exists a continuous function $\widetilde{b}: N \to \R$ such that 
\[
\lim_{t \to \infty}\frac{1}{t}\sup_{y \in N}\bigg|(a_t \circ h^{-1})(y) - \int_{0}^{t}(\widetilde{b} \circ \psi_s)(y)ds\bigg| = 0.
\]
Since for each $y \in N$ there exists a unique $x \in M$ such that $h(x) = y$, we have 
\[
\bigg|a_t(h^{-1}(y)) - \int_{0}^{t}(\widetilde{b} \circ \psi_s)(y)ds\bigg| = \bigg|a_t(x) - \int_{0}^{t}((\widetilde{b} \circ h) \circ \phi_s)(x)ds\bigg|.
\]
Hence
\[
\lim_{t \to \infty}\frac{1}{t}\sup_{x \in M}\bigg|a_t(x) - \int_{0}^{t}(b \circ \phi_s)(x)ds\bigg| = 0,
\]
where $b := \widetilde{b} \circ h: M \to \R$. 

\end{proof}

\begin{example}\label{TOR} Let $\mathbb{T}^{n}:= \R/\Z \times \cdots \times \R/\Z$ be the $n$-torus. Letting $\alpha \in \R^{n}$, the linear flow $\Phi = (\phi_t)_{t \in \R}$ on $\mathbb{T}^{n}$ in the direction $\alpha$ is defined by $\phi_t(x) = x + t\alpha \mod 1$. 

One can see that every linear flow on the $n$-torus is $C^{\infty}$- conjugate to a suspension flow (see for example \cite{FH20}). Then, it follows from Corollary \ref{corB} that asymptotically additive and additive families of continuous functions are physically equivalent with respect to the flow $\Phi$. \qed
\end{example}

\subsection{Hyperbolic flows and Markov partitions}\label{subsec32}
Let $\Phi =(\phi_t)_{t\in\R}$ be a $C^1$ flow on a smooth manifold~$M$.
A compact $\Phi$-invariant set $\Lambda \subset M$ is called a \emph{hyperbolic set for $\Phi$} if there exists a splitting
\[
T_{\Lambda}M = E^{s} \oplus E^{u} \oplus E^{0}
\]
and constants $c > 0$ and $\lambda \in (0,1)$ such that for each $x \in \Lambda$:	
\begin{enumerate}
	\item
	the vector $(d/dt)\phi_t(x)|_{t=0}$ generates $E^{0}(x)$;
	\item
	for each $t \in \R$ we have
	\[
	d_{x}\phi_t E^{s}(x) = E^{s}(\phi_t(x)) \quad \text{and} \quad d_{x}\phi_t E^{u}(x) = E^{u}(\phi_t(x));
	\]
	\item
	$\Vert d_x\phi_tv\Vert \le c\lambda^{t}\Vert v\Vert$ for $v \in E^{s}(x)$ and $t > 0$;
	\item
	$\Vert d_x\phi_{-t}v\Vert \le c\lambda^{t}\Vert v\Vert$ for $v \in E^{u}(x)$ and $t > 0$.
\end{enumerate}
Given a hyperbolic set $\Lambda$ for a flow $\Phi$, for each $x \in \Lambda$ and any sufficiently small $\epsilon > 0$ we define
\[
A^{s}(x) = \bigl\{y \in B(x,\epsilon): d(\phi_t(y),\phi_t(x)) \searrow 0 \ \text{when} \ t \to +\infty \bigr\}
\]
and
\[
A^{u}(x) = \bigl\{y \in B(x,\epsilon): d(\phi_t(y),\phi_t(x)) \searrow 0 \ \text{when} \ t \to -\infty \bigr\}.
\]
Moreover, let $V^{s}(x) \subset A^{s}(x)$ and $V^{u}(x) \subset A^{u}(x)$ be the largest connected components containing~$x$. These are smooth manifolds, called respectively \textit{(local) stable and unstable manifolds} of size $\epsilon$ at the point~$x$, satisfying:
\begin{enumerate}
	\item
	$T_xV^{s}(x) = E^{s}(x)$ and $T_xV^{u}(x) = E^{u}(x)$;
	\item
	for each $t > 0$ we have
	\[
	\phi_t(V^{s}(x)) \subset V^{s}(\phi_t(x)) \quad \text{and} \quad \phi_{-t}(V^{u}(x)) \subset V^{u}(\phi_{-t}(x));
	\]
	\item
	there exist $\kappa > 0$ and $\mu \in (0,1)$ such that for each $t > 0$ we have
	\[
	d(\phi_t(y),\phi_t(x)) \le \kappa \mu^{t}d(y,x) \quad \text{for} \ y \in V^{s}(x)
	\]
	and
	\[d(\phi_{-t}(y),\phi_{-t}(x)) \le \kappa \mu^{t}d(y,x) \quad \text{for} \ y \in V^{u}(x).
	\]
\end{enumerate}
We recall that a set $\Lambda$ is said to be \textit{locally maximal} (with respect to a flow~$\Phi$) if there exists an open neighborhood $U$ of $\Lambda$ such that
\[
\Lambda = \bigcap_{t \in \R}\phi_t(U).
\]
Given a locally maximal hyperbolic set $\Lambda$ and a sufficiently small $\epsilon > 0$, there exists $\delta > 0$ such that if $x, y \in \Lambda$ satisfy $d(x,y) \le \delta$, then there exists a unique $t=t(x,y) \in [-\epsilon,\epsilon]$ such that
\[
[x,y]:= V^{s}(\phi_t(x)) \cap V^{u}(x)
\]
is a single point in~$\Lambda$.

Now let us recall the notion of Markov partitions for continuous-time dynamical systems. Consider an open smooth disk $D \subset M$ of dimension $\dim M - 1$ that is transverse to~$\Phi$ and take $x \in D$. Let $U(x)$ be an open neighborhood of $x$ diffeomorphic to $D \times (-\epsilon,\epsilon)$. Then the projection $\pi_D\colon U(x) \to D$ defined by $\pi_D(\phi_t(y)) = y$ is differentiable.
We say that a closed set $R \subset \Lambda \cap D$ is a \textit{rectangle} if $R = \overline{\interior R}$ and $\pi_D([x,y]) \in R$ for $x, y \in R$.

Consider rectangles $R_1,\ldots,R_k \subset \Lambda$ (each contained in some open smooth disk transverse to the flow) such that
\[
R_i \cap R_j = \partial R_i \cap \partial R_j \quad \text{for} \ i\ne j.
\]
Let $Z=\bigcup_{i=1}^{k}R_i$. We assume that there exists $\epsilon > 0$ such that:
\begin{enumerate}
\item
$\Lambda = \bigcup_{t \in [0,\epsilon]}\phi_t(Z)$;
\item
whenever $i \ne j$, either
\[
\phi_t(R_i)\cap R_j = \emptyset \quad \text{for all} \ t \in [0,\epsilon]
\]
or
\[
\phi_t(R_j)\cap R_i = \emptyset \quad \text{for all} \ t \in [0,\epsilon].
\]
\end{enumerate}
Now define the function $\tau\colon \Lambda \to \R^+_0$ by
\[
\tau(x) = \min\{t>0: \phi_t(x) \in Z\},
\]
and the map $T\colon \Lambda \to Z$ by
\begin{equation}\label{tpx}
T(x) = \phi_{\tau(x)}(x).
\end{equation}
The restriction $T_Z$ of $T$ to $Z$ is invertible and we have $T^n(x) = \phi_{\tau_n(x)}(x)$, where
\[
\tau_n(x) = \sum_{i=0}^{n-1}\tau(T^{i}(x)).
\]
The collection $R_1,\ldots,R_k$ is said to be a \emph{Markov partition} for $\Phi$ on $\Lambda$ if
\[T(\interior (V^{s}(x) \cap R_i)) \subset \interior(V^{s}(T(x))\cap R_j)
\]
and
\[
T^{-1}(\interior(V^{u}(T(x)) \cap R_j)) \subset \interior(V^{u}(x)\cap R_i)
\]
for every $x \in \interior T(R_i)\cap \interior R_j$ and $i,j=1,\ldots,k$. By work of Bowen \cite{Bow73} and Ratner \cite{Rat73}, any locally maximal hyperbolic set $\Lambda$ has Markov partitions of arbitrarily small diameter and the function $\tau$ is H\"older continuous on each domain of continuity.

Given a Markov partition $R_1,\ldots,R_k$ for a flow $\Phi$ on a locally maximal hyperbolic set $\Lambda$, we consider the $k \times k$ matrix $A$ with entries
\[
a_{ij}= \begin{cases} 1 & \text{if} \ \interior T(R_i) \cap R_j \ne \emptyset,\\ 0 & \text{otherwise}, \end{cases}
\]
where $T$ is the map in \eqref{tpx}. We also consider the set
\[
\Sigma_A= \bigl\{(\cdots i_{-1} i_0 i_1 \cdots): a_{i_n i_{n+1}} = 1 \text{ for } n \in \Z\bigr\}\subset \{1,\ldots,k\}^{\Z}
\]
and the shift map $\sigma\colon \Sigma_A \to \Sigma_A$ defined by $\sigma(\cdots i_0 \cdots) = (\cdots j_0 \cdots)$, where $j_n = i_{n+1}$ for each $n \in \Z$.
We denote by $\Sigma_n$ the set of $\Sigma_A$-admissible sequences of length $n$, that is, the finite sequences $(i_1\cdots i_n)$ for which there exists $(\cdots j_0j_1j_2 \cdots) \in \Sigma_A$ such that $(i_1\dots i_n) = (j_1\cdots j_n)$.
Finally, we define a \emph{coding map} $\pi\colon\Sigma_A \to Z$ by
\[
\pi(\cdots i_0 \cdots) = \bigcap_{n \in \Z}R_{i_{-n} \cdots i_n},
\]
where $R_{i_{-n} \cdots i_n} = \bigcap_{l=-n}^n\overline{T_Z^{-l}\interior R_{i_l}}$.
The following properties hold:
\begin{enumerate}
	\item
	$\pi \circ \sigma = T \circ \pi$;
	\item
	$\pi$ is H\"older continuous and onto;
	\item
	$\pi$ is one-to-one on a full measure set with respect to any ergodic measure of full support and on a residual set.
\end{enumerate}
Given $\beta > 1$, we equip $\Sigma_A$ with the distance $d_\beta$ defined by
\[
d_\beta(\omega,\omega') = \begin{cases} \beta^{-n} &\mbox{if } \omega \ne \omega', \\ 0 &\mbox{if } \omega = \omega', \end{cases}
\]
where $n = n(\omega,\omega') \in \N\cup \{0\}$ is the smallest integer such that $i_n(\omega) \ne i_n(\omega')$ or $i_{-n}(\omega) \ne i_{-n}(\omega')$.
One can always choose $\beta$ so that $\tau \circ \pi$ is Lipschitz.

Now let $\nu$ be a $T_Z$-invariant probability measure on~$Z$. One can show that $\nu$ induces a $\Phi$-invariant probability measure $\mu$ on $\Lambda$ such that
\begin{equation}\label{eq13}
\int_{\Lambda}g\,d\mu = \frac{\int_{Z}\int_{0}^{\tau(x)}(g\circ\phi_s)(x)\,ds\,d\nu}{\int_{Z}\tau \,d\nu}
\end{equation}
for any continuous function $g\colon \Lambda \to \R$. In fact, any $\Phi$-invariant probability measure $\mu$ on $\Lambda$ is of this form for some $T_Z$-invariant probability measure $\nu$ on~$Z$. Abramov's entropy formula says that
\begin{equation}\label{eq14}
h_{\mu}(\Phi) = \frac{h_{\nu}(T_Z)}{\int_{Z}\tau \,d\nu}.
\end{equation}
By \eqref{eq13} and \eqref{eq14} we obtain
\begin{equation}\label{eq15}
h_{\mu}(\Phi) + \int_{\Lambda}g\,d\mu = \frac{h_\nu(T_Z) + \int_{Z}I_g\,d\nu}{\int_{Z}\tau \,d\nu},
\end{equation}
where $I_g(x) = \int_{0}^{\tau(x)}(g\circ\phi_s)ds$. Since $\tau>0$ on~$Z$, it follows from \eqref{eq15} that
\[
P_{\Phi}(g) = 0 \quad \text{if and only if} \quad P_{T_Z}(I_g) = 0,
\]
where $P_{\Phi}(g)$ is the topological pressure of $g$ with respect to $\Phi$ and $P_{T_Z}(I_g)$ is the topological pressure of $I_g$ with respect to the map~$T_Z$. When $P_{\Phi}(g) = 0$, this implies that $\mu$ is an equilibrium measure for $g$ if and only if $\nu$~is an equilibrium measure for~$I_g$.

As a direct consequence of the existence of Markov partitions for locally maximal hyperbolic sets together with Theorem \ref{main}, we obtain the following result:

\begin{corollary}\label{HYP}
Let $\Lambda$ be a locally maximal hyperbolic set for a $C^1$ flow $\Phi = (\phi_t)_{t \in \R}$ and suppose that $a = (a_t)_{t \ge 0}$ is an asymptotically additive family of continuous functions with respect to $\Phi$. Then, there exists a continuous function $b: \Lambda \to \R$ such that
\[
\lim_{t \to \infty}\frac{1}{t}\bigg\|a_t - \int_{0}^{t}(b\circ \phi_s)ds\bigg\|_{\infty} = 0.
\]
\end{corollary}

\subsection{Expansive flows}

Let $(X,d)$ be a metric space and $T: X \to X$ a dynamical system. $T$ is said to be \emph{expansive} if there exists $\delta > 0$ such that $d(T^{n}(x), T^{n}(y)) < \delta$ for all $n \in \mathbb{Z}$ implies $x = y$. 
 
The definition for continuous-time dynamical systems is more refined. A flow $\Phi$ on $X$ is said to be \emph{expansive} if for each $\epsilon > 0$ there exists $\delta > 0$ such that if $d(\phi_t(x), \phi_{s(t)}(y)) < \delta$ for all $t \in \R$, for points $x$ and $y$ and a continuous map $s: \R \to \R$ with $s(0) = 0$, then there exists a time $|t| < \epsilon$ such that $\phi_t(x) = y$ (see \cite{BW72}).

Let $X$ and $Y$ be metric spaces and consider the flows $\Phi$ on $X$ and $\Psi$ on $Y$. We say that $S \subset X$ is a $\emph{full set}$ when $\mu(S) = 0$ for all $\Phi$-invariant measure $\mu$. Suppose that $\pi: Y \to X$ is a \emph{topological extension} from $(X,\Phi)$ to $(Y,\Psi)$, that is, a surjective map topologically conjugating $(\Phi,X)$ and $(\Psi,Y)$. The extension $\pi: Y \to X$ is said to be \emph{strongly isomorphic} if there exists a full set $E \subset X$ such that $\pi: \pi^{-1}(E) \subset Y \to X$ is one-to-one (see for example \cite{Bur19}). 

Advancing the main results in \cite{Bur19}, recently Gutman and Shi proved the following: 

\begin{theorem}[{\cite[Theorem~B]{GS22}}]\label{B}
Let $X$ be a compact finite-dimensional space and $\Phi$ an expansive flow on $X$. Then, $(X,\Phi)$ is strongly isomorphic to a suspension flow over a subshift of finite type. 
\end{theorem}

This result together with Theorem \ref{main} gives the following:
\begin{theorem}\label{EXP}
	Let $\Phi = (\phi_t)_{t \in \R}$ be a continuous expansive flow on a compact finite-dimensional metric space $M$, and let $a = (a_t)_{t \ge 0}$ be an asymptotically additive family of continuous functions. Then, there exists a continuous function $b: M \to \R$ and a full set $N \subset M$ such that
	\begin{equation}\label{SYM}
	\lim_{t \to \infty}\frac{1}{t}\sup_{x \in N}\bigg|a_t(x) - \int_{0}^{t}(b\circ \phi_s)(x)ds\bigg| = 0. 
	\end{equation}
	Moreover, if $\Phi$ admits an invariant measure with full support then 
	\[
	\lim_{t \to \infty}\frac{1}{t}\bigg\|a_t - \int_{0}^{t}(b\circ \phi_s)ds\bigg\|_{\infty} = 0.
	\]
\end{theorem}
\begin{proof}
	By Theorem \ref{B}, there exists a full set $N \subset X$ such that $(N,\Phi)$ and $(Y,\Psi)$ are topologically conjugate, where $\Psi$ is a suspension flow over a subshift of finite type. Then, \eqref{SYM} follows directly from Corollary \ref{corB}. 
	
	Now suppose that $\nu \in \cM(\Phi)$ is a measure with full support. Then, one can see that $N$ is dense on the whole space $M$. Since the function
	\[
	x \mapsto D_t(x) := \bigg|a_t(x) - \int_{0}^{t}(b\circ \phi_s)(x)ds\bigg|
	\] 
	is continuous for every $t \ge 0$, we have $\sup D_t(M) = \sup D_t(\overline{N}) \le \sup \overline{D_t(N)} = \sup D_t(N)$ for every $t \ge 0$. This together with \eqref{SYM} yield
	\[
	\lim_{t \to \infty}\frac{1}{t}\bigg\|a_t - \int_{0}^{t}(b\circ \phi_s)ds\bigg\|_{\infty}=\lim_{t \to \infty}\frac{1}{t}\sup_{x \in M}D_t(x) \le \lim_{t \to \infty}\frac{1}{t}\sup_{x \in N}D_t(x) = 0,
	\]
	and the theorem is proved.
\end{proof}

The following result is a weaker notion of equivalence in the case of expansive flows.
\begin{corollary}\label{EEX}
Let $\Phi = (\phi_t)_{t \in \R}$ be a continuous expansive flow on a compact finite-dimensional metric space $M$, and let $a = (a_t)_{t \ge 0}$ be an asymptotically additive family of continuous functions. Then, there exists a continuous function $b: M \to \R$ such that
\[
\lim_{t \to \infty}\frac{1}{t}\int_{M}a_t d\mu = \int_{M} b d\mu 
\]
for every measure $\mu \in \cM(\Phi)$. 
\end{corollary}

\begin{proof}
It follows directly from Theorem \ref{EXP} and Birkhoff's ergodic theorem. 
\end{proof}

Observe that volume preserving continuous expansive flows on compact finite-dimensional manifolds satisfy all the hypotheses of Theorem \ref{EXP}.

\subsection{Some examples: conformal and non-conformal hyperbolic flows}\label{CON}

We will now introduce a source of asymptotically additive families of continuous potentials.

\subsubsection{Conformal flows}

We say that a $C^1$ flow $\Phi$ is \emph{conformal} on a hyperbolic set $\Lambda$ if there exist continuous functions $Q^s, Q^u\colon\Lambda \times \R \to \R$ such that
\[
d_x\phi_t|E^s(x) = Q^s(x,t)J^s(x,t)\quad \text{and} \quad d_x\phi_t|E^{u}(x) = Q^u(x,t)J^u(x,t)
\]
for every $x \in \Lambda$ and $t \in \R$, where
\[
J^s(x,t)\colon E^s(x) \to E^s(\phi_t(x))\quad\text{and}\quad J^u(x,t)\colon E^u(x) \to E^u(\phi_t(x))
\]
are isometries. For example, if
\[
\dim E^s(x) =\dim E^u(x) =1\quad\text{for} \ x \in \Lambda,
\]
then the flow is conformal on~$\Lambda$.
Proceeding as in \cite{PS01} we define:
\begin{equation}\label{XIS}
\Xi_s(x):
= \frac{\partial}{\partial t}\log |Q^s(x,t)|_{t=0} 
= \frac{\partial}{\partial t}\log \|d_x\phi_t|E^s(x)\|_{t=0} 
= \lim_{t \to 0}\frac{\log\|d_x\phi_t|E^s(x)\|}{t}
\end{equation}
and
\begin{equation}\label{XIU}
\Xi_u(x):
= \frac{\partial}{\partial t}\log |Q^u(x,t)|_{t = 0} 
= \frac{\partial}{\partial t}\log \|d_x\phi_t|E^u(x)\|_{t = 0}
= \lim_{t \to 0}\frac{\log\|d_x\phi_t|E^u(x)\|}{t}.
\end{equation}
Since the flow $\Phi$ is of class $C^1$, using $2$-norms one can write
\[
\begin{split}
\lim_{t \to 0}\frac{\log\|d_x\phi_t|E^s(x)\|}{t} &= \lim_{t \to 0}\frac{\log (\|d_x\phi_t|E^s(x)\|^2)}{2t} = \lim_{t \to 0}\frac{\langle d_x\phi_t|E^s(x), \frac{\partial}{\partial t}(d_x\phi_t|E^s(x))\rangle}{\|d_x\phi_t|E^u(x)\|^2}\\
& = \left\langle \text{Id}|E^s(x), \frac{\partial}{\partial t}(d_x\phi_t|E^s(x))|_{t=0}\right\rangle
\end{split}
\]
and, similarly,
\[
\lim_{t \to 0}\frac{\log\|d_x\phi_t|E^u(x)\|}{t} = \left\langle \text{Id}|E^u(x), \frac{\partial}{\partial t}(d_x\phi_t|E^u(x))|_{t=0}\right\rangle.
\]
In particular, the functions $\Xi_s$ and $\Xi_u$ are well defined. 
For an adapted norm $\lVert\cdot\rVert$ (that is, a norm for which one can take $c = 1$ in the definition of a hyperbolic set), we obtain
\[
\Xi_s(x) = \lim_{t \to 0^+}\frac{\log\|d_x\phi_t|E^s(x)\|}{t} \le \log\lambda < 0
\]
and
\[
\Xi_u(x) = \lim_{t \to 0^+}\frac{\log\|d_x\phi_t|E^u(x)\|}{t} \ge -\log\lambda>0
\]
for all $x \in \Lambda$. Moreover, for every $x \in \Lambda$ and $t \in \R$, it follows from \eqref{XIS} and \eqref{XIU} that
\begin{equation}\label{hju}
\|d_x\phi_t v\| = \|v\|\exp\biggl(\int_{0}^{t}\Xi_s(\phi_{\tau}(x))\,d\tau\biggr) \quad \text{for} \ v \in E^{s}(x)
\end{equation}
and
\begin{equation}\label{hjs}
\|d_x\phi_t v\| = \|v\|\exp\biggl(\int_{0}^{t}\Xi_u(\phi_{\tau}(x))\,d\tau\biggr) \quad \text{for} \ v \in E^{u}(x).
\end{equation}

In this case, notice that
\[
(\log\|d_x\phi_t|E^s(x)\|)_{t \ge 0} \quad \textrm{and} \quad (\log\|d_x\phi_t|E^u(x)\|)_{t \ge 0}
\]
are additive families of continuous functions with respect to $\Phi$. 

\subsubsection{Non-conformal flows with bounded distortion}\label{bdf}

Let $\Lambda$ be an hyperbolic set for a $C^1$ flow $\Phi = (\phi_t)_{t \in \R}$. Moreover, let $E^{s}(x)$ and $E^u(x)$ be the stable and unstable spaces at~$x$. We say that $\Phi$ has \emph{bounded distortion} (in the sense of \cite{PS01}) if there exist constants $C_1 > 0$, $C_2>0$ and H\"older continuous functions $b^s, b^u\colon \Lambda \to \R$ such that
\[
C_1 \|v\|\exp\int_0^t(b^s\circ\phi_{\tau})(x)\,d\tau \le \|d_x \phi_tv\| \le C_2\|v\|\exp\int_0^t(b^s\circ\phi_{\tau})(x)\,d\tau
\]
for $v \in E^s(x)$, and
\[
C_1 \|v\|\exp\int_0^t(b^u\circ\phi_{\tau})(x)\,d\tau \le \|d_x \phi_tv\| \le C_2\|v\|\exp\int_0^t(b^u\circ\phi_{\tau})(x)\,d\tau
\]
for $v \in E^u(x)$. In this case one can easily verify that the families $a^s= (a^s_t)_{t \ge 0}$ and $a^u= (a^u_t)_{t \ge 0}$ given by
\[
a^s_t(x) = \log\|d_x\phi_t|_{E^s(x)}\| \quad \text{and} \quad a^u_t(x) = \log\|d_x\phi_t|_{E^u(x)}\|
\]
are almost additive with respect to~$\Phi$ and satisfy
\[
\lim_{t \to \infty}\frac{1}{t}\biggl\|a^s_t -\int_0^t(b^s\circ\phi_{\tau})\,d\tau\biggr\|_{\infty} = 0
\]
and
\[
\lim_{t \to \infty}\frac{1}{t}\biggl\|a^u_t -\int_0^t(b^u\circ\phi_{\tau})\,d\tau\biggr\|_{\infty} = 0.
\]

\subsubsection{Quasiconformal flows}
Now let $\Lambda$ be an hyperbolic set for a $C^{1}$ flow $\Phi = (\phi_t)_{t \in \R}$, and consider the functions $K^{s}, K^{u}: \Lambda \times\R \to~\R$ given by 
\begin{equation}\label{KS}
K^{s}(x,t) = \frac{\max\{\|d_x \phi_t v\|: v \in E^{s}(x), \|v\| =1\}}{\min\{\|d_x\phi_t v\|: v \in E^{s}(x), \|v\| =1\}}
\end{equation}
and
\begin{equation}\label{KU}
K^{u}(x,t) = \frac{\max\{\|d_x \phi_t v\|: v \in E^{u}(x), \|v\| =1\}}{\min\{\|d_x\phi_t v\|: v \in E^{u}(x), \|v\| =1\}}.
\end{equation}

We say that $\Phi$ is \emph{uniformly quasiconformal} if the functions $K^{s}$ and $K^{u}$ are uniformly bounded for all $x \in \Lambda$ and $t \in \R$. Observe that when $\Phi$ is conformal on $\Lambda$, it follows directly from \eqref{hju} and \eqref{hjs} that $K^{s}(x,t) = 1$ and $K^{u}(x,t) = 1$ for all $x \in \Lambda$ and all $t \in \R$. The notion of a uniformly quasiconformal hyperbolic map is analogous (see \cite{Sad05}). Observe that conformal flows are quasiconformal and an Anosov diffeomorphism is quasiconformal if and only if its suspension flow is quasiconformal (see for example \cite{Fan05}). It follows from \eqref{KS} and \eqref{KU} that 
\[
\lim_{t \to \infty}\frac{1}{t}\|K^{s}(x,t)\|_{\infty} = 0 \quad \textrm{and} \quad \lim_{t \to \infty}\frac{1}{t}\|K^{u}(x,t)\|_{\infty} = 0,
\]
which readily implies that $(K^s(x,t))_{t \ge 0}$ and $(K^u(x,t))_{t \ge 0}$ are asymptotically additive families with respect to $\Phi$. 

\subsubsection{Average conformal flows} 

Inspired by previous work \cite{BCH10}, it was introduced in \cite{WWCZ20} a type of non-conformal hyperbolic maps, which can also be seen as a generalization of quasiconformal maps. Let $\Lambda$ be an hyperbolic set for a diffeomorphism $f: \Lambda \to \Lambda$. The set $\Lambda$ is called an \emph{average conformal} hyperbolic set for $f$ if it admits exactly two unique Lyapunov exponents, one strictly positive and the other strictly negative. 

If $\Phi = (\phi_t)_{t \in \R}$ is a suspension flow over an average conformal hyperbolic map, it follows from Lemma 2.3 in \cite{WWCZ20} that
\[
\begin{split}
&\lim_{t \to \infty}\frac{1}{t}\|\log \|d_x\phi_t|_{E^{s}(x)}\| -\log \|(d_x\phi_t|_{E^{s}(x)})^{-1}\|^{-1}\|_{\infty} = 0, \\
& \lim_{t \to \infty}\frac{1}{t}\|\log \|d_x\phi_t|_{E^{u}(x)}\| -\log \|(d_x\phi_t|_{E^{u}(x)})^{-1}\|^{-1}\|_{\infty} = 0 \quad \textrm{and}
\end{split}
\]
\[
\begin{split}
&\|(d_x\phi_t|_{E^{s}(x)})^{-1}\|^{-1} \le |\det (d_x\phi_t|_{E^s(x)})|^{\frac{1}{d_s}} \le \|d_x\phi_t|_{E^{s}(x)}\|,\\
& \|(d_x\phi_t|_{E^{u}(x)})^{-1}\|^{-1} \le |\det (d_x\phi_t|_{E^u(x)})|^{\frac{1}{d_u}} \le \|d_x\phi_t|_{E^{u}(x)}\|,
\end{split}
\]
where $d_s:= \dim E^{s}$ and $d^{u} := \dim E^{u}$. 
Since $(\det (d_x\phi_t|_{E^u(x)})_{t \ge 0}$ and $(\det (d_x\phi_t|_{E^s(x)})_{t \ge 0}$ are additive families with respect to $\Phi$, one can see that $a^{i}:= (\|\log \|d_x\phi_t|_{E^{i}(x)}\|)_{t \ge 0}$ and $b^{i}:=(\|\log \|d_x\phi_t|_{E^{i}(x)}\| - \log \|(d_x\phi_t|_{E^{i}(x)})^{-1}\|^{-1})_{t \ge 0}$ are asymptotically additive families of continuous functions with respect to $\Phi$ for $i \in \{u,s\}$.

\subsection{The general physical equivalence problem}\label{GER}

We recall that a family of functions $a = (a_t)_{t \ge 0}$ is said to have \textit{bounded variation} if there exist $\kappa> 0$ and $\epsilon > 0$ such that
\[
|a_t(x) - a_t(y)| < \kappa\quad\text{whenever} \ y\in B_t(x,\epsilon).
\]

Inspired by the examples in subsection \ref{CON} and by the work \cite{Cun20}, one can ask the following questions:

\textbf{Question A:} Given an asymptotically additive family of continuous functions $a=(a_t)_{t \ge 0}$ with respect to a continuous flow $\Phi = (\phi_t)_{t \in \R}$ on a compact metric space~$X$, is there any continuous function $b\colon X \to \R$ such~that
\begin{equation}\label{PEQ}
\lim_{t \to \infty}\frac{1}{t}\biggl\|a_t -\int_{0}^t(b\circ\phi_{s})ds\biggr\|_{\infty} = 0 \ ?
\end{equation}

\textbf{Question B:} When $\Phi|_{\Lambda}$ is a hyperbolic flow and the family $a$ is almost additive with bounded variation, is there any continuous function $b\colon \Lambda \to \R$ where the additive family 
\[
S_tb := \int_{0}^{t}(b\circ \phi_s)ds 
\]
has bounded variation and also satisfy \eqref{PEQ} ?

As we showed, Theorem \ref{main} answers the \textbf{Question A} in the case of suspension flows and, in particular, in the case of hyperbolic flows. Moreover, Theorem \ref{EXP} gives a setup where we also can positively answer the \textbf{Question A} in the case of expansive flows. Furthermore, observe that \textbf{Question B} is the more general Bowen regularity problem for flows and it is open even in the case of discrete-time dynamical systems (see \cite{Cun20} and \cite{HS24}). In general, as discussed before, positive answers to these questions are very useful for some extensions of the thermodynamic formalism and multifractal analysis for flows (see subsection \ref{APP}). On the other hand, considering the same problem with respect to H\"older regularity, \textbf{Question~B} has a negative answer in both discrete and continuous-time setups (see \cite{HS24} and subsection \ref{BHR}). 

Let us now give some directions on how to approach this problem in general. Let $\Phi = (\phi_t)_{t \in \R}$ be a continuous flow and let $a = (a_t)_{t \ge 0}$ be an asymptotically additive family of continuous functions with respect to $\Phi$. For any function $c: X \to \R$, we have 
\begin{equation}\label{EMB}
\int_{0}^{n}(c \circ \phi_s)ds = \sum_{k=0}^{n-1}(\widetilde{c} \circ \phi^{k}_1) \quad \textrm{for every $n \ge 1$}
\end{equation}
where $\widetilde{c} := \int_{0}^{1}(c \circ \phi_s)ds$, and one can check that $(a_n)_{n \ge 1}$ is an asymptotically additive sequence with respect to the map $\phi_1$. Then, by Theorem 1.2 in \cite{Cun20}, there exists a continuous function $\widetilde{b}: X \to \R$ such that 
\[
\lim_{n \to \infty}\frac{1}{n}\bigg\|a_n - \sum_{k=0}^{n-1}\widetilde{b} \circ \phi^{k}_1\bigg\|_{\infty} = 0.
\]

If there exists a continuous function $b: X \to \R$ such that $\widetilde{b} = \int_{0}^{1}(b \circ \phi_s)ds$, from \eqref{EMB} we readily obtain that 
\[
\lim_{n \to \infty}\frac{1}{t}\bigg\|a_t - \int_{0}^{t}(b \circ \phi_s)ds\bigg\|_{\infty} = 0.
\]

Just for illustration, here it goes a simple example: 

\begin{example}\label{EXA}
Let $\Phi = (\phi_t)_{t \in \R}$ be the flow $\phi_t(x) = e^{t}x$ on $\R^{+}:= (0,\infty]$. If $b(x) = 1/2 + \log x$, then for $a(x) = \log x$, we have 
\[
\int_{0}^{1}(a\circ \phi_s)(x)ds = b(x) \quad \textrm{for every $x \in \R^{+}$.}
\]
\qed
\end{example}

A less (but still) simple example:

\begin{example}\label{EXB}
Let $\Phi = (\phi_t)_{t \in \R}$ be the flow $\phi_t(x) = e^{t}x$ defined on $\R$. If $b: \R \to \R$ is given by 
\[
b(x) = p(x)\bigg(\frac{e^{d} - 1}{d}\bigg),
\]
where $p: \R \to \R$ is an homogeneous polynomial of degree $d$, then 
\[
\int_{0}^{1}(p\circ \phi_s)(x)ds = b(x) \quad \textrm{for every $x \in \R$.}
\]
\qed
\end{example}

Following this approach, we can ask:

\textbf{Question C.} Given a continuous flow $\Phi$ and a continuous function $\widetilde{b}: X \to \R$, is there any continuous function $b:X \to \R$ satisfying
\[
\widetilde{b}(x) = \int_{0}^{1}(b \circ \phi_s)(x)ds \quad \textrm{for every $x \in X$ ?}
\]

\textbf{Question D.} If we cannot give a positive answer to the previous question in general, what kind of functions and flows satisfy it ?

In order to extend the theory of ergodic optimization for flows, the \textbf{Question C} was also posted in \cite{BHVZ21}. 

The following example indicates a negative answer to \textbf{Question C} in the case where the function $\widetilde{b}: X \to \R$ is bounded. 
\begin{example}
Let $X$ be a compact metric space and $\Phi = (\phi_t)_{t \in \R}$ a continuous flow on $X$. Let $x^{*} \in X$ be a non-fixed point of $\Phi$ and consider $\widetilde{b}:= \mathbbm{1}_{x^*}: X \to \R$ the characteristic function of the set $\{x^{*}\}$. Suppose there exists a bounded function $b: X \to \R$ such that 
\begin{equation}\label{EQA}
\int_{0}^{1}(b\circ \phi_s)(x)ds = \widetilde{b}(x) \quad \textrm{for all $x \in X$}.
\end{equation}
This implies that 
\[
\int_{0}^{1}(b\circ \phi_s)(x^{*})ds = \widetilde{b}(x^{*}) = 1 \quad \textrm{and} \quad  \int_{0}^{1}(b\circ \phi_s)(x)ds = \widetilde{b}(x) = 0 \quad \textrm{for every $x \neq x^{*}$}.
\]

Given $\delta > 0$, consider the point $z:= \phi_{\delta}(x^{*})$. Since $x^{*}$ is not a fixed point, $z \neq x^{*}$ and we have
\[
0 = \int_{0}^{1}(b\circ \phi_s)(z)ds = \int_{\delta}^{\delta + 1}(b\circ \phi_s)(x^{*})ds 
\]
Then
\[
\begin{split}
\int_{0}^{1}(b\circ \phi_s)(x^{*})ds &= \int_{0}^{\delta}(b\circ \phi_s)(x^{*})ds + \int_{\delta}^{\delta+1}(b\circ \phi_s)(x^{*})ds + \int_{\delta +1}^{1}(b\circ \phi_s)(x^{*})ds \\
& = \int_{0}^{\delta}(b\circ \phi_s)(x^{*})ds + \int_{\delta +1}^{1}(b\circ \phi_s)(x^{*})ds,
\end{split}
\]
which readily implies that 
\[
\bigg|\int_{0}^{1}(b\circ \phi_s)(x^{*})ds\bigg| \le 2\delta\|b\|_{\infty}.
\]
Since $\delta$ is arbitrary and the function $b$ is bounded, taking $\delta < 1/(2\|b\|_{\infty})$ we obtain 
\[
1 = |\widetilde{b}(x^{*})|=\bigg|\int_{0}^{1}(b\circ \phi_s)(x^{*})ds\bigg| < 1,
\]
which is a contradiction. In particular, since $X$ is compact, there is no continuous function $b: X \to \R$ satisfying \eqref{EQA}. \qed
\end{example}

\begin{proposition}\label{BBP}
Let $\Phi$ be a continuous flow on a metric space $X$. Suppose that for a given continuous function $\widetilde{b}: X \to \R$ there exists a continuous function $b: X \to \R$ such that 
\[
\widetilde{b}(x) = \int_{0}^{1}(b \circ \phi_s)(x)ds \quad \textrm{for every $x \in X$}.
\]
Then, we have
\begin{equation}\label{DOT}
\lim_{t \to 0}\frac{\widetilde{b}(\phi_t(x)) - \widetilde{b}(x)}{t} = b(\phi_1(x)) - b(x)  \quad \textrm{for every $x \in X$}.
\end{equation}
Moreover,
\[
\int_{X}\bigg(\lim_{t \to 0}\frac{\widetilde{b}(\phi_t(x)) - \widetilde{b}(x)}{t}\bigg)d\mu(x) = 0 \quad \textrm{for every $\mu \in \cM(\phi_1)$}.
\]
\end{proposition}
\begin{proof}
Notice that 
\begin{equation}\label{BB}
\frac{\widetilde{b} \circ \phi_t - \widetilde{b}}{t} = \frac{1}{t}\bigg(\int_{t}^{t+1}(b \circ \phi_s)ds - \int_{0}^{1}(b \circ \phi_s)ds\bigg) \quad \textrm{for all $t > 0$}.
\end{equation}
One can check that the function
\[
t \mapsto I(t):= \int_{t}^{t+1}(b \circ \phi_s)ds - \int_{0}^{1}(b \circ \phi_s)ds
\]
is uniformly continuous on $[0,\infty)$, differentiable on $(0,\infty)$ and satisfies $\lim_{t \to 0}I(t) = 0$. Hence, by the L'H\^ospital's rule, we obtain that 
\[
\lim_{t \to 0}\frac{I(t)}{t} = \lim_{t \to 0}(b\circ\phi_{t+1} - b \circ \phi_t) = b \circ \phi_1 - b.
\]
This together with \eqref{BB} proves the result. 
\end{proof} 

\begin{remark}
One can check directly that the functions in Example \ref{EXA} and Example \ref{EXB} satisfy the condition \eqref{DOT}.
\end{remark}

Observe that if there exists a constant $\beta \in \R$ such that
\[
\lim_{t \to 0}\frac{\widetilde{b}(\phi_t(x)) - \widetilde{b}(x)}{t} = \beta \quad \textrm{for every $x \in X$},
\]
then Proposition \ref{BBP} says that $\beta = 0$. This fact is an inspiration for our \emph{negative} answer to \textbf{Question C}:  

\begin{example} (Counter-example)\label{CE1}. Let $\Phi = (\phi_t)_{t \in \R}$ be the linear flow on $\mathbb{T}^{2}$ given by $\phi_t(x,y) = (x + t\alpha_1 \bmod 1, y+t\alpha_2 \bmod 1)$, with $0 < \alpha_1 + \alpha_2 < 1$. Let $\widetilde{b}: \mathbb{T}^{2} \to \R$ be the continuous function given by $\widetilde{b}(x,y) = (x + y) \bmod 1$. One can check that
\[
\begin{split}
\lim_{t \to 0}\frac{\widetilde{b}(\phi_t(x,y)) - \widetilde{b}(x,y)}{t}	&= \lim_{t \to 0}\frac{(x + y + t(\alpha_1 + \alpha_2)) \bmod 1 - (x+y)\bmod 1}{t}\\
& = \lim_{t \to 0}\frac{t(\alpha_1 + \alpha_2)\bmod 1}{t} = \lim_{t \to 0}\frac{t(\alpha_1+\alpha_2)}{t} = \alpha_1 + \alpha_2
\end{split}
\]	
for every $x \in \mathbb{T}^{2}$. 	
If there exists a continuous function $b: \mathbb{T}^{2} \to \R$ such that $\widetilde{b} = \int_{0}^{1}(b \circ \phi_s)ds$, then Proposition~\ref{BBP} immediately gives that $\alpha_1 + \alpha_2 = 0$, which is a contradiction. Therefore, there is no continuous function $b: \mathbb{T}^{2} \to \R$ such that 
\[
(x + y) \bmod 1 = \int_{0}^{1}b(\phi_s(x,y))ds.
\]

\qed
\end{example}

\begin{remark}
Example \ref{CE1} is a negative answer to the general embedding problem presented in \cite{BHVZ21}. This counter-example can be naturally extended to the $n$-torus $\mathbb{T}^{n}$. Observe that even though we cannot guarantee in general a positive answer for \textbf{Question~C} in the case of linear flows on $\mathbb{T}^{n}$, in Example \ref{TOR} we were able to guarantee the physical equivalence between additive and asymptotically additive families of continuous potentials. 
\end{remark}

Notice that the converse of Proposition \ref{BBP} is false. In fact, still considering Example \ref{CE1}, if we take $b(x,y) = (x+y)\bmod 1$, then $b(\phi_1(x,y)) - b(x,y) = (\alpha_1 + \alpha_2) \bmod 1 = \alpha_1 + \alpha_2$ for every $x \in \mathbb{T}^{2}$. But 
\[
\int_{0}^{1}b(\phi_s(x,y))ds = \bigg[(x + y)\bmod1 + \frac{\alpha_1+\alpha_2}{2}\bigg]\bmod 1 \neq \widetilde{b}(x,y). 
\]

We will now recall the notions of cohomology for flows and maps.
\begin{definition}
Let $X$ be a metric space, $\Psi$ a flow on $X$ and $T: X \to X$ a map. We say that a function $g: X \to \R$ is \emph{$\Psi$-cohomologous} to a function $h:~X \to~\R$ if there exists a bounded measurable function $q: X \to \R$ such that 
\[
g(x) - h(x) = \lim_{t \to 0}\frac{q(\psi_t(x)) - q(x)}{t} \quad \textrm{for every $x \in X$}.
\]
Moreover, we say that $g$ is \emph{$T$-cohomologous} to $h$ if there exists a bounded measurable function $r: X \to \R$ such that
\[
g(x) - h(x) = r(T(x)) - r(x) \quad \textrm{for every $x \in X$}. 
\]
We also say that a function is a \emph{$\Psi$-coboundary} (\emph{$T$-coboundary}) when it is $\Psi$-cohomologous to the zero function ($T$-cohomologous to the zero function). The functions $q, r: X \to \R$ are usually called \emph{transfer functions}.
\end{definition}

The following result gives some directions for \textbf{Question D}: 

\begin{proposition}\label{COH}
Let $\Phi = (\phi_t)_{t \in \R}$ be a $C^{1}$ flow on a Riemmanian manifold $X$ and $\widetilde{b}:~X \to~\R$ a $C^{1}$ $\phi_1$-coboundary function admitting a $C^1$ transfer function. Then, there exists a continuous function $b: X \to \R$ such that 
\[
\int_{0}^{1}(b \circ \phi_s)ds = \widetilde{b}.
\] 
Moreover, the function $b$ is a $\Phi$-coboundary also admitting a $C^{1}$ transfer function. 
\end{proposition}

\begin{proof}
Suppose that  
\[
\widetilde{b} = g \circ \phi_1 - g \quad \textrm{for some $C^{1}$ function $g: X \to \R$.}
\]

Defining a function $b: X \to \R$ as 
\[
b(x) = d_x g(x) \bigg(\frac{d}{ds}\phi_s(x)|_{s = 0}\bigg) = \lim_{s \to 0}\frac{g(\phi_s(x)) - g(x)}{s},
\]
we can check that $b$ is continuous and
\[
\int_{0}^{t}(b \circ \phi_s)ds = g \circ \phi_t - g \quad \textrm{for every $t \ge 0$}.
\]
In particular, we have
\[
\int_{0}^{1}(b \circ \phi_s)ds = g \circ \phi_1 - g = \widetilde{b}.
\]

Notice that $\widetilde{b}$ is a $\phi_1$-coboundary and $b$ is a $\Phi$-coboundary, both admitting the same $C^{1}$ transfer function $g$.
\end{proof}
\begin{remark}
In Proposition \ref{COH}, if $X$ is a compact manifold, $\phi_1$ is a $C^{\infty}$ diffeomorphism and $\widetilde{b}$ is $C^{\infty}$ then one can always guarantee the existence of a $C^{\infty}$ function $g$ such that $\widetilde{b} = g\circ \phi_1 - g$ (see \cite{LMM86}). 
\end{remark}

A map $T$ is said to be \emph{uniquely ergodic} if it admits a unique $T$-invariant measure. In the same manner, we say that a flow $\Psi$ is \emph{uniquely ergodic} when there exists a unique $\Psi$-invariant measure. 
\begin{proposition}\label{UER}
Let $\Phi = (\phi_t)_{t \in \R}$ be a continuous flow on a compact metric space $M$ and such that $\phi_1: M \to M$ is uniquely ergodic. Then
\begin{enumerate}
\item $(a_t)_{t \ge 0}$ is an asymptotically additive family with respect to $\Phi$ if and only if $(a_n)_{n \ge 1}$ is an asymptotically additive sequence with respect to $\phi_1$;

\item for each asymptotically additive family of continuous functions $a = (a_t)_{t \ge 0}$ there exists a continuos function $c: M \to \R$ such that
\[
\lim_{t \to +\infty}\frac{1}{t}\bigg\|a_t - \int_{0}^{t}(c \circ \phi_s)ds\bigg\|_{\infty} = 0.
\] 
\end{enumerate}
\end{proposition}
\begin{proof}
Consider the linear operator $\mathcal{L}: C(X) \to C(X)$ given by 
\[
\mathcal{L}(a)(x) = \int_{0}^{1}(a\circ \phi_s)(x)ds.
\]
Then, one can see that $\|\mathcal{L}\| \le 1$, where we are considering the norm 
\[
\|\mathcal{L}\|:= \sup_{f \in C(X), \|f\|_{\infty} \neq 0} \frac{\|\mathcal{L}f\|_{\infty}}{\|f\|_{\infty}}.
\]

Since $\mathcal{L}$ is a bounded linear operator, its spectrum $\sigma(\mathcal{L})$ is compact and is such that
\[
\sigma(\mathcal{L}) \subset [-\|\mathcal{L}\|, \|\mathcal{L}\|] \subset [-1, 1].
\]
In particular, this implies that there exists $\lambda > 1$ such that the linear operator $\mathcal{L} - \lambda I: C(X) \to C(X)$ is a bijection. This readily implies that given a continuous function $b: X \to \R$ there exists a unique continuous function $a: X \to \R$ such that 
\begin{equation}\label{spec}
\int_{0}^{1}(a \circ \phi_s)ds - \lambda a = b.
\end{equation}

Now let $\nu$ be the unique $\phi_1$-invariant measure and consider the function
\[
c := a - \lambda\int_{M}a d\nu. 
\]
It follows from \eqref{spec} that 
\begin{equation}\label{EQ}
\bigg\|a_n - \int_{0}^{n}(c\circ \phi_s)ds\bigg\|_{\infty} \le \bigg\| a_n - \sum_{k=0}^{n-1}b \circ \phi_1^{k}\bigg\|_{\infty} + \bigg\|n\lambda\bigg(\int_{M}a d\nu - \sum_{k=0}^{n-1}a\circ \phi_1^{k}\bigg)\bigg\|_{\infty}.
\end{equation}
Since $\phi_1$ is uniquely ergodic, we have in particular that 
\[
\lim_{n \to \infty}\bigg\|\lambda\bigg(\int_{M}ad\nu - \frac{1}{n}\sum_{k=0}^{n-1}a\circ \phi_1^{k}\bigg)\bigg\|_{\infty} = 0
\]
for every $a \in C(X)$. 

Let us start proving item 1. By the arbitrariness of $b \in C(X)$, if $(a_n)_{n \ge 1}$ is asymptotically additive with respect to $\phi_1$ then it follows directly from \eqref{EQ} that $(a_t)_{t \ge 0}$ is asymptotically additive with respect to $\Phi$. The converse is immediate.

In order to prove item 2, notice that Theorem 1.2 in \cite{Cun20} guarantees the existence of a function $b \in C(X)$ such that 
\[
\lim_{n \to \infty}\frac{1}{n}\bigg\|a_n - \sum_{k=0}^{n-1} b\circ \phi_1^{k}\bigg\|_{\infty} = 0. 
\]
Therefore, it readily follows from \eqref{EQ} that there exists a function $c \in C(X)$ such that 
\[
\lim_{t \to \infty} \frac{1}{t}\bigg\|a_t - \int_{0}^{t}(c \circ \phi_s)ds\bigg\|_{\infty} = 0,
\]
as desired. 
\end{proof}

\begin{remark}
Notice that if the time-one map $\phi_1$ of a flow $\Phi$ is uniquely ergodic then the flow $\Phi$ itself is uniquely ergodic. From this, we can see that the setup of Proposition \ref{UER} is quite restrictive. We still don't know if the result remains true in full generality, without the unique ergodicity hypotheses.   
\end{remark}

\subsection{Some applications and consequences of Theorem \ref{main}}\label{APP}

In this section, using Theorem \ref{main}, we are going to show how to extend some results of the nonadditive thermodynamic formalism and multifractal analysis for flows. 

\subsubsection{Nonadditive thermodynamic formalism}

Let $M$ be a compact metric space and $\Lambda \subset M$ be a locally maximal hyperbolic set for a topologically mixing $C^{1}$ flow $\Phi$. Suppose that $a = (a_t)_{t \ge 0}$ is an asymptotically additive family of continiuous functions with tempered variation. Hence, Theorem \ref{main} guarantees the existence of a continuous function $b: \Lambda \to \R$ such that 
\[
\sup_{\mu \in \cM(\Phi)}\biggl\{h_{\mu}(\Phi) + \lim_{t \to \infty}\frac{1}{t}\int_{\Lambda}a_t d\mu\biggr\} = \sup_{\mu \in \cM(\Phi)}\biggl\{h_{\mu}(\Phi) + \int_{\Lambda}b d\mu\biggr\}.
\]
Moreover, by the definition of the nonadditive topological pressure introduced in \cite{BH20}, $P_{\Phi}(a) = P_{\Phi}\big((S_tb)_{t \ge 0}\big)= P^{classic}_{\Phi}(b)$, where the last one is the classical topological pressure for the function $b$. Therefore, the classical variational principle for continuous flows implies that
\[
P_{\Phi}(a) = \sup_{\mu \in \cM(\Phi)}\biggl\{h_{\mu}(\Phi) + \lim_{t \to \infty}\frac{1}{t}\int_{\Lambda}a_t d\mu\biggr\}.
\]
This is a variational principle for asymptotically additive families of continuous functions with respect to a hyperbolic flow. Notice that by Corollary \ref{EEX} the variational principle is also valid for expansive continuous flows. Therefore, in the case of locally maximal hyperbolic sets for flows, some suspension flows over subshifts of finite type and expansive flows admitting an invariant measure with full support, this result extends Theorem 9 in \cite{BH20} and Theorem 1.1 in \cite{BH21a}. 

Based on this variational principle, we can also define the notion of equilibrium measures for asymptotically additive families of potentials. We say that $\nu \in \cM(\Phi)$ is an equilibrium measure for $a$ with respect to $\Phi$ if 
\[
P_{\Phi}(a) = h_{\nu}(\Phi) + \lim_{t \to +\infty}\frac{1}{t}\int_{\Lambda}a_t d\nu.
\]

The following is a immediate consequence of the physical equivalence results for almost and asymptotically additive families. 
\begin{corollary}
Let $\Lambda \subset M$ be a locally maximal hyperbolic set for a $C^{1}$ flow $\Phi$ or suppose $\Phi$ is an expansive continuous flow admitting an invariant measure with full support. If $a$ is an asymptotically additive family of functions with respect to $\Phi$, then 
\begin{enumerate}
\item the set of equilibrium measures for $a$ is a non-empty compact and convex set;
\item each extreme point of the set of equilibrium measures is an ergodic measure.
\end{enumerate}
\end{corollary}
\begin{proof}
By Theorem \ref{main} and Theorem \ref{EXP}, respectively for hyperbolic flows and expansive flows, the space of equilibrium measures for $a$ is the same as the space of equilibrium measures for some continuous function.	
\end{proof}

\begin{remark}
Theorem 3.5 in \cite{BH21a} guarantees that for locally maximal hyperbolic sets for $C^{1}$ flows, every almost additive family with bounded variation admits a unique equilibrium measure. Since we don't know if \textbf{Question B} is true in general, we are not able to obtain this uniqueness result directly from the additive classical case originally proved in \cite{Fra77}. We will address some regularity issues of the physical equivalence relations in the second part of this work (see the subsections \ref{BHR} and \ref{BHB}). 
\end{remark}

\subsubsection{Multifractal analysis: beyond families with unique equilibrium measures}

Theorem 9 in \cite{BH21a} establishes a conditional variational principle for almost additive families of potentials with uniqueness of equilibrium measures. In this section, using the work of Climenhaga in \cite{Cli13}, Cuneo in \cite{Cun20}, we are going to show how to obtain a conditional variational principle including more classes of almost additive families of functions, without using the uniqueness of equilibrium assumption. Moreover, as a direct consequence of Theorem \ref{main}, we also can extend the conditional variational principle to include asymptotically additive families of continuous potentials. 

Let $M$ be a compact metric space and $\Lambda \subset M$ be a locally maximal hyperbolic set for a $C^1$ topologically mixing flow $\Phi = (\phi_t)_{t \in \R}$. Given two asymptotically additive families of continuous functions $a = (a_t)_{t \ge 0}$ and $b = (b_t)_{t \ge 0}$, we consider the level sets
\[
K_{\alpha}^{\Phi}(a,b):= \biggl\{x \in \Lambda: \lim_{t \to +\infty}\frac{a_t(x)}{b_t(x)} = \alpha\biggr\}, \quad \alpha \in \R.
\]

In this section, we consider $\cM_{erg}(\Phi)$ as the set of ergodic $\Phi$-invariant measures. We also use this same notation for maps.  

Let $\cA(M)$ be the set of all almost additive families of continuous functions and, respectively, $\cA\cA(M)$ the set of asymptotically additive families of continuous functions $a=(a_t)_{t\ge0}$ on $M$ with tempered variation such that
\[
\sup_{t \in [0,s]}\Vert a_t\Vert_{\infty} < \infty\quad\text{for some} \ s > 0,
\]
and $\cE(M) \subset \cA(M)$ the subset of families having a unique equilibrium measure (the existence of almost additive families with unique equilibrium measures is guaranteed by Theorem 3.5 in \cite{BH21a}). 

Now let us define the notion of \emph{$u$-dimension} for flows which was originally introduced in \cite{BS00}. Let $\Phi = (\phi_t)_{t \in \R}$ be a continuous flow on a compact metric space~$X$. Given a positive continuous function $u\colon X \to \R$, we consider the additive family of continuous functions $(u_t)_{t\ge0}$ defined by
\[
u_t(x) = \int_{0}^{t}(u \circ \phi_s)(x)ds
\]
for every $x \in X$ and $t > 0$.
For each $Z \subset X$ and $\alpha \in \R$, let
\[
N(Z,u,\alpha,\epsilon) = \lim_{T \to \infty}\inf_{\Gamma}\sum_{(x,t) \in \Gamma}e^{-\alpha u(x,t,\epsilon)},
\]
with the infimum taken over all countable sets
$\Gamma\subset X \times[T,+\infty)$ covering~$Z$.
Finally, we define
\[
\dim_{u,\epsilon}Z = \inf\bigl\{\alpha \in \R: N(Z,u,\alpha,\epsilon) = 0\bigr\}.
\]
The limit
\[
\dim_{u}Z := \lim_{\epsilon \to 0}\dim_{u,\epsilon}Z
\]
exists and is called the $u$-dimension of the set $Z$ with respect to the flow~$\Phi$. When $u \equiv 1$ the number $\dim_{u}Z$ coincides with the topological entropy $h(\Phi|_Z)$ of $\Phi$ on the set~$Z$.

\begin{theorem}\label{BFU}
Let $\Lambda \subset M$ be a locally maximal hyperbolic set for a $C^1$ topologically mixing flow $\Phi = (\phi_t)_{t \in \R}$ such that $h(\Phi) < \infty$. Let $a = (a_t)_{t \ge 0}$ and $b = (b_t)_{t \ge 0}$ be  asymptotically additive families of continuous functions in $\cA\cA(\Lambda)$ and such that 
\[
\lim_{t \to \infty}\frac{1}{t}\int_{\Lambda}b_t d\mu \ge 0 \quad \textrm{for all $\mu \in \cM(\Phi)$} 
\]
with equality only permitted when 
\[
\lim_{t \to \infty}\frac{1}{t}\int_{\Lambda}a_t d\mu \neq 0.
\]
If $\alpha \notin \cP(\cM(\Phi))$, then $K_{\alpha}^{\Phi}(a,b) =~\emptyset$. Moreover, if $\alpha \in \interior \cP(\cM(\Phi))$ then $K_{\alpha}^{\Phi}(a,b) \neq \emptyset$, and the following properties hold:
\begin{enumerate}
\item the level sets $K^{\Phi}_{\alpha}(a,b)$ satisfy the conditional variational principle
\[
\dim_{u}K_{\alpha}^{\Phi}(a,b) = \sup \biggl\{\frac{h_{\mu}(\Phi)}{\int_{\Lambda}u d\mu}: \mu \in \cM_{erg}(\Phi) \textrm{ and } \cP(\mu) = \alpha \biggr\};
\]
		
\item $\dim_{u}K_{\alpha}^{\Phi}(a,b) = \inf\{T_u(q): q \in \R\}$, where $T_u(q)$ is defined by $T_u(q) = \inf\{t \in \R: P(q(a-\alpha b) - tu) \le 0 \}$;
		
\item for each $\epsilon > 0$ there exists an ergodic measure $\mu_\alpha$ supported on $K_{\alpha}^{\Phi}(a,b)$ such that 
\[
\bigg|\dim_u \mu_{\alpha} = \frac{h_{\mu_{\alpha}}(\Phi)}{\int_{\Lambda}u d\mu_\alpha} - \dim_{u}K_{\alpha}^{\Phi}(a,b) \bigg| < \epsilon
\]
		
\item the spectrum $\alpha \mapsto \dim_{u}K_{\alpha}^{\Phi}(a,b)$ is continuous on $\interior \cP(\cM(\Phi))$.
\end{enumerate}
\end{theorem}

Notice that here we no longer require that $\Span\{a,b,\overline{u}\} \subset \cE(X)$ as in \cite{BH21b}. Moreover, this result extends Theorem 9 in \cite{BH21b} to asymptotically additive families. 

The proof will be divided in several steps.  

First, let us recall some definitions and introduce some nomenclature used in this section. Let $X$ be a compact metric space, $T: X \to X$ a continuous map. We write $\cM(T)$ to denote the space of $T$-invariant measures, $\cM_{erg}(T)$ the space of $T$-invariant ergodic measures, and $C(X)$ the space of continuous real valued functions defined on $X$. Let $\varphi, \psi \in C(X)$, where $\int_{X}\psi d\mu \ge 0$ for all $\mu \in \cM(T)$ with equality only permitted if $\int_{X}\varphi d\mu \neq 0$. We consider the map $\cP: \cM(T) \to \R$ given by 
\[
\cP_T(\mu) := \frac{\int_{X}\varphi d\mu}{\int_{X}\psi d\mu}.
\] 
We also define the level sets for $\varphi, \psi$ by 
\[
K_{\alpha}^{T}(\varphi,\psi):= \biggl\{x \in X: \lim_{n \to \infty}\frac{S_n\varphi(x)}{S_n\psi(x)} = \alpha\biggr\},
\]
where $S_nh = \sum_{k=0}^{n-1}h \circ T^{k}$ for every function $h: X \to \R$.
\begin{theorem}[{\cite[Theorem~3.3]{Cli13}}]\label{CLI}
Let $X$ be a compact metric space and $T: X \to X$ be a continuous map such that the map $\mu \mapsto h_{\mu}(T)$ is upper semicontinuous and $h(T) < \infty$. Suppose that there is a dense subspace $D \subset C(X)$ such that every $\varphi \in D$ has a unique equilibrium measure. Let $\varphi, \psi \in C(X)$ be such that $\int_{X}\psi d\mu \ge 0$ for all $\mu \in \cM(T)$ with equality only permitted if $\int_{X}\varphi d\mu \neq 0$. Then $K^{T}_{\alpha}(\varphi,\psi) = \emptyset$ for every $\alpha \notin \cP_T(\cM(T))$, while for every $\alpha \in \interior \cP_T(\cM(T))$ we have that
\begin{enumerate}
		
\item the level sets $K^{T}_{\alpha}(\varphi,\psi)$ satisfy the conditional variational principle
\[
\dim_u K^{T}_{\alpha}(\varphi,\psi) = \sup\biggl\{\frac{h_{\mu}(T)}{\int_{X}ud\mu}: \mu \in \cM(T) \quad \textrm{and} \quad \cP_T(\mu) = \alpha \biggr\};
\]
		
\item $\dim_u K^{T}_{\alpha}(\varphi,\psi) = \inf\{T_u(q):q \in \R\}$, where $T_u(q)$ is defined by $T_u(q) = \inf\{t \in \R: P_{classic}(q(\varphi - \alpha\psi) -tu) \le 0\}$;
		
\item For each $\epsilon > 0$ there exists an ergodic measure $\mu$ supported on $K_{\alpha}^{T}(\varphi,\psi))$ such that 
\[
\bigg|\frac{h_{\mu}(T)}{\int_{X}ud\mu} - \dim_u K^{T}_{\alpha}(\varphi,\psi)\bigg| < \epsilon.
\]
\end{enumerate}
\end{theorem}

Moreover, by Proposition 2.14 in \cite{Cli13}, one can see that the map
\[
\alpha \to \dim_u K^{T}_{\alpha}(\varphi,\psi)
\]
is continuous on $\interior \cP_{T}(\cM(T))$.

We recall that a sequence of functions $\cF = (f_n)_{n \ge 1}$ is asymptotically additive (with respect to $T$) if for each $\epsilon > 0$ there exists a function $f$ such that 
\[
\limsup_{n \to \infty}\frac{1}{n}\|f_n - S_nf\|_{\infty} < \epsilon.
\]

We also recall that a sequence $\cF = (f_n)_{n \ge 1}$ is almost additive (with respect to $T$) if there exists $C > 0$ such that 
\[
-C + f_{m}(x) + f_{n}(T^m(x)) \le f_{m+n}(x) \le f_{m}(x) + f_{n}(T^m(x)) + C
\]
for every $x \in X$ and all $m,n \ge 1$. 

Let $\cF = (f_n)_{n \ge 0}$ and $\cG = (g_n)_{n \ge 0}$ be two almost additive sequences of continuous functions where $\lim_{n \to \infty}\frac{1}{n}\int_{X}g_nd\mu \ge 0$ for all $\mu \in \cM(T)$ with equality only permitted if $\lim_{n \to \infty}\frac{1}{n}\int_{X}f_nd\mu \neq 0$.

Consider level sets
\[
K_{\alpha}^{T}(\cF,\cG):= \biggl\{x \in X: \lim_{n \to \infty}\frac{f_n(x)}{g_n(x)} = \alpha\biggr\},
\]
and also the map $\cQ_T: \cM(T) \to \R$ defined by 
\[
\cQ_T(\mu) = \lim_{n \to \infty}\frac{\int_{X}f_n d\mu}{\int_{X}g_n d\mu}.
\]

\begin{corollary}\label{NCVC}
Let $X$ be a compact metric space and $T: X \to X$ a map satisfying the conditions in Theorem \ref{CLI}. Then $K^{T}_{\alpha}(\cF,\cG) = \emptyset$ for every $\alpha \notin \cQ_T(\cM(T))$, while for every $\alpha \in \interior \cQ_T(\cM(T))$ we have that
\begin{enumerate}
		
\item the level sets $K^{T}_{\alpha}(\cF,\cG)$ satisfy the conditional variational principle
\[
\dim_u K^{T}_{\alpha}(\cF,\cG) = \sup\biggl\{\frac{h_{\mu}(T)}{\int_{X}ud\mu}: \mu \in \cM(T) \quad \textrm{and} \quad \cQ_T(\mu) = \alpha \biggr\};
\]
		
\item $\dim_u K^{T}_{\alpha}(\cF,\cG) = \inf\{T_u(q):q \in \R\}$, where $T_u(q)$ is defined by $T_u(q) = \inf\{t \in \R: P(q(\cF - \alpha\cG) -tu) \le 0\}$;
		
\item For each $\epsilon > 0$ there exists an ergodic measure $\mu$ supported on $K_{\alpha}^{T}(\cF,\cG)$ such that 
\[
\bigg|\frac{h_{\mu}(T)}{\int_{X}ud\mu} - \dim_u K^{T}_{\alpha}(\cF,\cG)\bigg| < \epsilon.
\]
		
\item the function $\alpha \to \dim_u K^{T}_{\alpha}(\cF,\cG)$ is continuous on $\interior \cQ_T(\cM(T))$.
\end{enumerate}
\end{corollary}

\begin{proof}
By Theorem 1.2 in \cite{Cun20} there exist $f, g \in C(X)$ such that 
\[
\lim_{n \to \infty}\frac{1}{n}\|f_n - S_nf\|_{\infty} = 0 \quad \textrm{and} \quad \lim_{n \to \infty}\frac{1}{n}\|g_n - S_ng\|_{\infty} = 0.
\]	
By the variational principle for almost additive sequence of continuous functions, one can see that  \[P(q(\cF-\alpha\cG)-su) = P_{classic}(q(f-\alpha g)-su)
\]
for every $q, \alpha$ and $s \in \R$. Moreover, $\cP_T(\mu) = \cQ_T(\mu)$ for all $\mu \in \cM(T)$, and $K_{\alpha}^T(\cF,\cG) = K_{\alpha}^T(f,g)$ for each $\alpha \in \R$. Hence, the result now follows directly from Theorem \ref{CLI}.  
\end{proof}

We say that a continuous map $T:X \to X$ on a compact metric space $X$ (or a continuous flow $\Phi$ on $X$) has \emph{entropy density of ergodic measures} if for every invariant measure $\mu$ there exist ergodic measures $\nu_n$ for $n\in\N$ such that $\nu_n \to \mu$ in the weak$^*$ topology and $h_{\nu_n}(T) \to h_\mu(T)$ (or $h_{\nu_n}(\Phi) \to h_{\mu}(\Phi)$) when $n\to\infty$.

In order to give some examples having entropy density of ergodic measures, we will first recall a few notions. 
Given $\delta>0$, we say that $T$~has \emph{weak specification at scale~$\delta$} if there exists $\gamma \in \N$ such that for every $(x_1,n_1),\ldots,(x_k,n_k) \in X \times \N$ there are $y \in X$ and times $\gamma_1,\ldots, \gamma_{k-1} \in \N$ such that $\gamma_i \le \gamma$ and
\[
d_{n_i}(T^{s_{i-1} + \gamma_{i-1}}(y),x_i) < \delta \quad \text{for} \ i=1,\ldots, k,
\]
where $s_i = \sum_{i=1}^{i}n_i + \sum_{i=1}^{i-1}\gamma_i$ with $n_0 = \gamma_0 = 0$.
When one can take $\gamma_i = \gamma$ for $i =1,\ldots,k-1$, we say that $T$~has \emph{specification at scale~$\delta$}. Finally, we say that $T$~has \emph{weak specification} if it has weak specification at every scale~$\delta$ and, analogously, we say that $T$~has \emph{specification} if it has specification at every scale~$\delta$. 

It was shown in \cite{EKW94} and \cite{PS05} that mixing subshifts of finite type and mixing locally maximal hyperbolic sets have entropy density of ergodic measures. More recently, it was shown in~\cite{CLT20} that a continuous map $T\colon X \to X$ on a compact metric space with the weak specification property such that the entropy map $\mu \mapsto h_\mu(T)$ is upper semicontinuous, has entropy density of ergodic measures. Some examples include expansive maps with specification or with weak specification, topologically transitive locally maximal hyperbolic sets for diffeomorphisms, and transitive topological Markov chains.

We can also introduce the same notions for flows. Let $\Phi = (\phi_t)_{t \in \R}$ be a continuous flow on a compact metric space $X$. We say that $\Phi$ has \emph{weak specification at scale} $\delta > 0$ if there exists $\gamma > 0$ such that for every finite collection of orbit segments $\{(x_i,t_i)\}_{i=1}^{k}$, there exists a point $y \in X$ and a sequence of \emph{transition times} $\gamma_1, ..., \gamma_{k-1} \in [0, \gamma]$ such that 
\[
d_{t_j}(\phi_{s_{j-1} + \gamma_{j-1}}(y),x_j) < \delta \quad \textrm{for} \quad j= 1,...,k,
\]
where $s_j = \sum_{i=1}^{j}t_i + \sum^{j-1}_{i=1}\gamma_{i}$ and $s_0 = \gamma_{0} = 0$. We say that $\Phi$ has \emph{weak specification} if it has weak specification at every scale $\delta$. When, for every scale $\delta > 0$, we can take the approximating orbit $y \in X$ to be periodic, and the transition times $\gamma_{i}$ close to $\gamma$, we say that $\Phi$ has \emph{specification}. For a proper definition, see for example \cite{Bow72}. It was also proved in \cite{CLT20} that every expansive flow with the weak specification property has entropy density of ergodic measures. In particular, locally maximal hyperbolic sets for $C^1$ topologically mixing flows are expansive and have the specification property, that is, they have entropy density of ergodic measures.

\emph{Proof of Theorem \ref{BFU}.} We are going to use the Markov partitions introduced in subsection \ref{subsec32}. Let us start proving the result for almost additive families $a$ and $b$. 

Let's consider the following level sets
\[
K^{\Phi}_\alpha(a,b) := \bigg\{x \in \Lambda: \lim_{t \to \infty}\frac{a_t(x)}{b_t(x)} = \alpha \bigg\}, \quad K^{T_Z}_\alpha(c,d) := \bigg\{x \in Z: \lim_{n \to \infty}\frac{c_n(x)}{d_n(x)} = \alpha \bigg\},
\]
where $c = (c_n)_{n \in \N}$ and $d = (d_n)_{n \in \N}$ are almost additive sequences of continuous functions given by $c_n(x) = a_{\tau_{n}(x)}(x)$ and $d_n(x) = b_{\tau_{n}(x)}(x)$ (see Lemma 3.1 in \cite{BH21a}). By Lemma 3.4 in \cite{BH21a}, for every $\Phi$-invariant measure $\mu$ induced by an ergodic $T_Z$-invariant measure $\nu$, we have that 
\begin{equation}\label{PQT}
\cP(\mu):= \lim_{t \to +\infty}\frac{\int_{X}a_t d\mu}{\int_{X}b_td\mu} = \lim_{n \to +\infty}\frac{\int_{X}c_n d\nu}{\int_{X}b_nd\nu}:= \cQ_{T_Z}(\nu).
\end{equation}

We also recall that a $T_Z$-invariant measure $\nu$ is ergodic if and only if the induced $\Phi$-invariant measure $\mu$ is ergodic (see identity \eqref{eq13}).  

It follows from Proposition 8, Lemma 14 and Lemma 15 in \cite{BH21b} that 
\begin{equation}\label{EQV}
\begin{split}
\dim_{u}K^{\Phi}_\alpha(a,b) &= \dim_{u}\{\phi_s(x) \in \Lambda : x \in K_{\alpha}(c,d) \textrm{ and } s \in [0,\tau(x)]\} \\
& = \inf\{\beta \in \R: N_{\beta}(K_{\alpha}(c,d)) = 0\} = \dim_{I_u}K^{T_Z}_\alpha(c,d).
\end{split}
\end{equation}
Now we are going to check the conditions to apply Corollary \ref{NCVC} for the map $T_Z$ and the sequences $c$ and $d$. In fact, since $\Phi$ is expansive, one can verify that $T_Z$ is also expansive, which implies that $\nu \mapsto h_{\nu}(T_Z)$ is upper semicontinuous. Moreover, by Abramov's entropy formula (identity \eqref{eq14}), we have 
\[
h_{\nu}(T_Z) \le h_{\mu}(\Phi)\sup \tau \le h(\Phi)\sup \tau < \infty
\]  
for every $\nu \in \cM(T_Z)$, and then $h(T_Z) < \infty$. Letting $D$ be the space of H\"older continuous functions, one can see that $D$ is dense in the space of continuous functions $\varphi: Z \to \R$, and every H\"older continuous function has a unique equilibrium measure with respect to the map $T_Z$. 

By hypothesis, $\lim_{t \to \infty}\frac{1}{t}\int_{\Lambda}b_t d\mu \ge 0$ for all $\mu \in \cM(\Phi)$ with equality only permitted when $\lim_{t \to \infty}\frac{1}{t}\int_{\Lambda}a_t d\mu \neq 0$. It follows again from Lemma 3.4 in \cite{BH21a} that 
\[
\frac{\lim_{n \to \infty}\frac{1}{n}\int_{Z}c_n d\nu}{\int_{Z}\tau d\nu} = \lim_{t \to \infty}\frac{1}{t}\int_{\Lambda}a_t d\mu \quad \textrm{and}
\]
\[
\frac{\lim_{n \to \infty}\frac{1}{n}\int_{Z}d_n d\nu}{\int_{Z}\tau d\nu} = \lim_{t \to \infty}\frac{1}{t}\int_{\Lambda}b_t d\mu.
\]
for every $\Phi$-invariant ergodic measure $\mu$ induced by an ergodic $T_Z$-invariant measure $\nu$. Since $\int_{Z}\tau d\nu >0$ for all $\nu \in \cM(T_Z)$, we obtain that $\lim_{n \to \infty}\frac{1}{n}\int_{Z}d_n d\nu \ge 0$ for all $\nu \in \cM_{erg}(T_Z)$ and, in particular, that
\begin{equation}\label{ERG}
\lim_{n \to \infty}\frac{1}{n}\int_{Z}d_n d\nu = 0 \quad \textrm{implies} \quad \lim_{n \to \infty}\frac{1}{n}\int_{Z}c_n d\nu \neq 0 \quad \textrm{for all  $\nu \in \cM_{erg}(T_Z)$.}
\end{equation}
Now observe that $T_Z: Z \to Z$ is topologically mixing and is conjugated to a topological Markov chain $\sigma: \Sigma_A \to \Sigma_A$ (see subsection \ref{subsec32}). From this, one can see that $T_Z$ has entropy density of ergodic measures which, in particular, implies that $\cM_{erg}(T_Z)$ is dense in $\cM(T_Z)$. Let $\eta \in \cM(T_Z)$. By the density of $\cM_{erg}(T_Z)$ in $\cM(T_Z)$, given $\epsilon > 0$ there exists $\nu \in \cM_{erg}(T_Z)$ such that
\[
\int_{Z}\frac{1}{n}d_n d\eta > \int_{Z}\frac{1}{n}d_n d\nu - \epsilon \quad \textrm{for every $n \ge 1$.}
\]  
Then
\[
\lim_{n \to \infty}\frac{1}{n}\int_{Z}d_n d\eta \ge \lim_{n \to \infty}\frac{1}{n}\int_{Z}d_n d\nu - \epsilon \ge -\epsilon.
\]
Since the measure $\eta$ and $\epsilon> 0$ are arbitrary, we conclude that $\lim_{n \to \infty}\frac{1}{n}\int_{Z}d_n d\nu \ge 0$ for every $\nu \in \cM(T_Z)$. 

Now fix any measure $\eta \in \cM(T)$. Given $\epsilon > 0$, the entropy density of ergodic measures implies the existence of $\nu \in \cM_{erg}(T_Z)$ such that 
\begin{equation}\label{TIV}
\bigg|\lim_{n \to \infty}\frac{1}{n}\int_{Z}d_n d\eta - \lim_{n \to \infty}\frac{1}{n}\int_{Z}d_n d\nu\bigg| \le \epsilon \quad \textrm{and} \quad \bigg|\lim_{n \to \infty}\frac{1}{n}\int_{Z}c_n d\eta - \lim_{n \to \infty}\frac{1}{n}\int_{Z}c_n d\nu\bigg| \le \epsilon.
\end{equation}
Now suppose that 
\[
\lim_{n \to \infty}\frac{1}{n}\int_{Z}d_n d\eta = 0 \quad \textrm{and} \quad \lim_{n \to \infty}\frac{1}{n}\int_{Z}c_n d\eta = 0.
\]
Then it follows from \eqref{TIV} that 
\[
\bigg|\lim_{n \to \infty}\frac{1}{n}\int_{Z}d_n d\nu\bigg| \le \epsilon \quad \textrm{and} \quad \bigg|\lim_{n \to \infty}\frac{1}{n}\int_{Z}c_n d\nu\bigg| \le \epsilon.
\]
By the arbitrariness of $\epsilon$, this contradicts \eqref{ERG}. Hence, since $\eta \in \cM(T_Z)$ is also arbitrary, we conclude in particular that
\[
\lim_{n \to \infty}\frac{1}{n}\int_{Z}d_n d\eta = 0 \quad \textrm{implies} \quad \lim_{n \to \infty}\frac{1}{n}\int_{Z}c_n d\eta \neq 0 \quad \textrm{for all  $\eta \in \cM(T_Z)$,}
\]
as desired. Now we finally are in the conditions of applying Corollary \ref{NCVC}.  

Let $\mu \in \cM(T_Z)$ with $\cQ_{T_Z}(\mu) = \alpha$. By applying item (1) of Corollary \ref{NCVC} to the map $T_Z: Z \to Z$ and the sequences $c$ and $d$, we have 
\[
\dim_{I_u} K_{\alpha}^{T_Z}(c,d) \ge \frac{h_{\mu}(T_Z)}{\int_{Z}I_u d\mu}.
\]

Given any $\epsilon > 0$, by item (3) of Corollary \ref{NCVC}, there exists $\nu \in \cM_{erg}(T_Z)$ with $\cQ_{T_Z}(\nu) = \alpha$ such that 
\[
\frac{h_{\nu}(T_Z)}{\int_{Z}I_u d\nu} > \dim_{I_u} K_{\alpha}^{T_Z}(c,d) - \epsilon \ge \frac{h_{\mu}(T_Z)}{\int_{Z}I_u d\mu} - \epsilon.
\]
Since the measure $\mu$ and $\epsilon > 0$ are arbitrary, we obtain that

\[
\begin{split}
\sup\biggl\{\frac{h_{\nu}(T_Z)}{\int_{X}I_ud\nu}: \nu \in \cM_{erg}(T_Z) \quad \textrm{and} \quad \cQ_{T_Z}(\nu) = \alpha \biggr\} \ge \dim_{I_u} K_{\alpha}^{T_Z}(c,d). 
\end{split}
\]
Then, from the fact that $\cM_{erg}(T_Z) \subset \cM(T_Z)$, we conclude that 
\begin{equation}\label{DIU}
\dim_{I_u} K_{\alpha}^{T_Z}(c,d) = \sup\biggl\{\frac{h_{\nu}(T_Z)}{\int_{X}I_ud\nu}: \nu \in \cM_{erg}(T_Z) \quad \textrm{and} \quad \cQ_{T_Z}(\nu) = \alpha \biggr\}.
\end{equation}

Now it follows from \eqref{eq13}, \eqref{eq14} and \eqref{PQT} that for each $\nu \in \cM_{erg}(T_Z)$ with $\cQ_{T_Z}(\nu) = \alpha$, the induced measure $\eta \in \cM_{erg}(\Phi)$ is such that $\cP(\eta) = \alpha$, and 
\begin{equation}\label{PTN}
\begin{split}
\frac{h_{\eta}(\Phi)}{\int_{\Lambda}u d\eta} = \bigg(\frac{h_{\nu}(T_Z)}{\int_{Z}\tau d\nu}\bigg)\bigg(\frac{1}{\int_{\Lambda}ud\eta}\bigg) = \bigg(\frac{h_{\nu}(T_Z)}{\int_{Z}\tau d\nu}\bigg)\bigg(\frac{\int_{Z}\tau d\nu}{\int_{Z}I_u d\nu}\bigg) = \frac{h_{\nu}(T_Z)}{\int_{Z}I_u d\nu}.
\end{split}
\end{equation}
Hence, it follows from \eqref{EQV} and \eqref{DIU} that 
\[
\dim_{u}K^{\Phi}_\alpha(a,b) = \sup\biggl\{\frac{h_{\mu}(\Phi)}{\int_{\Lambda}ud\mu}: \mu \in \cM_{erg}(\Phi) \quad \textrm{and} \quad \cP(\mu) = \alpha \biggr\},
\]
and this proves the item (1) of Theorem \ref{BFU}. 

Now let us prove the item (2). For each $\mu \in \cM_{erg}(\Phi)$ induced by $\nu \in \cM_{erg}(T_Z)$, we have
\[
\begin{split}
& h_{\mu}(\Phi) + \lim_{t \to \infty}\frac{1}{t}\int_{\Lambda}q(a_t - \alpha b_t)d\mu - s\int_{\Lambda}u d\mu \\
&= \frac{h_{\nu}(T_Z)}{\int_{Z}\tau d\nu} + \frac{\lim_{n \to \infty}\int_{Z}q(c_n-\alpha d_n)d\nu}{\int_{Z}\tau d\nu} - \frac{s\int_{Z}I_{u}d\nu}{\int_{Z}\tau d\nu}
\end{split}
\]
for every $\alpha, q, s \in \R$. Since $0<\inf \tau := \inf_{x \in \Lambda}\tau(x) \le \inf_{x \in Z}\tau(x)$, we have 
\[
\begin{split}
&h_{\mu}(\Phi) + \lim_{t \to \infty}\frac{1}{t}\int_{\Lambda}q(a_t - \alpha b_t)d\mu - s\int_{\Lambda}u d\mu \\
& \le h_{\nu}(T_Z) + \lim_{n \to \infty}\int_{Z}[q(c_n-\alpha d_n) - sI_u]d\nu\bigg(\frac{1}{\inf \tau}\bigg)\\
& \le P(q(c-\alpha d) - sI_u)\bigg(\frac{1}{\inf \tau}\bigg),
\end{split}
\]
which implies that 
\[
P(q(a -\alpha b) - su) \le P(q(c-\alpha d) - sI_u)\bigg(\frac{1}{\inf \tau}\bigg)
\]
for every $\alpha, q, s \in \R$.
Hence, $T_{u}(q) \le T_{I_u}(q)$ for every $q \in \R$. 

Since $\sup \tau := \sup_{x \in \Lambda} \tau(x) < \infty$, we can follow in an analogous way to obtain that 
\[
P(q(a -\alpha b) - su) \ge P(q(c-\alpha d) - sI_u)\bigg(\frac{1}{\sup \tau}\bigg)
\]
for every $\alpha, q, s \in \R$. Then, $T_{u}(q) \ge T_{I_u}(q)$ for every $q \in \R$. Since it follows from \eqref{EQV} that $\dim_u K_{\alpha}^{\Phi}(a,b) = \dim_{I_u}K_{\alpha}^{T_Z}(c,d)$, the item (2) of Corollary \ref{NCVC} gives that 
\[
\dim_u K_{\alpha}^{\Phi}(a,b) = \dim_{I_u}K_{\alpha}^{T_Z}(c,d) = \inf\{T_{I_u}(q): q \in \R\} = \inf\{T_{u}(q): q \in \R\}
\]
for each $\alpha \in \interior\cP(\cM(\Phi))$, as desired. 

In order to prove item (3), we just observe again that 
\[
\dim_u K_{\alpha}^{\Phi}(a,b) = \dim_{I_u}K_{\alpha}^{T_Z}(c,d)
\]
and that we can use \eqref{PTN} together with item (3) of Corollary \ref{NCVC}.  

Now let us show how to obtain the item (4) using entropy density of ergodic measures. Since $\Lambda$ is a locally maximal hyperbolic set for the $C^1$ topologically mixing flow $\Phi$, we have that $\Phi|_{\Lambda}$ has entropy density of ergodic measures. In particular, the set $\cM_{erg}(\Phi)$ is dense in $\cM(\Phi)$. Additionally, recall that $T_Z$ also has entropy density of ergodic measures.  

By \eqref{PQT} we already know that $\cP(\cM_{erg}(\Phi)) = \cQ_{T_Z}(\cM_{erg}(T_Z))$. Since $\overline{\cM_{erg}(\Phi)} = \cM(\Phi)$, $\overline{\cM_{erg}(T_Z)} = \cM(T_Z)$, and the maps $\mu \mapsto \cP(\mu)$, $\nu \mapsto \cQ_{T_Z}(\nu)$ are continuous, we have that 
\begin{equation}\label{MME}
\begin{split}
\cQ_{T_Z}(\cM(T_Z)) = \cQ_{T_Z}(\overline{\cM_{erg}(T_Z)}) &= \overline{\cQ_{T_Z}(\cM_{erg}(T_Z))} \\
& = \overline{\cP(\cM_{erg}(\Phi))} = \cP(\overline{\cM_{erg}(\Phi)}) = \cP(\cM(\Phi)). 
\end{split}
\end{equation}

It follows from item (4) of Corollary \ref{NCVC} that $\alpha \mapsto \dim_{I_u}K_{\alpha}^{T_Z}(c,d)$ is continuous on $\interior \cQ_{T_Z}(\cM(T_Z))$. Since $\dim_u K_{\alpha}^{\Phi}(a,b) = \dim_{I_u}K_{\alpha}^{T_Z}(c,d)$, by \eqref{MME} we conclude that the map $\alpha \mapsto \dim_u K_{\alpha}^{\Phi}(a,b)$ is also continuous on $\interior \cP(\cM(\Phi))$, and the Theorem \ref{BFU} is proved for almost additive families. Since, in particular, the result holds for the additive case, now we can use directly Corollary \ref{HYP} to complete the proof for asymptotically additive families. \qed

\begin{remark}
Notice that one could start proving Theorem \ref{BFU} for the additive case (without the hypothesis on density of ergodic measures) and after that, apply Theorem \ref{main} to obtain the full result for asymptotically and, consequently, almost additive families. The choice to start proving the result already in the almost additive case is due to some connections to the aforementioned results in the previous works \cite{BH21a} and \cite{BH21b}. 
\end{remark}

\begin{remark}
Observe that since we cannot guarantee the uniqueness of equilibrium measures for asymptotically additive families (even for hyperbolic flows), the extension of Climenhaga's results in \cite{Cli13} to include continuous potentials is crucial for our extension to asymptotically additive families. 
\end{remark}

\section{Part II: a nonadditive Liv\v{s}ic theorem and some regularity relations}

The second part of this work is inspired and based on the recent developments presented in \cite{HS24} for the case of discrete-time dynamical systems. The main focus here is to extend the notions and ideas to the case of flows, allowing us to study regularity equivalence issues for asymptotically and almost additive families with respect to suspension flows and, in particular, hyperbolic flows. 

\subsection{Cohomology for asymptotically additive families}\label{chm}

In this subsection, we are going to introduce some notions of cohomology for asymptotically additive families of continuous functions with respect to flows. The definitions are inspired by the physical equivalence results obtained in Part I (see Theorem \ref{main}, Corollary~\ref{corB}, Corollary \ref{HYP} and Theorem \ref{EXP}). We start recalling some basic concepts and tools in the additive (classical) setup.

We say that a function $\psi\colon X \to \mathbb{R}$ is \emph{Walters} (with respect to a flow $\Phi$) if for each $\epsilon > 0$ there exists a $\delta > 0$ such that for $x, y \in X$ and $t \ge 0$, we have that
\[
d_t(x,y) < \delta \textrm{ $\implies$ } \left|S_t\psi(x) - S_t\psi(y)\right| < \epsilon,
\]
where $d_t(x,y) = \max \{d(\phi_s(x), \phi_s(y)): s \in [0,t]\}$. In this case, we also say that the additive family $(S_t\psi)_{t \ge 0}$ satisfies the \emph{Walters property}. Moreover, we say that a function $\xi\colon X \to \mathbb{R}$ is \emph{Bowen} (with respect to $\Phi$) if there exist two numbers $L > 0$ and $\delta > 0$ such that for $x, y \in X$ and $t \ge 0$, we have that 
\[
d_t(x,y) < \delta \textrm{ $\implies$} |S_t\xi(x) - S_t\xi(y)| \le L.
\]
Clearly every Walters function is also Bowen. In the hyperbolic framework, the H\"older continuous functions are always Walters and, consequently, Bowen (see Proposition 7.3.1 in \cite{FH20}).  

A continuous flow $\Phi$ on a compact metric space $X$ is said to satisfy the \emph{Closing Lemma} if for every $\epsilon > 0$ there exists $\delta > 0$ such that if $x \in X$ and $t \ge 0$ satisfying $d(\phi_t(x),x)<\delta$, then there exists a periodic orbit $\{\phi_s(y): 0 \le s \le T \}$ with $|T - t| < \epsilon$ such that $d_t(x,y) < \epsilon$ (see for example Theorems 5.3.11 and 6.2.4 in \cite{FH20}).

Let us recall the notion of cohomology for functions with respect to flows. A continuous function $a\colon X \to \mathbb{R}$ is said to be \emph{$\Phi$-cohomologous to zero} if there exists a continuous function $q\colon X \to \mathbb{R}$ such that 
\[
a(x) = \lim_{t \to 0}\frac{q(\phi_t(x)) - q(x)}{t} \quad \textrm{for every $x \in X$}.
\] 

We say that a point $x \in X$ has a \emph{forward dense orbit} if $\overline{\{\phi_s(x): s \ge 0\}} = X$. When $\overline{\{\phi_s(x): s \in \mathbb{R}\}} = X$, we say that $x \in X$ has a \emph{dense orbit}. We say that a flow is \emph{topologically transitive} if there exists at least one point with a forward dense orbit.

The next result is a slightly more general version of the celebrated \emph{Liv\v{s}ic theorem} for flows (\cite{Liv72}). Unless explicitly stated, $X$ is always assumed to be a compact metric space. 
\begin{definition}\label{CCL} Let $\Phi = (\phi_t)_{t \in \R}$ be a continuous flow on a compact metric space X. $\Phi$ is said to satisfy the Closing Lemma if for every $\epsilon > 0$ there exists $\delta > 0$ such that if $x \in X$ and $t \ge 0$ satisfying $d(\phi_t(x),x)<\delta$, then there exists a periodic orbit $\{\phi_s(y): 0 \le s \le T \}$ with $|T - t| < \epsilon$ such that $d(\phi_s(x),\phi_s(y)) < \epsilon$ for all $0 \le s \le t$.
\end{definition}

Let us recall the notion of cohomology for functions. A continuous function $a: X \to \R$ is said to be $\Phi$-cohomologous to zero if there exists a continuous function $q: X \to \R$ such that 
\[
a(x) = \lim_{t \to 0}\frac{q(\phi_t(x)) - q(x)}{t} \quad \textrm{for every $x \in X$}.
\] 

We say that a point $x \in X$ has a \emph{forward dense orbit} if $\overline{\{\phi_s(x): s \ge 0\}} = X$. When $\overline{\{\phi_s(x): s \in \R\}} = X$, we say that $x \in X$ has a \emph{dense orbit}. We say that a flow is \emph{topologically transitive} if there exists at least one point with a forward dense orbit.

The following proposition is a more general version of the so-called \emph{Liv\v{s}ic Theorem} originally obtained in \cite{Liv72}.

\begin{theorem}\label{LIV}
Let $\Phi = (\phi_t)_{t \in \R}$ be a topologically transitive continuous flow satisfying the Closing Lemma, and $a\colon X \to \R$ a continuous function satisfying the Walters property. Then $a$ is cohomologous to zero if and only if for every periodic point $x = \phi_T(x)$ we have $S_Ta(x) = 0$.
\end{theorem}
\begin{proof}
See the proof of Theorem 5.3.23 in \cite{FH20}. 
\end{proof}
We also can obtain a characterization of additive families generated by coboundary functions.
\begin{proposition}\label{EQC}
Under the conditions of Theorem \ref{LIV}, a function $a: X \to \R$ is $\Phi$-cohomologous to zero if and only if 
\[
\lim_{t \to \infty}\frac1t\|S_ta\|_{\infty}=0.
\]
In particular, 
\[
\lim_{t \to \infty}\frac1t\|S_ta\|_{\infty}=0 \textrm{ if and only if } \sup_{t \ge 0}\|S_ta\|_{\infty} < \infty.
\]
\end{proposition}

\begin{proof}
Suppose $a$ is $\Phi$-cohomologous to zero. This implies the existence of a continuous function $q: X \to \R$ such that  $S_ta = q \circ \phi_t - q$ for all $t \ge 0$. Consequently, one has $\|S_ta\|_{\infty} \le 2\|q\|_{\infty} < \infty$ for every $t \ge 0$. Then, $\lim_{t \to \infty}\frac1t\|S_ta\|_{\infty} =~0$. 

Conversely, let $\lim_{n \to \infty}\frac1t\|S_ta\|_{\infty} = 0$. The Lebesgue's dominated convergence theorem together with the Birkhoff's ergodic theorem for flows gives that
\begin{equation}\label{0mu}
0 = \int_{X}\lim_{n \to \infty}\frac1tS_ta~d\mu = \int_{X}a~d\mu \quad \textrm{for all $\mu \in \cM(\Phi)$.}
\end{equation}
For all $x \in X$ with $x = \phi_T(x)$, the measure $(\int_{0}^{T}\delta_{\phi_s(x)}ds)/T$ is $\Phi$-invariant. In particular, identity \eqref{0mu} gives that $S_Ta(x) = 0$ for all $x \in X$ with $x = \phi_T(x)$.  Hence, by Theorem \ref{LIV} we conclude that $a$ is $\Phi$-cohomologous to zero, as desired. 
\end{proof}

Based on Proposition~\ref{EQC}, we give a definition of cohomology for asymptotically additive families of continuous potentials.

\begin{definition}\label{CHM} We say that an asymptotically (or almost) additive family of continuous functions $\cA = (a_t)_{t \ge 0}$ is $\Phi$-cohomologous to a constant if there exists a continuous function $a$ which is $\Phi$-cohomologous to a constant and such that
\[
\lim_{t \to \infty}\frac1t\|a_t - S_ta\|_{\infty}=0.
\]
\end{definition}

One can easily check that a family $\cA = (a_t)_{t \ge 0}$ is $\Phi$-cohomologous to a constant if and only if the sequence $(a_n/n)_{n \in \N}$ is uniformly convergent to a constant. In particular, $\cA$ is $\Phi$-cohomologous to zero if and only if $\lim_{t \to \infty}\frac1t\|a_t\|_{\infty}=0$. 

Observe that, in general, the classical definition of cohomology for a function is way stronger than the one we are suggesting for nonadditive sequences in Definition \ref{CHM}. On the other hand, Proposition \ref{EQC} also hints into a new definition which is still weaker than the classical one but stronger than Definition \ref{CHM}:

\begin{definition}\label{SCH} An asymptotically (or almost) additive family of continuous functions $\cA = (a_t)_{t \ge 0}$ is $\Phi$-cohomologous to a constant if there exists a continuous function $a$ $\Phi$-cohomologous to a constant and such that
\[
\sup_{n \in \N}\|a_t - S_ta\|_{\infty} < \infty.
\]
\end{definition}  
In this case, $\cA$ is $\Phi$-cohomologous to zero if and only if $\cA$ is uniformly bounded. 

In the next section, our main result gives a scenario where the definitions \ref{CHM} and \ref{SCH} are, in fact, equivalent for almost additive families (see Theorem \ref{BOU}). 

\subsection{A Liv\v{s}ic-type theorem for almost additive families}\label{mr}

Let $\Phi = (\phi_t)_{t \in \R}$ be a continuous flow on a compact metric space $X$. We say that a family of functions $\cA = (a_t)_{t \ge 0}$ has \emph{bounded variation} if there exists $\epsilon > 0$ such that
\[
\sup_{t \ge 0}\sup\{|a_t(x) - a_t(y)|: d_t(x,y) < \epsilon \} < \infty,
\]
where $d_t(x,y) := \max \{d(\phi_s(x), \phi_s(y)): s \in [0,t]\}$. Moreover, we say that $\cA$ has \emph{tempered variation} if 
\[
\limsup_{\epsilon \to 0}\lim_{t \to \infty}\frac{\gamma_t(\epsilon)}{t} = 0,
\] 
where $\gamma_t(\epsilon):= \sup\{|a_t(x) - a_t(y)|: d_t(x,y) < \epsilon \}$.

We note that if $\phi: X \to \R$ satisfies the Bowen property then the additive sequence $(S_t\phi)_{t \ge 0}$ has bounded variation. The Walters property for functions (additive families) also can be extended naturally to the nonadditive case. A family of functions $\cA = (a_t)_{t \ge 0}$ satisfies the \emph{Walters property} if for each $\kappa > 0$ there exists $\epsilon > 0$ such that for $x, y \in X$ and $t \ge 0$, we have that $d(\phi_s(x),\phi_s(y)) < \epsilon$ for every $s \in [0,t]$ implies $|a_t(x) - a_t(y)| < \kappa$. It is clear from the definitions that a family satisfying the Walters property also satisfies the Bowen property.

The following theorem is our main result in this part.
\begin{theorem}\label{BOU}
Let $\Phi = (\phi_t)_{t \ge 0}$ be a topologically transitive continuous flow on a compact metric space $X$ and satisfying the Closing Lemma. Let $\cB = (b_t)_{t \ge 0}$ be an almost additive family of continuous functions (with respect to $\Phi$) with bounded variation. Then, the following are equivalent:
\begin{enumerate}
\item $\lim_{t \to \infty}\|b_{t}\|_{\infty}/t = 0$;
		
\item $\sup_{t \ge 0}\|b_t\|_{\infty} < \infty$;
		
\item there exists $K > 0$ such that  $|b_t(p)| \le K$ for all $p \in X$ and $t \ge 0$ with $\phi_t(p) = p$. 
\end{enumerate}
\end{theorem}

The following result is a direct consequence of Theorem \ref{BOU}.

\begin{corollary}\label{EQ2}
Let $\cA = (a_t)_{t \ge 0}$ be an almost additive family of continuous functions with bounded variation. Then, a continuous function $a: X \to \R$ such that $(S_ta)_{t \ge 0}$ has bounded variation satisfies 
\[
\lim_{t \to \infty}\frac{1}{t}\|a_t - S_ta\|_{\infty} = 0 \quad \textrm{if and only if} \quad \sup_{t \ge 0}\|a_t - S_ta\|_{\infty} < \infty.
\]
	
In particular, if $(S_ta)_{t \ge 0}$ does not have bounded variation we must have
\[
\sup_{t \ge 0}\|a_t -S_ta\|_{\infty} = \infty.
\] 
\end{corollary}

Notice that Corollary \ref{EQ2} readily implies that the Bowen regularity problem is in fact equivalent to the uniform bound problem for topologically transitive flows satisfying the Closing Lemma. We also note that Theorem \ref{BOU} is an extension of Proposition \ref{EQC} to the case of almost additive families of functions. 

\begin{remark}
Observe that Proposition \ref{EQC} asks for the sequence of potentials to have the Walters property, which is stronger than the bounded variation condition. This is because the classical cohomology result obtained for a single potential is also stronger than the uniformly bounded one obtained in Theorem \ref{BOU}. In addition to that, as we shall see on subsection \ref{Liv-sec}, Theorem \ref{BOU} is also particularly related to Theorem 1.2 in \cite{Kal11}, where control over the periodic data implies control over the full data. 
\end{remark}

In order to prove Theorem \ref{BOU}, let us first obtain a more general auxiliary result. 

\begin{lemma}\label{UBI}
Let $\Phi = (\phi_t)_{t \in \R}$ be a continuous flow on a compact metric space $X$ and let $\cC = (c_t)_{t \ge 0}$ be an almost additive family of continuous functions with uniform constant $C > 0$ and such that $\lim_{t \to \infty}\| c_t \|_{\infty}/t =~0$. Then 
	
\begin{enumerate}
\item for every $\tau$-periodic point $x_0 \in X$, we have $\sup_{q \in \mathbb{N}}|c_{q\tau}(x_0)| \leq C$;
		
\item for every $\tau$-periodic point $x_0$, there exists a constant $L:=L(\tau) \ge 0$ (only depending on the period of $x_0$) such that $\sup_{t \ge 0}|c_t(x_{0})| \le L$;
		
\item we have
\[
\sup_{\mu \in \cM(\Phi)}\bigg|\int_X c_{t} d\mu\bigg| \leq C \quad \textrm{for all $t \ge 0$}. 
\]
\end{enumerate}
\end{lemma}

\begin{proof}
We proceed as in the  proof of Lemma 4 in \cite{HS24}. Since the family $\cC$ is almost additive with uniform constant $C>0$, one can see that
\begin{equation}\label{asfh}
\sum_{k=0}^{p-1}c_t \circ \phi_{kt} - (p-1)C \le c_{pt} \le \sum_{k=0}^{p-1}c_t \circ \phi_{kt} + (p-1)C
\end{equation}
for all $t \ge 0$ and $p \in\N$. 
	
Now suppose $x_0$ is a $\tau$-periodic point, that is, $\phi_{\tau}(x_0) = x_0$. If $t = q\tau$ for some $q \in \N$, then $\phi_{kt}(x_0) = \phi_{kq\tau}(x_{0}) = \underbrace{(\phi_{\tau} \circ \phi_{\tau} \circ \cdots \circ \phi_{\tau})}_{\text{\textrm{$kq$ times}}}(x_0) = x_0$ for all $k \in \N$. In particular, this implies that
\begin{equation}\label{ctp}
\lim_{p \to \infty} \frac{1}{p}\sum_{k=0}^{p-1}c_t(\phi_{kt}(x_{0})) = c_t(x_{0}).
\end{equation}
Since $\lim_{t \to \infty}\| c_t \|_{\infty}/t = 0$, it follows directly from \eqref{asfh} and \eqref{ctp} that
\begin{equation}\label{CC}
-C \le c_t(x_{0}) = c_{q\tau}(x_0) \le C.
\end{equation}
Since $q \in \N$ is arbitrary, we conclude the proof of item 1. 
	
In order to obtain item 2, let $x_0$ be a $\tau$-periodic point, consider $t = q\tau + r$ with $r \in (0,\tau)$ and fix the numbers
\begin{equation}\label{AB}
A(\tau) := \inf\bigg\{\inf_{x \in X}c_{s}(x): s \in [0,\tau]\bigg\} \quad \text{and} \quad B(\tau) := \sup\bigg\{\sup_{x \in X}c_{s}(x): s \in [0,\tau]\bigg\}.
\end{equation}
By almost additivity together with \eqref{CC} and \eqref{AB}, we have that
\[
-2C + A(\tau) \le -C + c_{qt}(x_0) + c_r(\phi_{q\tau}(x_0)) \le c_t(x_0)
\]
and
\[
c_t(x_0) \le c_{qt}(x_0) + c_r(\phi_{q\tau}(x_0)) + C \le 2C + B(\tau).
\]	
Therefore,
\[
L_1(\tau):= \min\{A(\tau) - 2C, -C\} \le c_t(x_{0}) \le \max\{B(\tau) + 2C, C\}:= L_2(\tau)
\]
for all $t \ge 0$. Taking $L = L(\tau):= \max\{|L_1(\tau)|,|L_2(\tau)|\}$, the item 2 is proved.
	
Now let us prove item 3. Suppose $\mu$ is a $\Phi$-invariant measure. Then, in particular, $\mu$ is also $\phi_{t}$-invariant for every $t \ge 0$. By applying that $\lim_{t \to \infty}\| c_{t} \|_{\infty}/t = 0$ in the inequalities \eqref{asfh} together with the Birkhoff's ergodic theorem applied to the map $\phi_{t}$, we obtain that
\[
-C \leq \int_{X} \lim_{p \to \infty} \frac{1}{p}\sum_{k = 0}^{p -1} c_{t}(\phi_{kt}(x)) d \mu(x) = \int_{X} c_{t} d \mu \leq C
\]
for all $t \ge 0$. Since the measure $\mu \in \cM(\Phi)$ is arbitrary, the lemma is proved. 
\end{proof}

\emph{Proof of Theorem \ref{BOU}.} 
The core idea in the proof is that any future position of any point in the space can always be well approximated by a point contained in a fixed finite piece of orbit. 

\begin{figure}[h!]
\centering
\begin{tikzpicture}[scale = 0.6]
	
%\filldraw[color=black!60, fill = black!5,very thick] (0,0) circle (1.5cm);
\filldraw[color=black!60, fill=white!5, dashed] (5.5,0) circle (2.0cm);
\filldraw[color=black!60, fill=white!5, dashed] (13.7,1.7) circle (2.0cm);
%\filldraw[color=black!60, fill=black!5, very thick] (11,0) circle (1.5cm);
	
%\filldraw[color=black!60, fill=black!5, very thick] (18,0) circle (1.5cm);
	
\draw[thick] plot [smooth,tension=1] coordinates{(-0.5,1) (3,4) (6,0) (13,2)};
\filldraw [black] (-0.5,1) circle (2pt);
\filldraw [black] (6,0) circle (2pt);
\filldraw [black] (13.0,2.0) circle (2pt);
	
\draw[thick] plot [smooth,tension=1] coordinates{(-0.5,-0.5) (1,1) (3,3) (4,0) (9,-1.6) (11.5,0) (13,1)};
\filldraw [black] (-0.5,-0.5) circle (2pt);
\filldraw [black] (4.6,-0.42) circle (2pt);
\filldraw [black] (12.73,0.8) circle (2pt);
	
%\draw[thick] plot [smooth,tension=1] coordinates{(14,4) (16,2) (16,0) (18,1)};
%\filldraw [black] (18,1) circle (2pt);
	
%\draw[thick] plot [smooth,tension=1] coordinates{(14,-2) (16,-1) (18.5,-1)};
%\filldraw [black] (18.5,-1) circle (2pt);
	
\node (x1) at (-0.5,-1.0) {$x$};
\node (y1) at (-0.5,0.5) {$z$};
\node (x2) at (4.8,-1.2) {$\phi_t(x)$};
\node (y2) at (6.1,0.5) {$\phi_s(z)$};
\node (x3) at (13.5,0.5) {$\phi_{\tau}(x)$};
\node (y3) at (14.0,2.5) {$\phi_{\Delta(z,\delta)}(z)$};
%\node (xn+1) at (18.4,-0.3) {$T^n_Z(x)$};
%\node (yn+1) at (19,1.9) {$T^n_Z(y) = \phi_{\tau_n(y)}(y)$};
%\node (phix) at (3,-3.6) {$\phi_s(x)$};
%\filldraw [black] (3,-3.0) circle (2pt);
%\node (phiy) at (4.2,4.4) {$\phi_s(z)$};
%\filldraw [black] (3,4.0) circle (2pt);
	
%\draw [black, -triangle 60] (1.55,1) -- (1.65,1.5) [right];
\draw [black, -triangle 60] (8.55,0.15) -- (8.65,0.17) [right];
%\draw [black, -triangle 60] (14,4) -- (14.7,3.5) [right];
%\draw [black, -triangle 60] (14,-2) -- (14.7,-1.48) [right];
%\draw [black, -triangle 60] (1.55,-1.5) -- (1.9,-2) [right];
\draw [black, -triangle 60] (8.55,-1.65) -- (8.9,-1.6) [right];
		
\end{tikzpicture}
\caption{Approximating any point in the space by a finite piece of a dense orbit.}\label{orb}
\end{figure}
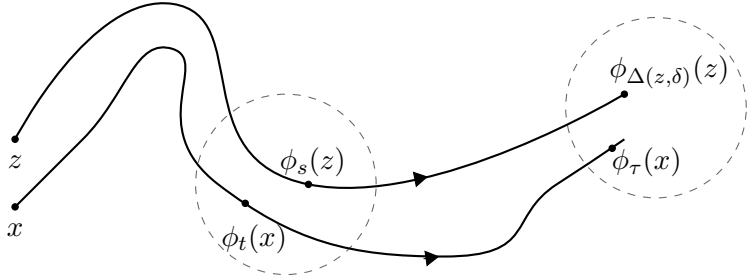
Let us start fixing the uniform bound given by the bounded variation property. In fact, since $\cB$ has bounded variation, there exists $\epsilon> 0$ such that
\begin{equation}\label{BVA}
Q:= \sup_{t \ge 0}\sup\{|b_t(x) - b_t(y)|: d_t(x,y) < \epsilon\} < \infty.
\end{equation}

Let's show that 3 implies 2. Suppose that there exists a uniform constant $K > 0$ such that  $|b_t(p)| \le K$ for all $p \in X$ and $t \ge 0$ with $\phi_t(p) = p$. Since $\Phi$ is topologically transitive there exists a point $z \in X$ with a dense forward orbit. Now let $\delta > 0$ be the number given by the Closing Lemma (Definition \ref{CCL}). By the density of the forward orbit, there exists a number $\Delta(z,\delta) \in \R$ such that for each $x \in X$ and $t \in \R$ there exists some $s \in [0,\Delta(z,\delta)]$ with $d(\phi_t(x),\phi_s(z)) < \delta$ (see Figure \ref{orb}). For $t > \Delta(z,\delta)$, in particular, there exists $s' \in [0,\Delta(z,\delta)]$ such that $d(\phi_t(z), \phi_{s'}(z)) < \delta$, which is the same as $d(\phi_{t-s'}(\phi_{s'}(z)), \phi_{s'}(z)) < \delta$. 

By the Closing Lemma,  there exists a point $p \in X$ with $\phi_{T}(p) = p$ with $|T - t + s'| < \epsilon$ and such that $d_{t-s'}(\phi_{s'}(z),p) < \epsilon$. From almost additivity, there exists a uniform constant $L = L(\epsilon) > 0$ such that $\|b_{T} - b_{t-s'}\|_{\infty} \le L$. Applying the bounded variation condition \eqref{BVA} we have that
\[
|b_{t-s'}(\phi_{s'}(z)) - b_{t-s'}(p)| \le Q,
\]
which gives that  $|b_{t-s'}(\phi_{s'}(z)| \le Q + |b_{t-s'}(p)| \le Q + |b_{T}(p)| + L \le Q + K + L$. By using almost additivity again, we get
\[
\begin{split}
|b_t(z)| = |b_{(t-s')+s'}(z)| &\le |b_{s'}(z)| + |b_{t-s'}(\phi_{s'}(z))| + C \\ 
& \le \sup_{s \in [0,\Delta(z,\delta)]}|b_s(z)| + Q + K + L + C := \widetilde{K}.
\end{split}
\]
Since the time $t > \Delta(z,\delta)$ was arbitrarily chosen and also taking into consideration times $0 < t \le \Delta(z,\delta)$, there exists $K_0 \ge \widetilde{K}$ such that $|b_t(z)| \le K_0$ for all $t \ge 0$. Notice that the constant $K_0 > 0$ depends only on $z$, $\delta > 0$ and $\epsilon > 0$.

Applying the almost additivity property one more time, we have
\[
|b_t(\phi_{s}(z)| \le |b_s(z)| + |b_{t+s}(z)| + C \le 2K_0 + C  \quad \textrm{for all $t, s \ge 0$}.
\]
Now consider any point $x \in X$. Since $\overline{\{\phi_t(z): t \ge 0\}} = X$, there exists a sequence of points $(z_q)_{q \ge 1} \subset \{\phi_t(z): t \ge 0\}$ such that $\lim_{q \to \infty}z_q = x$. Since every function $b_t: X \to \R$ is continuous, we obtain that
\[
|b_t(x)| = \lim_{q \to \infty}|b_t(z_q)| \le 2K_0 + C.
\]
Hence, by the arbitrariness of $x$, we conclude that $\sup_{t \ge 0}\|b_t\|_{\infty} \le 2K_0+ C < \infty$, as desired. It is obvious that 2 implies 1. Moreover, by Lemma \ref{UBI}, we have that 1 implies 3 and the theorem is proved.  \qed

\begin{remark}
Notice that, with the exception of the core idea coming from the counterpart result for maps in \cite{HS24}, the proof of Theorem \ref{BOU} does not use any direct tool from the discrete-time case and also does not depend on any physical equivalence results with respect to maps or flows. 
\end{remark}

Theorem \ref{BOU} does not hold for asymptotically additive nor subadditive families in general. In fact, let $\Phi$ be any continuous flow on a compact metric space $X$ and consider the family $\cA = (a_t)_{t \ge 0}$ given by $a_t(x) := \sqrt t$ for all $t \ge 0$ and $x \in X$. The family $\cA$ has bounded variation and is asymptotically additive and also subadditive with respect to $\Phi$. Moreover, it is clear that $\lim_{t \to \infty}\|a_t\|_{\infty}/t = 0$ but $\sup_{t \ge 0}\|a_t\|_{\infty} = \infty$. 
In fact, proceeding as in Example 1 in \cite{HS24}, one can show that
\[
\sup_{t \ge 0}\bigg\|a_t - \int_{0}^{t}(a \circ \phi_s)ds\bigg\|_{\infty} = \infty \quad \textrm{for every continuous function $a\colon X \to \R$}.
\]
	
This simple example shows that Theorem \ref{BOU} has the optimal nonadditive setup in the sense that it cannot be extended to more general classes of families. Moreover, it also indicates that definitions \ref{CHM} and \ref{SCH} are not equivalent for asymptotically additive families in general.  

\subsubsection{An application to linear cocycles}\label{Liv-sec}

In this section we follow closely some definitions and notions of \cite{BH21b}. 

Let $\Phi = (\phi_t)_{t \in \R}$ be a continuous flow on a compact metric space~$M$. Moreover, let $GL(d,\R)$ be the set of all invertible $d \times d$ matrices. A continuous map $\cA\colon \R \times M \to GL(d,\R)$ is called a \emph{linear cocycle} over $\Phi$ if for all $t,s\in\R$ and $x\in M$ we have:
\begin{enumerate}
\item $\cA(0,x) =$ Id;	
\item $\cA(t+s,x) = \cA(s,\phi_t(x))\cA(t,x)$.
\end{enumerate}
We shall always assume that all entries $a_{ij}(t,x)$ of $\cA(t,x)$ are positive for every $(t,x) \in \R \times M$. Moreover, for definiteness we shall consider the norm on $GL(d,\R)$ defined by $\|B\| = \sum_{i,j=1}^d |b_{ij}|$, denoting by $b_{ij}$ the entries of the matrix $B$.

Now we consider the family of continuous functions $a := (a_t)_{t \ge 0}$ defined by
\[
a_t(x)= \log\|\cA(t,x)\| \quad \textrm{for all $t \ge 0$ and $x \in M$.}
\]
It follows from Proposition 12 in \cite{BH21b} that the family $a$ is almost additive with respect to the flow $\Phi$. Notice that for a general linear cocycle, the family $a$ is only subadditive. 

We say that the cocycle $\cA$ has \emph{tempered distortion} if 
\[
\limsup_{t \to \infty}\frac{1}{t}\log \sup\left\{\|\cA(t,x) \cA(t,y)^{-1}\|: z \in M \ \text{and} \ x, y \in B_{t}(z,\epsilon)\right\} = 0
\]
for some $\epsilon > 0$. Moreover, we say that $\cA$ has \emph{bounded distortion} if 
\[
\sup\left\{\|\cA(t,x) \cA(t,y)^{-1}\|: z \in M \ \text{and} \ x, y \in B_{t}(z,\epsilon)\right\} < \infty
\]
for some $\epsilon > 0$. Clearly, bounded distortion implies tempered distortion. 

Now observe that
\[
\|\cA(t,x) \cA(t,x)^{-1}\| = \|\Id\| = d
\]
for every $(t,x) \in \R\times M$, which implies that
\[
\|\cA(t,x)^{-1}\| \ge d \|cA(t,x)\|^{-1}.
\]
Then,
\[
\|\cA(t,x) \cA(t,y)^{-1}\| 
\ge \frac{K}{d}\|\cA(t,x)\|\cdot\|\cA(t,y)^{-1}\| \ge K\|\cA(t,x)\|\cdot\|\cA(t,y)\|^{-1}
\]
and so
\[
\bigl\lvert\log \|\cA(t,x)\| - \log \|\cA(t,y)\|\bigr\rvert \le - \log K + \log \|\cA(t,x) \cA(t,y)^{-1}\|.
\]
In particular, for $z \in M$ and $\epsilon > 0$ we have
\[
\sup_{x, y \in B_t(z,\epsilon)}|a_t(x) - a_t(y)| \le -\log K + \log \sup_{x, y \in B_t(z,\epsilon)}\|\cA(t,x) \cA(t,y)^{-1}\|.
\]
Hence, if $\cA$ has tempered distortion, then the family $a$ has tempered variation, and if $\cA$ has bounded distortion, then $a$ has bounded variation.

For a concrete example of linear cocycles over flows, one can consider a $C^1$ flow $\Phi$ on a compact set $M \subset \R^d$ such that for every $t \in\R$ and $x \in M$ the matrix $d_x\phi_t$ has only positive entries. Then $\cA(t,x)=d_x\phi_t$ is a linear cocycle over $\Phi$ and the family~$a$ defined by $a_t(x)= \log\|d_x \phi_t\|$ is an almost additive family of continuous functions with respect to~$\Phi$.

Now let $GL^{+}(d,\R) \subset GL(d,\R)$ be the set of all matrices with strictly positive entries. We have the following direct application of Theorem \ref{BOU} to the case continuous-time cocycles:

\begin{proposition}\label{KLN}
Let $\Phi = (\phi_t)_{t \in \R}$ be a topologically transitive continuous flow on a compact metric space $M$ satisfying the Closing Lemma. Let $\cA: \R \times M \to GL^{+}(d,\R)$ be a linear cocycle over $\Phi$ and with bounded distortion. Suppose there exists a compact set $\Omega \subset GL^{+}(d,\R)$ such that $\cA(t,p) \subset \Omega$ for all $t \ge 0$ and $p \in M$ with $\phi_t(p) = p$. Then there exists a compact set $\widetilde{\Omega}$ such that  $\cA(t,x) \subset \widetilde{\Omega}$ for all $t \ge 0$ and $x \in M$. 
\end{proposition}

\begin{proof}
By the hypotheses, the family of continuous functions $a := (a_t)_{t \ge 0}$ given by $a_t(x):= \log\|\cA(t,x)\|$ is almost additive with respect to $\Phi$. Moreover, since the cocycle $\cA$ has bounded distortion, the family $a$ has bounded variation. Now suppose there exists a compact $\Omega \subset GL^{+}(d,\R)$ where $\cA(t,p) \subset \Omega$ for all $t \ge 0$ and $p \in M$ with $\phi_t(p) = p$. Since the map $\cA(t,p) \mapsto \log\|\cA(t,p)\|$ is continuous, there exists $K > 0$ such that $|a_t(p)|\le K$ for all $t \ge 0$ and all $p \in M$ with $\phi_t(p) = p$. By Theorem \ref{BOU}, there exists a constant $\widetilde{K} > 0$ such that $\sup_{t \ge 0}\|a_t\|_{\infty} \le \widetilde{K}$. In particular, this implies that $e^{-\widetilde{K}} \le \|\cA(t,x)\| \le e^{\widetilde{K}}$ for all $t \ge 0$ and all $x \in M$.
Hence, we conclude that
\[
\|\cA(t,x) - \Id\| \le \|\cA(t,x)\| + \|\Id\| \le e^{\widetilde{K}} + d \quad \textrm{for all $t \ge 0$ and $x \in M$},
\]
as desired.
\end{proof}

\begin{remark}
Proposition \ref{KLN} is a particular continuous-time counterpart of Theorem 1.2 in \cite{Kal11}, where uniformly bounded periodic data guarantee a uniform bound for the entire phase space.
\end{remark}

\subsubsection{Nonadditive notions of (weak) Gibbs measures for maps and flows}

In this section, we are going to consider families of functions derived from (not necessarily invariant) measures. We will see later that those types of natural examples play an important role in the understanding of some regularity problems arising from the physical equivalence relation for flows. Let us first recall some relevant concepts and notions of the nonadditive thermodynamic formalism for maps.

Let $T\colon X \to X$ be a continuous map of a compact metric space, and let $\cF = (f_n)_{n \in \N}$ be an almost additive sequence of continuous functions with tempered variation (a condition weaker than bounded variation). We have the following variational principle (see \cite{Bar06} and \cite{Mum06}: 
\[
P_T(\cF) = \sup_{\mu \in \cM(T)}\bigg(h_{\mu}(T) + \lim_{n \to \infty}\frac1n\int_{X}f_n d\mu\bigg),
\]
where $P_T(\cF)$ is the nonadditive topological pressure of $\cF$ with respect to $T$, and $h_{\mu}(T)$ is the Kolmogorov-Sinai entropy. A measure $\nu \in \cM_T$ is said to be an \emph{equilibrium measure} for $\cF$ (with respect to $T$) if the supremum is attained in $\nu$, that is,
\[
P_T(\cF) = h_{\nu}(T) + \lim_{n \to \infty}\frac1n\int_{X}f_n d\nu.
\]
When $\cF$ is an additive sequence generated by a single continuous function $f: X \to \R$, one can easily see that
\[
P_T(\cF) = P_T(f) \quad \textrm{and} \quad \lim_{n \to \infty}\frac1n\int_{X}f_n d\mu = \int_{X}fd\mu \quad \textrm{for all $\mu \in \cM(T)$},
\]
where $P_T(f)$ is the classical (additive) topological pressure of $f$ with respect to the continuous map $T$. For asymptotically additive sequences in general, Theorem 1.2 in \cite{Cun20} allows us to obtain the notion of nonadditive topological pressure and variational principle directly from the classical theory. 

\begin{definition}\label{GB1}
We say that a measure $\mu$ on $X$ (not necessarily $T$-invariant) is a \emph{Gibbs measure} with respect to an asymptotically additive (almost additive) sequence $\cF = (f_n)_{n \in \N}$ if for any sufficiently small $\epsilon > 0$ there exists a constant $K(\epsilon) \ge 1$ such that
\[
K(\epsilon)^{-1} \le \frac{\mu(B_n(x,\epsilon))}{\exp[-nP_T(\cF) + f_n(x)]}\le K(\epsilon)
\]
for all $x \in X$ and $n\in\N$, where $B_n(x,\epsilon)$ is the Bowen ball given by
\[
B_n(x,\epsilon):= \{y \in X: d(T^{k}(x),T^{k}(y)) < \epsilon \textrm{ for all $0 \le k \le n-1$}\}.
\]
\end{definition}

It was introduced in \cite{Bar06} a definition of Gibbs measures with respect to almost additive sequences using Markov partitions. In this case, one can show that this notion of Gibbs measures and the one in Definition \ref{GB1} are equivalent when the system admits Markov partitions with arbitrarily small diameter, as in the case of locally maximal hyperbolic sets or repellers for $C^{1}$ diffeomorphisms (see \cite{Bow75a}). 

In the case of shifts $\sigma: \Sigma^{\N} \to \Sigma^{\N}$, the definition is simpler. A measure $\mu$ on $\Sigma^{\N}$ is said to be \emph{Gibbs} with respect to $\cF = (f_n)_{n \in \N}$ when there exists a constant $K \ge 1$ such that
\[
K^{-1} \le \frac{\mu_n(C_{i_1 \ldots i_n})}{\exp(-nP_{\sigma}(\cF) + f_n(x))} \le K
\]

for all $x \in C_{i_1 \ldots i_n}$ and $n \in \N$, where $C_{i_1\ldots i_n}$ is the \emph{cylinder} set 
\[
C_{i_1 \ldots i_n}:= \{y = (j_1 j_2 \cdots) \in \Sigma^{\N}: j_1 = i_1,\ldots, j_n = i_n\}.
\]

In the same way, we also have the following

\begin{definition}
We say that a measure $\mu$ on $X$ (not necessarily $T$-invariant) is a \emph{weak Gibbs measure} with respect to an almost additive sequence $\cF = (f_n)_{n \in \N}$ if for any sufficiently small $\epsilon > 0$ there exists a sequence $(K_n(\epsilon))_{n \in \N} \subset [1, \infty)$ with $\lim_{n \to \infty}\log K_n(\epsilon)/n = 0$ such that
\[
K_n(\epsilon)^{-1} \le \frac{\mu(B_n(x,\epsilon))}{\exp[-nP_T(\cF) + f_n(x)]}\le K_n(\epsilon)
\]
for all $x \in X$ and $n\in\N$
\end{definition}

These definitions are natural nonadditive versions of the classical Gibbs definitions with respect to a single continuous function.

Let us now consider the case of continuous-time dynamical systems. Let $X$ be a compact metric space, $T\colon X \to X$ a continuous map and $\tau\colon X \to (0, \infty)$ a positive continuous function. Consider the space
\[
W = \bigl\{(x,s) \in X \times \R: 0 \le s \le \tau(x)\bigr\}
\]
and let $Y$ be the set obtained from $W$ identifying $(x,\tau(x))$ with $(T(x),0)$ for each $x \in X$. By the Bowen-Walters distance (see \cite{BW72}), the space $Y$ is a compact metric space. Recall that the \emph{suspension semi-flow} over $T$ with height function $\tau$ is the semi-flow $\Phi = (\phi_t)_{t \ge 0}$ on $Y$ with the maps $\phi_t\colon Y \to Y$ defined by $\phi_t(x,s) = (x, s+t)$ whenever $t + s \in [0, \tau(x)]$.

Let $X$ be a Riemannian manifold, $T:X \to X$ a $C^{1}$ map, and let $\Lambda \subset X$ be a compact $T$-invariant set, that is, $T^{-1}(\Lambda) = \Lambda$. The map $T$ is said to be \emph{uniformly expanding} on $\Lambda$ if there exist constants $c > 0$ and $\lambda > 1$ such that 
\[
\|d_xf^{n}v\| \ge c\lambda^{n}\|v\| \quad \textrm{for all $x \in \Lambda$, $n \in \N$ and $v \in T_xM$} 
\]
In this case, the set $\Lambda$ is called a \emph{repeller} for the map $T$. For instance, $T: [0,1] \to [0,1]$ given by $T(x) = 2x \mod 1$ is a uniformly expanding map and the whole closed interval $[0,1]$ is a repeller for $T$.

Let $M$ be any compact metric space. A map $S: M \to M$ is said to satisfy the \emph{property P} if for each H\"older continuous function $\xi: M \to \R$ there exists a constant $D > 0$ such that if $x \in M$, $n \in \N$ and $d(T^{k}(x),T^{k}(y)) < \epsilon$ for all $k \in \{0,...,n-1\}$, then 
\[
\bigg|\sum_{k=0}^{n-1}\xi(S^{k}(x)) - \sum_{k=0}^{n-1}\xi(S^{k}(y))\bigg| < D\epsilon.
\]

Notice that full shifts, subshifts of finite type, uniformly expanding and hyperbolic maps all satisfy the property P (see for example \cite{Wal78} and \cite{Bou02}). 

\begin{remark}
In \cite{BS00} the property P is called bounded distortion. We intentionally changed the name in order not to cause too much confusion with the definitions given in subsections \ref{bdf} and \ref{Liv-sec}.  
\end{remark}

Now consider $\Phi$ the suspension semi-flow over a continuous map $T: X \to X$ satisfying the property P and with H\"older continuous height function $\tau$. Proposition 19 in \cite{BS00} guarantees that for each sufficiently small $\epsilon > 0$ there exists a constant $\kappa > 0$ such that 
\begin{equation}\label{BY1}
	B^{Y}_{\tau_m(x)}((x,s),\epsilon) \subset B^{X}_m(x,\kappa\epsilon) \times (s-\kappa\epsilon, s+\kappa\epsilon),
\end{equation}
\begin{equation}\label{BY2}
	B^{X}_m(x,\epsilon/\kappa) \times (s-\epsilon/\kappa, s+\epsilon/\kappa) \subset B^{Y}_{\tau_m(x)}((x,s),\epsilon)
\end{equation}
for every $x \in X$, $0 < s < \tau(x)$ and $m \in \mathbb{N}$, where 
\[
\tau_{n}(x) = \sum_{k=0}^{n-1}r(T^{k}(x)) \quad \textrm{for all $x \in X$},
\]
\[
B^{Y}_t(y, \delta) := \{z \in Y: d(\phi_{s}(z),\phi_{s}(y)) < \delta \textrm{ for all $s \in [0,t]$}\},
\] 
\[
B^{X}_n(x,\delta) := \{w \in X: d(T^{k}(w),T^{k}(x)) < \delta \textrm{ for all $0 \le k \le n-1$}\}.
\] 

Let $a = (a_t)_{t \ge 0}$ be a family of almost additive continuous functions with respect to the flow $\Phi$. Following as in the proof of Lemma 3.1 in \cite{BH21a}, the sequence $c = (c_n)_{n \in \N}$ given by $c_n(x) = a_{\tau_n(x)}(x)$ is almost additive with respect to $T:X \to X$. Now consider $\mu$ a Gibbs measure for the sequence $c$ on $X$ and let $\nu$ be the measure on $Y$ induced by $\mu$ (see identity \eqref{eq13})\footnote{For the almost additive thermodynamics with respect to maps, proper definitions of Gibbs and weak-Gibbs measures for sequences, see \cite{Bar06,Mum06} and the review \cite{Bar10}).}. In particular, $\nu = (\mu \times \lambda)/(\int_{X}\tau d\mu)$, where $\lambda$ is the Lebesgue measure on $\R$. Since $\mu$ is Gibbs, for any sufficiently small $\epsilon > 0$ there exist $K_1(\epsilon) > 0$ and $K_2(\epsilon) > 0$ such that 
\begin{equation}\label{K1}
K_1(\epsilon)^{-1}\exp[-mP_T(c) + c_m(x)] \le \mu(B^{X}_m(x,\kappa \epsilon)) \le K_1(\epsilon)\exp[-mP_T(c) + c_m(x)]
\end{equation}
\begin{equation}\label{K2}
K_2(\epsilon)^{-1}\exp[-mP_T(c) + c_m(x)] \le \mu(B^{X}_m(x,\epsilon/\kappa)) \le K_2(\epsilon)\exp[-mP_T(c) + c_m(x)]
\end{equation}
for all $x \in X$ and $m \in \N$.

Observe that by identity \eqref{eq15} and the definition of $\tau_m$, we have that 
\[
\bigg(\frac{1}{\sup \tau}\bigg)P_T(c) \le P_{\Phi}(a) \le \bigg(\frac{1}{\inf \tau}\bigg)P_T(c) \quad \textrm{and} \quad m\inf \tau \le \tau_m(x) \le m \sup\tau.
\]
Moreover, one can check that for all $t > 0$ there exists $m \in \N$ such that $\tau_{m}(x) \le t \le \tau_{m+1}(x)$ with $t-\tau_m(x) \in [0,\sup\tau]$, which gives $|a_t(x) - a_{\tau_m(x)}(x)| \le \sup_{s \in [0,\sup\tau]}\|a_s\|_{\infty}=:~q$. Hence, it follows from \eqref{BY1} and \eqref{K1} that  
\[
\begin{split}
\nu(B^{Y}_{t}(\phi_s(x),\epsilon)) &\le \nu(B^{Y}_{\tau_m(x)}(\phi_s(x),\epsilon))\\ &\le \frac{2\kappa\epsilon K_1(\epsilon)}{\inf \tau}\exp[-\tau_m(x)P_{\Phi}(a) + a_{\tau_m(x)}(x)]\\
& \le \frac{2\kappa\epsilon K_1(\epsilon)}{\inf \tau}\underbrace{\exp[\sup\tau P_{\Phi}(a) + q]}_{\text{$=:L_1$}}\exp[-tP_{\Phi}(a) + a_{t}(x)]\\
&=  \frac{2\kappa\epsilon L_1K_1(\epsilon)}{\inf \tau}\exp[-tP_{\Phi}(a) + a_{t}(x)]
\end{split}
\]
for all $x \in X$ and $s \in [0, \tau(x)]$. By the almost additivity of the family $a$, we also have that $|a_t(x) - a_t(\phi_s(x))| \le 2 \sup_{s \in [0,\sup\tau]}\|a_s\|_{\infty}=: \widetilde{q}$. Since for each $y \in Y$ there exist $x \in X$ and $s \in [0,\sup\tau]$ such that $y = \phi_s(x)$, we finally get that 
\[
\nu(B^{Y}_{t}(y,\epsilon)) \le  \frac{2\kappa\epsilon K_1(\epsilon)e^{\widetilde{q}}}{\inf \tau}\exp[-tP_{\Phi}(a) + a_{t}(y)] = \widetilde{K_1}(\epsilon)\exp[-tP_{\Phi}(a) + a_{t}(y)]
\]
for all $y \in Y$ and $t > 0$, where $\widetilde{K_1}(\epsilon):= (2\epsilon L_1 K_1(\epsilon)e^{\widetilde{q}})/\inf \tau$ only depends on $\epsilon > 0$ and the strictly positive function $\tau$. 

Proceeding in the same way, the identities \eqref{BY2} and \eqref{K2} guarantee the existence of a constant $\widetilde{K_2}(\epsilon) > 0$ such that 
\[\nu(B^{Y}_{t}(y,\epsilon))\ge \widetilde{K_2}(\epsilon)\exp[-tP_{\Phi}(a) + a_{t}(y)] \quad \textrm{for all $y \in Y$ and $t > 0$.}
\]

So, we just showed that Gibbs measures for the almost additive sequence on the base space induce a Gibbs property for the almost additive family with respect to the flow. Analogously, one can show that weak Gibbs measures for the asymptotically additive sequence on the base also induce measures which satisfy a weak Gibbs property for the asymptotically additive family with respect to the flow. These relations involving (weak) Gibbs properties between the map on base space and the flow also hold for suspension flows over maps with the property P, including locally maximal hyperbolic sets for topologically mixing $C^{1}$ flows. Based on this, we have the following definitions:

\begin{definition}\label{GBf}
Let $\Phi = (\phi_t)_{t \in \R}$ be a continuous flow on a compact metric space $M$. We say that a measure $\mu$ on $M$ (not necessarily $\Phi$-invariant) is a \emph{Gibbs measure} for an asymptotically additive family of continuous functions $a = (a_t)_{t \ge 0}$ (with respect to $\Phi$) if for any sufficiently small $\epsilon > 0$ there exists a constant $K(\epsilon) \ge 1$ such that 
\[
K(\epsilon)^{-1} \le \frac{\mu(B_t(x,\epsilon))}{\exp[-tP_{\Phi}(a) + a_t(x)]}\le K(\epsilon)
\]
for all $x \in M$ and $t > 0$, where $B_t(x,\epsilon)$ is the usual Bowen ball with respect to $\Phi$ and $P_{\Phi}(a)$ is the nonadditive topological pressure for $a$ with respect to the flow $\Phi$. 
\end{definition}

\begin{definition}\label{WGBf}
Let $\Phi = (\phi_t)_{t \in \R}$ be a continuous flow on a compact metric space $M$. We say that a measure $\mu$ on $M$ (not necessarily $\Phi$-invariant) is a \emph{weak Gibbs measure} for an asymptotically additive family of continuous functions $a = (a_t)_{t \ge 0}$ (with respecto to $\Phi$) if for any sufficiently small $\epsilon > 0$ there exists a sequence $(K_t(\epsilon))_{t > 0} \subset [1, \infty)$ with $\lim_{t \to \infty}\log K_t(\epsilon)/t = 0$ such that
\[
K_t(\epsilon)^{-1} \le \frac{\mu(B_t(x,\epsilon))}{\exp[-tP_{\Phi}(a) + a_t(x)]}\le K_t(\epsilon)
\]
for all $x \in M$ and $t > 0$. 
\end{definition}

\begin{remark}
In the case of hyperbolic flows or suspension flows over subshifts of finite type, uniformly expanding or hyperbolic maps in general, Definition \ref{GBf} is a generalization of the classical notion of Gibbs measures to the nonadditive setup (see for example Definition 4.3.25 in \cite{FH20}). 
\end{remark}

The following result guarantees the existence of Gibbs measures for almost additive families of functions with respect to hyperbolic flows.
\begin{proposition}[{\cite[Theorem~3.5]{BH21a}}]\label{H20}
Let $\Lambda$ be a hyperbolic set for a topologically mixing $C^1$ flow~$\Phi$ and let $a$ be an almost additive family of continuous functions on $\Lambda$ with bounded variation.  Then:
\begin{enumerate}
\item
there exists a unique equilibrium measure for $a$;
\item
there exists a unique invariant Gibbs measure for $a$;
\item
the two measures are equal and are ergodic.
\end{enumerate}
\end{proposition}

\begin{remark}
We note that Proposition \ref{H20} is also valid for appropriate versions of suspension flows over subshifts of finite type, uniformly expanding or hyperbolic maps with H\"older continuous height functions. 
\end{remark}

The existence of Gibbs measures in Proposition \ref{H20} was originally obtained for the equilibrium measure on the base space using Markov partitions for the flow. As we saw above, this is the same as having the induced $\Phi$-invariant measure satisfying the Definition~\ref{GBf}. Moreover, in the case of asymptotically additive families in the same setup of Proposition~\ref{H20}, we cannot guarantee uniqueness of equilibrium states and one can easily check that the measures lifted from the base are only weak Gibbs in the sense of Definition \ref{WGBf}. \\

For the sake of completeness, now we will show other natural sources of almost and asymptotically additive families of functions. Proceeding as in the proof of Proposition 12 in \cite{HS24}, one can check that if a measure $\eta$ on a compact metric space $M$ is Gibbs for some almost additive family of continuous functions with respect to a flow $\Phi = (\phi_t)_{t \in \R}$ on $M$, then for every sufficiently small $\delta > 0$ there exists a constant $K(\delta) > 1$ such that 
\[
\frac{1}{K(\delta)} \le \frac{\eta(B_{t+s}(x,\delta))}{\eta(B_t(x,\delta))\eta(B_s(\phi_t(x),\delta))} \le K(\delta) \quad \textrm{for all $x \in M$ and $t,s > 0$.}
\]

In particular, for each $\delta > 0$, the well defined family of functions $a^{\delta}:= (a^{\delta}_t)_{t \ge 0}$ given by $a^{\delta}_t(x) = \log\eta(B_t(x,\delta))$ is almost additive (defining $a_0 \equiv 0$). Since every family of functions admitting a Gibbs measure clearly has bounded variation (following the same argument as in the proof of Proposition 15 in \cite{HS24}), the family $a^{\delta}$ also satisfies the bounded variation condition. On the other hand, notice that the functions $a_t^{\delta}: M \to \R$ are not necessarily continuous. In fact, using the Gibbs property of $\eta$, in general one can only guarantee the existence of constants $K_1(\delta) \ge K_{2}(\delta) > 0$ such that 
\[
K_2(\delta) + \limsup_{x \to x_0}a_t^{\delta}(x) \le a_t^{\delta}(x_0) \le K_1(\delta) + \liminf_{x \to x_0}a_t^{\delta}(x) 
\]
for all $x_0 \in M$ and all $t\ge0$. In particular, $x \mapsto a^{\delta}_t(x)$ is upper semicontinuous for each $t \ge 0$.  Interestingly, the thermodynamic formalism for discontinuous additive and nonadditive families of potentials with respect to flows was initiated quite recently in \cite{TLS25,TLS24}, including a more detailed treatment of the upper semicontinuous case.

\begin{proposition}\label{vol}
Let $\Phi = (\phi_t)_{t \in \R}$ be a continuous flow on a compact metric space $M$, and let $\nu$ be any measure on $M$. Then:
\begin{enumerate}

\item if for some $\delta > 0$ there exist constants $A(\delta) \ge B(\delta) > 0$ and an almost additive family of continuous functions $g = (g_t)_{t \ge 0}$ such that 
\[
B(\delta)e^{g_t(x)}\le \nu(B_t(x,\delta)) \le A(\delta)e^{g_t(x)} \quad \textrm{for all $x \in M$ and $t \ge 0$},
\] 
then there exists an almost additive family of H\"older continuous functions $h=(h_t)_{t \ge 0}$ satisfying $\sup_{t \ge 0}\sup_{x \in M}|\log\nu(B_t(x,\delta)) - h_t(x)| < \infty$;

\item if for some $\delta > 0$ there exist sequences $(C_t(\delta))_{t \ge 0}, (D_t(\delta))_{t \ge 0} \subset [1,\infty)$ such that $\log C_t(\delta)/t \to 0$, $\log D_t(\delta)/t \to 0$ and an asymptotically additive family of continuous functions $f = (f_t)_{t \ge 0}$ such that 
\[
D_t(\delta)e^{f_t(x)}\le \nu(B_t(x,\delta)) \le C_t(\delta)e^{f_t(x)} \quad \textrm{for all $x \in M$ and $t \ge 0$}
\] 
then there exists an asymptotically additive family of H\"older continuous functions  $j= (j_t)_{t \ge 0}$ satisfying 
\[
\lim_{t \to \infty}\frac{1}{t}\sup_{x \in M}|\log\nu(B_t(x,\delta)) - j_t(x)| = 0.
\]
\end{enumerate}
\end{proposition}

\begin{proof}
By using the fact that the space of H\"older continuous functions are dense in the space of continuous functions on compact spaces, we can guarantee the existence of families $h = (h_{t})_{t \ge 0}$ and $j = (j_t)_{t \ge 0}$ composed of H\"older continuous functions and such that $\sup_{x \in M}|g_t(x) - h_t(x)| \le 1$ and  $\sup_{x \in M}|f_t(x) - j_t(x)| \le 1$ for all $t \ge 0$. It is clear that $h$ is almost additive and $j$ is asymptotically additive. 
\end{proof}

\begin{example}
Observe that Gibbs measures satisfy the conditions of the first item and weak Gibbs measures satisfy the conditions of the second item of Proposition \ref{vol}. 
Now let $M$ be a compact Riemannian manifold and $\Lambda \subset M$ a hyperbolic set for a topologically mixing (at least) $C^1$ flow~$\Phi$. For each $t > 0$, consider the continuous function $J_t: \Lambda \to \R$ given by $J_t(x) = -\log\|d_x\phi_t|_{E^{u}(x)}\|$. Since $x \mapsto E^{u}(x)$ is H\"older continuous, we also have that $x \mapsto J_t(x)$ is H\"older (see subsection \ref{subsec32}). Now let $\lambda$ be the Lebesgue measure on $M$. Assuming that $\Phi$ is at least $C^{2}$, the \emph{Volume Lemma} (see Proposition 7.4.3 in \cite{FH20}) says that for any sufficiently small $\delta > 0$ there exist constants $C_{\delta}, D_{\delta} > 0$ such that 
\[
D_{\delta}J_t(x) \le \lambda(B_t(x,\delta)) \le C_{\delta}J_t(x) \quad \textrm{for all $x \in \Lambda$ and $t \ge 0$.}
\]
Moreover, one can also check that the family $(J_t)_{t \ge 0}$ is additive with respect $\Phi$. Notice that, in this case, the measure $\lambda$ satisfy the conditions in the first item of Proposition~\ref{vol}. Hence, the family $Leb^{\delta} = (Leb^{\delta}_{t})_{t \ge 0}$ given by $Leb^{\delta}_{t}(x):= \log \lambda(B_t(x,\delta))$ is almost additive with bounded variation and  physically equivalent to an almost additive family of Lipschitz continuous functions. In this particular case, we actually have that $Leb^{\delta}$ is physically equivalent to the additive family of H\"older continuous functions  $(J_t)_{t \ge 0}$. \qed 
\end{example}

\subsubsection{A classification of nonadditive equilibrium states}

In this section, we will apply Theorem \ref{BOU} to see how we can compare families with the same unique equilibrium measure just based on the information about the periodic data of the system. 

\begin{theorem}\label{CEM}
Let $\Lambda$ be a hyperbolic set for a topologically mixing $C^1$ flow~$\Phi$ and let $a = (a_t)_{t \ge 0}$ and $b = (b_t)_{t \ge 0}$ be two almost additive families of continuous functions with bounded variation. Then $a$ and $b$ have the same equilibrium measure if and only if there exists a constant $K > 0$ such that $|a_t(p)-b_t(p) -t(P_{\Phi}(a)-P_{\Phi}(b))| \le K$ for all $p \in \Lambda$ and $t \ge 0$ with $\phi_t(p) = p$.  
\end{theorem}

\begin{proof} 

First, let us suppose $|a_t(p)-b_t(p)-t(P_{\Phi}(a) - P_{\Phi}(b))| \le K$ for all $p \in \Lambda$ and $t \ge 0$ with $\phi_t(p) = p$. It follows from Theorem \ref{BOU} that 
\[
\sup_{t \ge 0}\sup_{x \in \Lambda}|a_t(x) - b_t(x) -t(P_{\Phi}(a) - P_{\Phi}(b))|< \infty.
\] 

Now consider the almost additive family $d= (d_t)_{t \ge 0}$ given by $d_t := b_t + t(P_{\Phi}(a) - P_{\Phi}(b))$. By the definition of nonadditive topological pressure and now the fact that the sequence $(a_t - d_t)/t$ converges uniformly for zero on $\Lambda$, we have respectively 
\[
P_{\Phi}(a) = P_{\Phi}(d) \quad \textrm{and} \quad \lim_{t \to \infty}\frac{1}{t}\int_{\Lambda}a_t d\nu  =  \lim_{t \to \infty}\frac{1}{t}\int_{\Lambda}d_t d\nu \quad \textrm{for all $\nu \in \cM(\Phi)$}.
\]
This readily implies that the families $a$ and $d$ share the same equilibrium measures. Moreover, since $P_{\Phi}(d) = P_{\Phi}(b) + (P_{\Phi}(a) - P_{\Phi}(b))$ and 
\[
\sup_{\mu \in \cM(\Phi)}\bigg(h_{\mu}(\Phi) + \lim_{t \to \infty}\frac{1}{t}\int_{\Lambda}d_td\mu\bigg) = \sup_{\mu \in \cM(\Phi)}\bigg(h_{\mu}(\Phi) + \lim_{t \to \infty}\frac{1}{t}\int_{\Lambda}b_td\mu\bigg) + P_{\Phi}(a) - P_{\Phi}(b),
\]
$d$ and $b$ also have the same equilibrium states. Hence, $a$ and $b$ do share the same equilibrium measures, as desired.

Now let's prove the converse. By Proposition \ref{H20}, both families $a$ and $b$ have unique equilibrium measures, each one of them satisfying the Gibbs property with respect to $\Phi$. Now, by assumption, suppose these equilibrium measures are the same unique equilibrium state $\eta \in \cM(\Phi)$. By the Gibbs property (Definition \ref{GBf}), for each sufficiently small $\epsilon > 0$ there exist constants $K_1(\epsilon) \ge 1$ and $K_2(\epsilon) \ge 1$ such that 
\[
K_1(\epsilon)^{-1} \le \frac{\eta(B_t(x,\epsilon))}{\exp[-tP_{\Phi}(a) + a_t(x)]}\le K_1(\epsilon) \quad \textrm{and} \quad K_2(\epsilon)^{-1} \le \frac{\eta(B_t(x,\epsilon))}{\exp[-tP_{\Phi}(b) + b_t(x)]}\le K_2(\epsilon)
\]
for all $x \in \Lambda$ and $t \ge 0$. This readily gives that 
\[
K_1(\epsilon)^{-1}K_2(\epsilon)^{-1} \le \exp[a_t(x)-b_t(x) - t(P_{\Phi}(a)-P_{\Phi}(b))]\le K_1(\epsilon)K_2(\epsilon)
\]
for all $x \in \Lambda$ and $t \ge 0$, which implies
\[
\sup_{x \in \Lambda}|a_t(x) - b_t(x) - t(P_{\Phi}(a)-P_{\Phi}(b))| \le \log(K_1(\epsilon)K_2(\epsilon)) \quad \textrm{for all $t \ge 0$}.
\]
In particular, we finally have that $|a_t(p)-b_t(p) -t(P_{\Phi}(a)-P_{\Phi}(b))| \le \log(K_1(\epsilon)K_2(\epsilon))$ for all $p \in \Lambda$ and $t \ge 0$ with $\phi_{t}(p) = p$.  
\end{proof}

Observe that Theorem \ref{BOU} together with Theorem \ref{CEM} shows that two almost additive families with bounded variation share the same unique equilibrium state if and only if they are cohomologous to each other (modulo some uniform constant) in the sense of definitions \ref{CHM} and \ref{SCH}. In this context, Theorem \ref{CEM} is a nonadditive counterpart of some classical results for hyperbolic flows (see for example Theorem 7.3.24 in \cite{FH20}).

\subsection{Regularity of almost and asymptotically additive families of potentials}

In this section, we will consider the regularity problem for nonadditive families, specially the asymptotically and almost additive ones based on Theorem \ref{main}. We will start investigating some natural simple examples in the non-hyperbolic context. After this, we will address the regularity issues for hyperbolic suspension flows and related setups.  

\subsubsection{Regularity in some non-hyperbolic scenarios}

Observe that one of the main ingredients in the proof of Theorem \ref{BOU} is the simultaneous existence of periodic and transitive points. A simple reasonable question is to ask what would happen in a system with no periodic points or no transitive data. In this regard, the most natural examples seem to be the linear flows with different types of directions on compact spaces. 

Let us start with an example of a setup where the periodic data is everywhere and with the same period.

\begin{example}\label{luis}
Let $\mathbb{T}^{n} = \R^{n}/\Z^{n}$ be the $n$-torus and consider $\alpha = (\alpha_1,...,\alpha_n) \in \R^{n}$ to be linear dependent, that is, there exist (not all zero) $k_j \in \Z$ such that $\sum_{j = 1}^{n}k_j\alpha_j = 0$. The linear flow $\Phi^{\alpha} = (\phi_t)_{t \in \R}$ on $\mathbb{T}^{n}$ in the direction $\alpha$ is defined by $\phi_t(x) = x + t\alpha \mod 1$. Letting $a := (a_t)_{t \ge 0}$ be any almost additive family of continuous functions with respect to $\Phi^{\alpha}$, Theorem \ref{main} together with Example \ref{TOR} guarantees the existence of a continuous function $b: \mathbb{T}^{n} \to \R$ such that $\lim_{t \to \infty}\|a_t - S_tb\|_{\infty}/t = 0$, where $\|.\|_{\infty}$ is the supremum norm on $\mathbb{T}^{n}$. 
Since $\Phi^{\alpha}$ is a periodic flow, by Lemma \ref{UBI} there exists a constant $L > 0$ (depending only on the period) such that $\sup_{t \ge 0}\|a_t - S_tb\|_{\infty} \le~L$. Observe that the uniform bound exists even if the family $a$ does not have the bounded variation property.  Moreover, it is clear that the additive family $(S_tb)_{t \ge 0}$ has bounded variation if and only if $a$ also has it. \qed
\end{example}

Let us now check what happens in the opposite extreme: transitive systems without periodic points.

\begin{example}\label{mama}
Let $\mathbb{T}^{n} = \R^{n}/\Z^{n}$ be the $n$-torus and consider $\alpha = (\alpha_1,...,\alpha_n) \in \R^{n}$ to be linear independent. In this case, the linear flow $\Phi^{\alpha} = (\phi_t)_{t \in \R}$ on $\mathbb{T}^{n}$ in the direction $\alpha$ given by $\phi_t(x) = x + t\alpha \mod 1$ is \emph{minimal}, that is, every orbit is dense in $\mathbb{T}^{n}$. Now let $a = (a_t)_{t \ge 0}$ be any almost (asymptotically) additive family of continuous functions. In particular, letting $\nu$ be the Lebesgue measure on $\mathbb{T}^{n}$ and $b: \mathbb{T}^{n} \to \R$ the continuous function given by Theorem \ref{main}, we have
\[
\lim_{t \to \infty}\frac{1}{t}\bigg\|a_t - t\int_{\mathbb{T}^{n}} b d\nu\bigg\|_{\infty} \le \lim_{t \to \infty}\frac{1}{t}\|a_t - S_tb\|_{\infty} + \lim_{t \to \infty}\frac{1}{t}\bigg\|S_tb - t\int_{\mathbb{T}^{n}}b d\nu\bigg\|_{\infty} = 0.
\]
In this case, the additive family of Theorem \ref{main} can be taken as the one generated by the constant function $\int_{\mathbb{T}^n} b d\nu$, which always satisfies the bounded variation property.  

On the other hand, the classical Gottshalk and Hedlund theorem for flows (see for example Theorem C in \cite{McC99}), guarantees that $\sup_{t \ge 0}\|S_tg - t\int_{\mathbb{T}^{n}}g d\nu\|_{\infty} = \infty$ for every continuous function $g: \mathbb{T}^{n} \to \R$ not $\Phi^{\alpha}$-cohomologous to a constant. Therefore, for these types of linear flows, Theorem \ref{BOU} fails even in the additive classical case  assuming functions with any strong form of regularity. \qed
\end{example}

\begin{remark}
Example \ref{luis} does not satisfy the hypotheses of Theorem \ref{BOU} but all the equivalences there are satisfied, even without asking for the bounded variation property on the families. On the contrary, Example \ref{mama} also does not satisfy the hypotheses of Theorem \ref{BOU} but the uniform bound property cannot be satisfied even asking for any type of regularity on the families or functions. 
\end{remark}

We will address now a more natural and richer setup for investigating H\"older potentials and families of functions with the bounded variation condition.  

\subsubsection{H\"older regularity of almost and asymptotically additive families}\label{BHR}

In this section, we are going to show how to construct almost and asymptotically additive families of (H\"older) continuous functions having the bounded variation condition with respect to some suspension flow but not physically equivalent to any additive family generated by a H\"older continuous function. Our general approach is based on the following examples, which are obtained from the ones showed in Theorem 11 in \cite{HS24}. 

\begin{theorem}\label{NHD}
Let $\sigma: \Sigma^{\Z} \to \Sigma^{\Z}$ be the two-sided full shift. Then
	
\begin{itemize}	
\item there exists an almost additive sequence of continuous functions with bounded variation and which is not physically equivalent to any additive sequence generated by a H\"older continuous function.
\item there exist almost additive sequences of H\"older continuous functions with bounded variation and which are not physically equivalent to any additive sequence generated by a H\"older continuous function.
\end{itemize}
\end{theorem}

\begin{proof}
	
Let $\sigma_1: \Sigma^{\N} \to \Sigma^{\N}$ be the left-sided full shift of finite type and $\sigma_2: \Sigma^{\Z} \to \Sigma^{\Z}$ be the two-sided full shift of finite type. Consider the canonical projection $\pi: \Sigma^{\Z} \to \Sigma^{\N}$, that is, $\omega = (...\omega_{-2}, \omega_{-1}, \omega_{0}, \omega_{1}, \omega_{2}, ... ) \longmapsto \pi(\omega) = (\omega_{1}, \omega_{2}, \omega_{3} ... )$. 

Let $\beta > 1$ and define $s:= s(\omega,\widetilde{\omega})$ on $\Sigma^{\N} \times \Sigma^{\N}$ as the smallest positive number $s$ such that $\omega_{s} \neq \widetilde{\omega_{s}}$. In this case, we consider the distance on $\Sigma^{\N}$ to be $d_1(\omega, \widetilde{\omega}) = \beta^{-s(\omega,\widetilde{\omega})}$ if $\omega \neq \widetilde{\omega}$ and $d_1(\omega, \widetilde{\omega}) = 0$ if $\omega = \widetilde{\omega}$. 
Following in an analogous way, we define $q:= q(\omega, \omega')$ on $\Sigma^{\Z} \times \Sigma^{\Z}$ as the smallest positive number $q$ such that $\omega_{-q} \neq \omega_{-q}'$ or $\omega_{q} \neq \omega_{q}'$. From this, we consider the distance on $\Sigma^{\Z}$ to be $d_2(\omega,\omega') = \beta^{-q(\omega,\omega')}$ if $\omega \neq \omega'$ and $d_2(\omega,\omega') = 0$ whenever $\omega = \omega'$. 

Now let $\cF = (f_n)_{n \in \N}$ be any almost additive sequence of continuous functions on $\Sigma^{\N}$ with respect to $\sigma_1$, satisfying the bounded variation condition and not physically equivalent to any additive sequence generated by a H\"older continuous function (for example, the sequence generated by the quasi-Bernoulli measure in Theorem 11 in \cite{HS24}). Let $\cG = (g_n)_{n \in}$ be the sequence on $\Sigma^{\N}$ given by $g_n = f_n \circ \pi$. Since $d_1(\pi(\omega),\pi(\omega')) \le d_2(\omega, \omega')$ for all $\omega, \omega' \in \Sigma^{\Z}$, $g_n$ is continuous for each $n \in \N$.  Since $(\sigma_1 \circ \pi)(\omega) = (\pi \circ \sigma_2)(\omega)$ for all $\omega \in \Sigma^{\Z}$, one can easily see that $\cG$ is almost additive. Moreover, since $\cF$ has bounded variation, we have
\[
\begin{split}
&\sup_{n \in \N}\{|g_n(\omega) - g_n(\widetilde{\omega})|: \omega, \widetilde{\omega} \in C_{i_1...i_n}\}\\
&\le \sup_{n \in \N}\{|f_n(\pi(\omega)) - f_n(\pi(\widetilde{\omega}))|: \pi(\omega), \pi(\widetilde{\omega}) \in C_{i_1...i_n}\cap\Sigma^{\N}\} < \infty. 
\end{split}
\]
That is, $\cG$ also has the bounded variation condition. 

Now suppose  $\phi: \Sigma^{\Z} \to \R$ is a H\"older continuous function such that $\cG$ is phisically equivalent to $(S_n\phi)_{n \in \N}$ with respect to $\sigma_2: \Sigma^{\Z} \to \Sigma^{\Z}$. Lemma 1.6 in \cite{Bow75a} (see also Section 3 in \cite{Sin72}) guarantees the existence of a H\"older continuous function $\psi: \Sigma^{\Z} \to \R$ cohomologous to $\phi$ and such that $\psi((...\omega_{-2}, \omega_{-1}, \omega_{0}, \omega_{1}, \omega_{2}, ... )) = \psi((\omega_{1}, \omega_{2}, ... ))$. That is, there exists a continuous functions $q: \Sigma^{\Z} \to \R$ satisfying $\phi - \psi = q \circ \sigma_{2} - q$ and $\psi$ is such that $\psi \circ \pi = \psi$. Based on this, for all $\omega \in \Sigma^{\Z}$ we get
\[
\begin{split}
\bigg|g_n(\omega) - \sum_{k=0}^{n-1}(\phi\circ \sigma^{k}_2)(\omega)\bigg| &= \bigg|f_n(\pi(\omega)) - \sum_{k=0}^{n-1}(\psi\circ \pi \circ \sigma^{k}_2)(\omega) + q - q\circ\sigma^{n}_{2}\bigg|\\
& \ge  \bigg|f_n(\pi(\omega)) - \sum_{k=0}^{n-1}(\psi\circ\sigma^{k}_1)(\pi(\omega))\bigg| - 2\|q\|_{\infty}.
\end{split}
\]
Since $\pi(\Sigma^{\Z}) = \Sigma^{\N}$ and $\cG$ is physically equivalent to $(S_n\phi)_{n \in \N}$, we have

\[
\lim_{n \to \infty}\frac{1}{n}\sup_{\omega'\in \Sigma^{\N}}\bigg|f_n(\omega') - \sum_{k=0}^{n-1}(\psi\circ\sigma^{k}_1)(\omega')\bigg| \le \lim_{n \to \infty}\frac{1}{n}\sup_{\omega \in \Sigma^{\Z}}\bigg|g_n(\omega) - \sum_{k=0}^{n-1}(\phi\circ \sigma^{k}_2)(\omega)\bigg| = 0,
\]
which contradicts Theorem 11 in \cite{HS24}. Therefore, the almost additive sequence $\cG:= \cF \circ \pi$ is not physically equivalent to any additive sequence generated by a H\"older continuous function. Moreover, fixing any real number $\gamma > 0$, by the density of H\"older functions on the space of continuous functions on $\Sigma^{\Z}$, one can see that for each $n \in \N$ there exists a H\"older function $h^{\gamma}_n$ such that $\|g_n - h^{\gamma}_n\|_{\infty} \le \gamma$. Notice that since $\cG$ has bounded variation, the sequence $\cH^{\gamma}:= (h^{\gamma}_n)_{n \in \N}$ also satisfies the bounded variation condition. Furthermore, Lemma 6 in \cite{HS24} guarantees that the sequence $\cH^{\gamma}$ also is almost additive. Hence, there exists an uncountable number of almost additive sequences of H\"older continuous functions on $\Sigma^{\Z}$, with bounded variation and which are not physically equivalent to any additive sequence generated by a H\"older function, as desired.
\end{proof}

Now let us start addressing how we can pass things from discrete to continuous time dynamical systems.  First, an auxiliary result.

\begin{lemma}\label{AAQ}
Let $M$ be a compact metric space. Every almost additive sequence of continuous functions $\cQ = (q_n)_{n \in \N}$ with respect to a continuous map $T: M \to M$ satisfy 
\begin{equation}\label{QFQ}
\sup_{n \in \N}\|q_n \circ T - q_n\|_{\infty} < \infty.
\end{equation}
\end{lemma}
\begin{proof}
Since $\cQ$ is almost additive, there exists a constant $K > 0$ such that
\[
-K + q_1(x)+ q_{n-1}(T(x)) \le q_{n}(x) \le q_{n-1}(T(x)) + q_1(x) + K
\]
for every $x \in M$ and every $n \in \N$. 
From this, we obtain that 
\begin{equation}\label{QN1}
\left|q_n(x) - q_{n-1}(T(x))\right| \le K + \|q_1\|_{\infty} := K_1 < \infty.
\end{equation}
In a similar way, we also have that 
\[
-K + q_1(T^{n-1}x)+ q_{n-1}(x) \le q_{n}(x) \le q_{n-1}(x) + q_1(T^{n-1}x) + K,
\]
and then
\begin{equation}\label{QN2}
\left|q_n(x) - q_{n-1}(x)\right| \le K_1
\end{equation}
for every $x \in M$ and every $n \in \N$.
	
It follows from \eqref{QN1} and \eqref{QN2} that 
\[
|q_n(T(x)) - q_n(x)| \le |q_{n}(T(x)) - q_{n+1}(x) | + |q_{n+1}(x) - q_n(x)|\le 2K_1
\]
for every $x \in M$ and every $n \in \N$.
	
Therefore, $\cQ$ satisfies property \eqref{QFQ} as desired. 
\end{proof}

\begin{proposition}\label{edu}
Let $\Phi = (\phi_t)_{t \in \R}$ be a suspension flow on $Y$ over a continuous invertible map $T: X \to X$ with continuous height function $\tau: X \to (0,\infty)$. Let $c = (c_n)_{n \in \N}$ be an almost additive sequence of continuous functions with respect to $T$ on $X$ and satisfying the bounded variation condition. Then, there exists an almost additive family of continuous functions $a = (a_t)_{t \ge 0}$ with respect to $\Phi$ on $Y$, satisfying the bounded variation condition and such that $a_n(x) = c_n(x)$ for all $x \in X$ and $n \in \N$. The same result holds for the asymptotically additive case.
\end{proposition}

\begin{proof}
Consider the floor function $\lfloor \cdot \rfloor : \R \to \Z$ given by $\lfloor x \rfloor = \max\{m \in \Z: m \le x\}$. Now for each $t > 0$ define the function $a_t: Y \to \R$ as 
\[
a_t(y) = a_t(\phi_s(x)) = a_{\lfloor t \rfloor}(\phi_s(x)) = c_{\lfloor t \rfloor}(x) \quad \textrm{and} \quad a_0 := c_0 \equiv 0.
\]

For simplicity, let us consider the height function to be constant $\tau \equiv 1$. Notice that by construction, $\phi_1 = T$ on $X$. The sequence $(a_n)_{n \in \N}$ is almost additive with respect to $\phi_1$ on $Y$. In fact, by definition and the almost additivity of $c$ on $X$, for all $y \in Y$ and $m,n \in \N$ we have
\[
\begin{split}
a_{m+n}(y) = a_{m+n}(\phi_s(x)) = c_{m+n}(x) &\le c_m(x) + c_n(T^{m}(x)) + C\\
& = a_m(\phi_s(x)) + a_n(\phi_1^{m}(x)) + C \\
& = a_m(y) + (a_n \circ \phi_s \circ \phi_1^{m})(x) + C\\
& = a_m(y) + a_n(\phi_1^{m}(\phi_s(x))) + C\\
& = a_m(y) + a_n(\phi_1^{m}(y)) + C
\end{split}
\]
Proceeding in the same manner, we also have $a_{m+n}(y) \ge a_m(y) + a_n(\phi_1^{m}(y)) - C$ for all $y \in Y$ and $m,n \in \N$, as desired. 

Let us now show that the family $a = (a_t)_{t \ge 0}$ is almost additive with respect to the flow~$\Phi$ on $Y$. For each $y \in Y$, $m \le t < m+1$ and $n \le s < n+1$, it follows from the definition of $(a_t)_{t > 0}$ that $a_{t+s}(y) = a_{m+n+\ell}(y)$, where $\ell \in \{0,1\}$. In this case, by the almost additivity of $(a_n)_{n \in \mathbb{N}}$ with respect to $\phi_1$ on $Y$, we get
\begin{equation}\label{AAD}
\begin{split}
a_{t+s}(y) = a_{m + n + \ell}(y) &\le a_{m+n}(y) + a_{\ell}(\phi_{m+n}(y)) + C  \\
& \le a_{m}(y) + a_n(\phi_m(y)) + a_{\ell}(\phi_{m+n}(y)) + 2C \\
& \le a_t(y) + a_s(\phi_m(y)) + \underbrace{\|a_{\ell}\|_{\infty} + 2C}_{=: C_1}\\
& = a_t(y) + a_s(\phi_t(y)) + \big[a_n(\phi_m(y)) - a_n(\phi_t(y))\big] + C_1.
\end{split}
\end{equation}
On the other hand, letting $y = \phi_w(x)$ for some $x \in X$ and $w \in [0,1)$ and $m = t + u$ with $u \in [0,1)$, we also have
\begin{equation}\label{DOT}
\begin{split}
\big|a_n(\phi_m)(y) - a_n(\phi_t(y))\big|&= \big|a_n(\phi_m(\phi_w(x))) - a_n(\phi_t(\phi_w(x)))\big|\\
& = \big|a_n(\phi_w(\phi_m(x))) - a_n(\phi_{u+w}(\phi_m(x)))\big| \\
&\le \big|c_n(\phi_m(x)) - c_n(\phi_1(\phi_m(x)))\big|.
\end{split}
\end{equation}
Since $c$ is almost additive with respect to $\phi_1$ on $X$, Lemma \ref{AAQ} guarantees the existence of a uniform constant $K > 0$ such that $\sup_{n \in \N}\sup_{x \in X}|c_n(\phi_m(x)) - c_n(\phi_1(\phi_m(x)))| \le K$. This readily implies that $|a_n(\phi_m)(y) - a_n(\phi_t(y))| \le K$ for all $m, n \in \N$, $t > 0$ and $y \in Y$. Hence, it follows now from \eqref{AAD} that 
$a_{t+s}(y) \le a_t(y) + a_s(\phi_t(y)) + K + C$. The other inequality can be obtained in a similar way, and we conclude that $a$ is almost additive with respect to $\Phi$ on $Y$. 

Considering the Bowen-Walters distance $d_Y$ on $Y$ (see \cite{BW72} or Appendix A in \cite{BS00}) and using the continuity of $c_n: X \to \R$ for all $n \in \N$, one can easily check that the function $a_t: Y \to \R$ is continuous for each $t \ge 0$. 

Let us conclude the proof showing that the family $a$ also has the bounded variation property. Take two arbitrary points $y, z \in Y$ such that $d_Y(\phi_{\tau}(y),\phi_{\tau}(z)) < \epsilon$ for $\tau \in [0,t]$ with $m \le t < m+1$. Writing $y = \phi_{u}(x)$ and $z = \phi_{r}(x')$ for $x, x' \in X$ and $u, r \in [0,1)$, we have
\begin{equation}\label{BVP}
|a_t(y) - a_t(z)| = |a_t(\phi_{u}(x)) - a_t(\phi_{r}(x'))| = |c_m(x) - c_m(x')|
\end{equation}
In particular, we have 
\[
\begin{split}
& d_X(x, x') \le d_Y(y,z) < \epsilon\\
& d_X(\phi_1(x),\phi_1(x')) \le d_Y(\phi_u(\phi_1(x)),\phi_r(\phi_1(x'))) = d_Y(\phi_1(y),\phi_1(z)) < \epsilon\\
&\cdots \cdots \cdots  \\
& d_X(\phi_{m-1}(x),\phi_{m-1}(x')) \le d_Y(\phi_u(\phi_{m-1}(x)),\phi_r(\phi_{m-1}(x'))) = d_Y(\phi_{m-1}(y),\phi_{m-1}(z)) < \epsilon,
\end{split}
\]
where $d_Y$ is the Bowen-Walters distance on $Y$ and $d_X$ is any given distance on the base space $X$. 
Since the sequence $c$ has bounded variation (with respect to $T = \phi_1$ on $X$), there exists a constant $L = L(\epsilon) > 0$ such that $|c_m(x) - c_m(x')| \le L$. We conclude now from \eqref{BVP} that the family $a$ also has the bounded variation property with the same parameters $\epsilon > 0$ and $L > 0$ as the sequence $c$ on the base space $X$. 

For the general case where $r$ is any positive continuous function, we have $T(x) = \phi_{\tau(x)}(x)$ and $T^{m}(x) = \phi_{\tau_m(x)}(x)$ for all $m \in \mathbb{N}$ and $x \in X$, with $\tau_m = \sum_{k=0}^{m-1}\tau \circ T^{k}$. In this case, we define $a_t\colon Y \to \mathbb{R}$ as 
\begin{equation}\label{BVQ}
a_t(y) = a_t(x,s) := a_{\tau_n(x)}(\phi_s(x)) := c_n(x) \quad \textrm{with} \quad a_0 = c_0 := 0,
\end{equation}
when $t,s \ge 0$ and $\tau_n(x) \le t + s \le \tau_{n+1}(x)$. Making the necessary modifications and proceeding as in the case with $\tau \equiv 1$, one can easily check that the family $a$ is almost additive with respect to $\Phi$ on $Y$. The continuity of $c_m: X \to \R$ for all $m \ge 1$ together with definition \eqref{BVQ} directly implies that $a_t: Y \to \R$ is continuous for each $t \ge 0$.  Moreover, since the sequence $c$ has bounded variation with respect to the map $T$ on $X$, the same relation between the distance on $X$ and the Bowen-Walters distance on $Y$ guarantees the bounded variation property for the family $a$ with respect to $\Phi$ on $Y$.

Now suppose $c = (c_n)_{n \in \N}$ is asymptotically additive with respect to $T: X \to X$ and consider again the family $a = (a_t)_{t \ge 0}$ defined in \eqref{BVQ}. By the asymptotically additivity of $c$, given any $\epsilon > 0$ there exists a continuous function $h_{\epsilon}: X \to \R$ such that
\begin{equation}\label{asy}
\limsup_{n \to \infty}\frac{1}{n}\sup_{x \in X}\bigg|c_n(x) - \sum_{k = 0}^{n-1}(h_{\epsilon}\circ T)(x)\bigg| < \epsilon.
\end{equation}
On the other hand, by Lemma \ref{LE2}, there exists a continuous function $g_{\epsilon}: Y \to \R$ such that $I_{g_{\epsilon}}|_X = h_{\epsilon}$. By the definition of $a$, for each $y = \phi_u(x)$ with $u \in [0, \sup\tau)$, and $\tau_{n}(x) \le t < \tau_{n+1}(x)$ we have
\[
\begin{split}
\bigg|a_t(y) - \int_{0}^{t}(g_{\epsilon}\circ \phi_s)(y)ds\bigg| &= \bigg|a_t(\phi_u(x)) - \int_{u}^{t+u}(g_{\epsilon}\circ \phi_s)(x)ds\bigg|\\
& \le \bigg|c_n(x) - \sum_{k=0}^{n-1}(h_{\epsilon} \circ T^{k})(x)\bigg| + \sup \tau \sup g_{\epsilon} + \sup h_{\epsilon}. 
\end{split}
\]

Since $n \to \infty$ implies $t \to \infty$, we conclude from \eqref{asy} that $a$ is asymptotically additive with respect to $\Phi$ on $Y$. The continuity and bounded variation condition of $a$ clearly follow from the same arguments presented in the almost additive case.   
\end{proof}

Suspension flows over two-sided subshifts of finite type with H\"older continuous height functions are also called \emph{hyperbolic symbolic flows} (see \cite{FH20}). One can check that additive families generated by H\"older continuous functions satisfy the bounded variation condition with respect to hyperbolic symbolic flows, and Proposition \ref{H20} clearly also can be extended to these types of flows in general. On the other hand, just as in the discrete-time case, one can find asymptotically additive families having bounded variation with respect to an hyperbolic symbolic flow, but admitting more than one equilibrium measure. 

The following result is a continuous-time counterpart of Theorem \ref{NHD}:

\begin{theorem}\label{SFT}
Let $\Phi = (\phi_t)_{t \in \R}$ be a suspension flow over $\sigma:\Sigma^{\Z} \to~\Sigma^{\Z}$ and with a H\"older continuous height function $\tau: \Sigma^{\Z} \to (0, \infty)$. Then:

\begin{enumerate}
	
\item there exist almost additive families of H\"older continuous functions with respect to $\Phi$, satisfying the bounded variation condition and not physically equivalent to any additive family generated by a H\"older continuous function; 

\item there exist asymptotically additive families of H\"older continuous functions with respect to $\Phi$, satisfying the bounded variation condition, admitting a unique equilibrium measure but not physically equivalent to any additive family generated by a H\"older continuous function.

\end{enumerate}

\end{theorem}

\begin{proof}
First, let $c = (c_n)_{n \in \N}$ be any sequence from Theorem \ref{NHD}, that is, $c$ is an almost additive sequence of continuous functions with respect to $\sigma: \Sigma^{\Z} \to \Sigma^{\Z}$, satisfying the bounded variation condition and not physically equivalent to any additive sequence generated by a H\"older continuous function. By Proposition \ref{edu} there exists an almost additive family of continuous functions $a = (a_t)_{t \ge 0}$ with respect to $\Phi$ on $Y$, with bounded variation and such that $a_n(x) = c_n(x)$ for all $x \in \Sigma^{\Z}$ and $n \in \N$. Suppose that 
\[
\lim_{t \to \infty}\frac{1}{t}\sup_{y \in Y}\bigg|a_t(y) - \int_{0}^{t}(b \circ \phi_s)(y)ds\bigg| = 0 \quad \textrm{where} \quad \textrm{$b: Y \to \R$ is H\"older}.  
\]
In particular, this implies that 
\begin{equation}\label{pou}
\lim_{t \to \infty}\frac{1}{t}\sup_{x \in \Sigma^{\Z}}\bigg|a_t(x) - \int_{0}^{t}(b \circ \phi_s)(x)ds\bigg| = 0.
\end{equation}

By the proof of Lemma 15 in \cite{BH21b}, for each $t > 0$ there exists a unique $n \in \N$ with $t = \tau_n(x) + \kappa$ for some $\kappa \in [0, \sup\tau]$ such that 
\[
\bigg|\int_{0}^{t}(b\circ \phi_s)(x)ds - \sum_{k=0}^{n-1}(I_{b}\circ \sigma^{k})(x)\bigg| \le \sup b \sup\tau,
\]
where $I_b(x) = \int_{0}^{\tau(x)}(b \circ \phi_s)(x)ds$. 

Thus, it follows from \eqref{pou} that 
\[
\lim_{n \to \infty}\frac{1}{n}\sup_{x \in \Sigma^{\Z}}\bigg|c_n(x) - \sum_{k=0}^{n-1}(I_{b}\circ \sigma^{k})(x) \bigg| = \lim_{n \to \infty}\frac{1}{n}\sup_{x \in \Sigma^{\Z}}\bigg|a_n(x) - \int_{0}^{n}(b \circ \phi_s)(x)ds\bigg| = 0.
\]

Since $b: Y \to \R$ is H\"older, by Proposition 18 in \cite{BS00} we also have that $I_b: X \to \R$ is H\"older. Hence, the sequence $c$ is physically equivalent to the additive sequence generated by $I_b$, which is a contradiction. 

Now fix a number $\gamma > 0$. By the density of H\"older functions on the space of continuous functions, for each $t \ge 0$ there exists a H\"older continuous function $b^{\gamma}_t: Y \to \R$ such that $\sup_{y \in Y}|b^{\gamma}_t(y) - a_t(y)| \le \gamma$. It is easy to check that the family $b^{\gamma}:= (b^{\gamma}_t)_{t \ge 0}$ is almost additive and satisfy the bounded variation condition with respect to the flow $\Phi$ on $Y$. Moreover, since $b^{\gamma}$ is physically equivalent to $a$, it cannot be physically equivalent to any additive family generated by a H\"older continuous function, as desired. 

Now let's prove the second item. It was showed in \cite{HS24} the existence of asymptotically additive sequences of continuous functions $c = (c_n)_{n \in \N}$ with respect to the left-sided full shift $\sigma: \Sigma^{\N} \to \Sigma^{\N}$ satisfying the bounded variation condition, with a unique equilibrium measure but not physically equivalent to any additive sequence generated by a H\"older continuous function. Proceeding as in the proof of Theorem \ref{NHD} and with some abuse of notation, one can also assume that the sequence $c$ is asymptotically additive with respect to the two-sided full shift $\sigma: \Sigma^{\Z} \to \Sigma^{\Z}$. By Proposition \ref{edu}, following as in the proof of item 1 and also using density of H\"older functions, one can guarantee the existence of a family of asymptotically additive H\"older continuous functions $a = (a_t)_{t \ge 0}$ with respect to $\Phi$ on $Y$ and satisfying the bounded variation condition but never physically equivalent to any additive family generated by a H\"older continuous function. Moreover, it is clear from \eqref{eq15} (adapted for suspension flows) that $a$ also admits a unique equilibrium measure, which is induced by the unique equilibrium measure for $c$ (see also the proof of Theorem 3.5 in \cite{BH21a}). 

\end{proof}

\begin{remark}
By Theorem 11 in \cite{HS24}, one can obtain the same class of examples in Theorem~\ref{SFT} to the case of suspension semi-flows over left-sided full shifts $\sigma: \Sigma^{\N} \to \Sigma^{\N}$. 
\end{remark}

Notice that the counter-examples in Theorem \ref{SFT} show that the physical equivalence result of Theorem \ref{main} does not always allow us to reduce the study of almost additive families with bounded variation to the case of classical potentials with H\"older regularity. Since the thermodynamic formalism and multifractal phenomena for H\"older continuous potentials are suitable and well understood for hyperbolic setups in general, Theorem \ref{SFT} indicates, as in the case of maps, a crucial barrier considering the exchange of information between the additive, almost additive and asymptotically additive worlds with respect to continuous-time dynamical systems as well. 

\subsubsection{Bowen regularity of almost and asymptotically additive families}\label{BHB}

In this section, based on the physical equivalence relations obtained in the first part of this work, we briefly address more general regularity aspects of almost and asymptotically additive families of continuous functions. In general, equilibrium measures satisfying the Gibbs property play a relevant role in our approach. 

Recall that a function $\xi: M \to \R$ is Bowen (with respect to a flow $\Phi$ on a topological space $M$) if there exist $\kappa > 0$ and $\epsilon > 0$ such that for $x, y \in M$ and $t \ge 0$, we have that $d(\phi_s(x),\phi_s(y)) < \epsilon$ for every $s \in [0,t]$ implies $|S_t\xi(x) - S_t\xi(y)| \le \kappa$. This is also the same as saying that the additive family $(S_t\xi)_{t \ge 0}$ has the bounded variation property with respect to $\Phi$ on $M$.

\begin{theorem}\label{BKMH}
Let $\Phi = (\phi_t)_{t \in \R}$ be a suspension flow over $\sigma: \Sigma^{\Z} \to~\Sigma^{\Z}$ and with a H\"older continuous height function $\tau: \Sigma^{\Z} \to (0, \infty)$ and let $a = (a_t)_{t \ge 0}$ be an almost additive family of continuous functions with respect to $\Phi$ on $Y$ with bounded variation. Then, the following statements are equivalent:
\begin{enumerate}
\item the equilibrium measure for the family $a$ is Gibbs with respect to some continuous Bowen function;
		
\item there exists a continuous Bowen function $b: Y \to \R$ such that
\[
\lim_{t \to \infty}\frac{1}{t}\sup_{y \in Y}|a_t(y) - S_tb(y)| = 0.
\]
		
\item there exists a continuous Bowen function $b: Y \to \R$ such that
\[
\sup_{t \ge 0}\sup_{y \in Y}|a_t(y) - S_tb(y)| < \infty.
\]
\end{enumerate} 
\end{theorem}

\begin{proof}
Let us start proving that 1 implies 3. By the appropriate version of Proposition \ref{H20} for hyperbolic symbolic flows, the family $a$ has a unique equilibrium measure $\nu$, which is Gibbs with respect to $a$. By hypothesis, $\nu$ is also Gibbs with respect to some continuous Bowen function $b: Y \to \R$. Then, for some sufficiently small $\delta > 0$, there exist constants $K_1 = K_1(\delta) \ge 1$ and $K_2 = K_2(\delta)\ge 1$ such that 
\begin{equation}\label{GBB}
K_1^{-1} \le \frac{\nu(B_t(y,\delta))}{\exp[-tP_{\Phi}(a) + a_t(y)]}\le K_1 \quad \textrm{and} \quad K_2^{-1} \le \frac{\nu(B_t(y,\delta))}{\exp[-tP^{\rm{classic}}_{\Phi}(b) + S_tb(y)]}\le K_2
\end{equation}
for all $y \in Y$ and $t \ge 0$, where $P^{\rm{classic}}_{\Phi}(b)$ is the classical (additive) topological pressure for $b$ with respect to $\Phi$. This clearly implies that $|a_t(y) - S_t\widetilde{b}(y)| \le \log{K_1K_2}$ for all $y \in Y$ and $t \ge 0$, where $\widetilde{b} := b + P^{\rm{classic}}_{\phi}(b) - P_{\Phi}(a)$. Since $\widetilde{b}$ is also a continuous Bowen function, we just obtained item 3. 
	
Now suppose 3 holds, that is, there exist a uniform constant $K_3 > 0$ and a continuous Bowen function $b: Y \to \R$ such that $|a_t(y) - S_tb(y)| \le K_3$ for all $y \in Y$ and $t \ge 0$. Then, by the Gibbs property with respect to the family $a$ in \eqref{GBB} and the fact that $P_{\Phi}(a) = P^{\rm{classic}}_{\Phi}(b)$ in this case, for all sufficiently small $\delta > 0$ we get that
\[
(K_1e^{K_3})^{-1}= K_1^{-1}e^{-K_3} \le \frac{\nu(B_t(y,\delta))}{\exp[-tP^{\rm{classic}}_{\Phi}(b) + S_tb(y)]}\le K_1e^{K_3}
\]
for all $y \in Y$ and $t \ge 0$, which is item 1, as desired. Since every hyperbolic symbolic flow is topologically transitive (Proposition 1.6.30 in \cite{FH20}) and satisfy the hypothesis of the Closing Lemma (see for example Corollary 18.1.8 in \cite{KH12}), Corollary \ref{EQ2} immediately gives that items 2 and 3 are equivalent, and the result is proved. 
\end{proof}

Recall that in the case of hyperbolic symbolic flows or locally maximal hyperbolic sets for $C^1$ topologically mixing flows, an almost additive family satisfies the bounded variation condition if and only if it admits a Gibbs measure (see Definition \ref{GBf} and Proposition 15 in \cite{HS24}). In Theorem \ref{BKMH}, the equivalence between items 1 and 2 indicates a possible way of classifying almost additive families with bounded variation with respect to hyperbolic symbolic flows or locally maximal hyperbolic sets for $C^1$ topologically mixing flows. Based on Theorem \ref{BKMH} and proceeding as in \cite{HS24} (and using the same nomenclature), we also can propose the following classification of families with respect to hyperbolic symbolic flows:

\begin{itemize}
\item \textbf{Type 1:} Almost additive families with bounded variation and admitting Gibbs measures with respect to some Bowen continuous function;

\item \textbf{Type 2:} Almost additive families with bounded variation and but not admitting Gibbs measures with respect to any Bowen continuous function;

\item \textbf{Type 3:} Almost additive families without the bounded variation condition but having a unique equilibrium measure;

\item \textbf{Type 4:} Almost additive families having more than one equilibrium measure. 
\end{itemize}

\begin{remark}
As discussed in \cite{HS24}, we can easily construct examples of types 1, 3 and 4. On the other hand, examples of families of type 2 seem to be much more complicated to produce or they actually don't exist. In the discrete-time setup, the existence of sequences of type 2 is connected with the problem of relating Gibbs and quasi-Bernoulli measures and, as far as we know, is still an open question even in the case of full shifts of finite type.
\end{remark}

Let us now finish our work showing how to treat the Bowen regularity problem for asymptotically additive families.

\textbf{Asymptotically additive families.} Let $\cG = (g_n)_{n \in \N}$ be an asymptotically additive sequence of continuous functions with bounded variation with respect to $\sigma: \Sigma^{\Z} \to \Sigma^{\Z}$ and with a unique equilibrium measure but never physically equivalent to any additive sequence generated by a Bowen function (such an example was given in \cite{HS24}). Now, as in Proposition \ref{edu}, consider an asymptotically additive family $a = (a_t)_{t \ge 0}$ with respect to the hyperbolic symbolic flow $\Phi$ on $Y$ (with height function $\tau$) and such that $a_n(x) = g_n(x)$ for all $x \in \Sigma^{\Z}$ and all $n \in \N$. 

Suppose the existence of a continuous Bowen function $b: Y \to \R$ such that $a$ is physically equivalent to the additive family $(S_tb)_{t \ge 0}$. By the appropriate versions of the Lemmas 3.1 and 3.3 in \cite{BH21a} for hyperbolic symbolic flows, the sequence $c = (c_n)_{n \in \N}$ given by $c_n(x) = \int_{0}^{\tau_n(x)}(b \circ \phi_s)(x)ds$ is additive and satisfy the bounded variation condition with respect to $\sigma: \Sigma^{\Z} \to \Sigma^{\Z}$. By the physical equivalence relation between $a$ and $(S_tb)_{t \ge 0}$, we have in particular that
\begin{equation}\label{oda}
\lim_{n \to \infty}\frac{1}{n}\sup_{x \in \Sigma^{\Z}}\bigg|g_n(x) - c_n(x)\bigg| = \lim_{n \to \infty}\frac{1}{n}\sup_{x \in \Sigma^{\Z}}\bigg|a_n(x) - c_n(x)\bigg| \le \lim_{t \to \infty}\frac{1}{t}\sup_{y \in Y}|a_t(y) - S_tb(y)| = 0.
\end{equation}

Since, by the proof of Lemma 15 in \cite{BH21b}, $c_n(x) = \sum_{k = 0}^{n-1}I_b \circ \sigma^{k}(x)$ for all $x \in \Sigma^{\Z}$ and all $n \in \N$, the sequence $(S_nI_b)_{n \in \N}$ has the bounded variation condition. Hence, it follows from \eqref{oda} that the sequence $\cG$ is physically equivalent to the additive sequence $(S_nI_b)_{n \in \N}$ generated by the continuous Bowen function $I_b : \Sigma^{\Z} \to \R$, which is a contradiction. Hence, by construction, the asymptotically additive family $a$ satisfies the bounded variation condition and has a unique equilibrium measure (with respect to $\Phi$ on $Y$), but cannot be physically equivalent to any additive family generated by a Bowen continuous function. 

Therefore, as it also happened for the asymptotically additive case with respect to discrete-time dynamical systems, the physical equivalence result obtained in Theorem \ref{main} does not always preserve the Bowen regularity property. \qed

\subsection{Concluding remarks}

Observe that all the results in the regularity sections are developed for hyperbolic symbolic flows. There is a deeper reason for that, which comes all the way from \cite{BKM20}. In this last work, studying almost additivity in the context of planar matrix cocycles, the authors showed an example of a quasi-Bernoulli measure that is not Gibbs for any H\"older continuous function with respect to the left-sided full shift map. In view of the nonadditive versions of the Liv\v{s}ic theorem for maps and flows (Theorem 5 in \cite{HS24} and Theorem \ref{BOU}, respectively), this particular example plays a fundamental role in the production of the counter-examples in \cite{HS24} for the left-sided full shift map and, consequently, the ones in Theorem \ref{SFT} for symbolic flows. Based on this, morally speaking, all the counter-examples and results discussed here in the regularity section can be adapted to the case of hyperbolic flows and, more generally, to suspension flows over topologically mixing subshifts of finite type. Finally, let us mention the still open problem of the existence of sequences and families of type 2, which is related to the understanding of how some results in \cite{BKM20} could accommodate Bowen continuous functions, going beyond the H\"older regularity case previously considered in it. A positive answer in this direction would finally reveal the existence of quasi-Bernoulli measures that do not satisfy the Gibbs property for any continuous function, consequently giving examples of sequences of type 2 with respect to the full shift of finite type. By our constructions in this note, we could as well produce examples of families of type 2 with respect to hyperbolic symbolic flows and hyperbolic flows (via Markov partitions).  \\

\textbf{Acknowledgments:} The author was partially supported by NSF of China, grant no.~12222110 and is supported by the National Council for Scientific and Technological Development, CNPq, process 446515/2024-8.

%\textbf{Data availability statement:} Data sharing is not applicable to this article as no data sets were generated or analyzed during the current study. \\

%All authors declare that they have no conflicts of interest and no relevant financial or non-financial interests to disclose. 

\bibliographystyle{alpha}

\end{document}